\documentclass[11pt,a4paper]{amsart}
\usepackage[margin=27mm]{geometry}
\usepackage{amsmath,amssymb,amsthm,mathtools,mathrsfs}
\usepackage[T1]{fontenc}
\usepackage{lmodern,microtype,booktabs,array,enumitem,xcolor}
\usepackage[hidelinks]{hyperref}
\usepackage{xurl}
\hypersetup{pdftitle={Spectral summation and critical asymptotics of Newton knot invariants},pdfauthor={Honghuai Fang},pdfsubject={Newton interpolation, spectral summation, critical asymptotics and logarithmic monodromy}}
\setlist[enumerate]{itemsep=3pt,topsep=4pt}
\allowdisplaybreaks[2]
\numberwithin{equation}{section}
\newtheorem{theorem}{Theorem}[section]
\newtheorem{proposition}[theorem]{Proposition}
\newtheorem{lemma}[theorem]{Lemma}
\newtheorem{corollary}[theorem]{Corollary}
\newtheorem*{theoremA}{Theorem A}
\newtheorem*{theoremB}{Theorem B}
\newtheorem*{theoremC}{Theorem C}
\newtheorem*{theoremD}{Theorem D}
\newtheorem*{theoremE}{Theorem E}
\theoremstyle{definition}
\newtheorem{definition}[theorem]{Definition}

\theoremstyle{remark}
\newtheorem{remark}[theorem]{Remark}
\newcommand{\C}{\mathbb C}\newcommand{\Z}{\mathbb Z}
\newcommand{\cB}{\mathcal B}
\newcommand{\cL}{\mathcal L}\newcommand{\cS}{\mathcal S}

\newcommand{\cO}{\mathcal O}
\newcommand{\Res}{\operatorname{Res}}

\newcommand{\T}{\mathsf T}
\newcommand{\sRes}{\mathfrak R}
\newcommand{\W}{\mathscr W}
\newcommand{\ord}{\operatorname{ord}}

\newcommand{\qb}[2]{\genfrac{[}{]}{0pt}{}{#1}{#2}_Q}
\newcommand{\Gauss}[2]{\genfrac{[}{]}{0pt}{}{#1}{#2}_{Q}}
\title[Spectral summation and critical asymptotics]{Spectral summation and critical asymptotics of Newton knot invariants}
\author{Honghuai Fang}
\address{Institute for Theoretical Sciences, Westlake University, 600 Dunyu Road, Xihu District, Hangzhou, Zhejiang 310030, China}
\email{fanghonghuai@westlake.edu.cn}
\subjclass[2020]{Primary 33D15, 39A13, 57K14; Secondary 30D15, 41A60}
\keywords{colored HOMFLY--PT polynomial, Newton interpolation, spectral summation, Stokes data, critical asymptotics, growing rank, geometric periods}
\date{}

\begin{document}

\begin{abstract}
The symmetric colored invariants of a knot are defined at integer colors,
and Newton interpolation assembles them into a formal function of a
continuous color variable. We sum it by
$q$-Borel--Laplace summation and use the resulting analytic function to
study the invariants as $q\to1$. For every knot and every rank, the critical summation level is a rational number determined by the growth of the colored invariants; when it is positive, the critical Borel transform is meromorphic in the
whole plane. When the Newton coefficients have analytic tails, the
Borel transform is assembled from convergent sums of simple poles along
geometric progressions, and for such sums the Stokes jump of the Laplace
sum determines every pole and residue. For an infinite family of twist
knots containing the figure-eight knot, a positive Markov kernel gives
an exact formula for the remainder of the Newton series. With it we
show that, near a root of the Alexander polynomial, the rescaled
analytic interpolation converges to Gaussian and Scorer-type profiles,
uniformly for complex values of the rescaled color, and we determine
how fast the rank may grow before these limits change. Finally,
continuation in $\log q$ produces a logarithmic monodromy whose
coefficient is a convergent pairing of the same Borel residues that
give the Stokes jump. For the figure-eight knot, the asymptotics of this
coefficient are governed by the hyperbolic volume at every rank.
\end{abstract}
\maketitle
\setcounter{tocdepth}{1}
\tableofcontents

\section{Introduction}\label{merged:introduction}

The colored quantum invariants of a knot are indexed by an integer
color, but many of their most interesting features appear in limits in
which the color and the quantum parameter vary together: the
Melvin--Morton--Rozansky expansion near the trivial parameter
\cite{GLAsymptotics,GLAnalytic}, and the exponential growth governed by
hyperbolic geometry, with the resurgent structure behind it
\cite{Witten,GGM,GGMW}. In such limits one would like to treat the color
as a continuous variable, and for this one needs an analytic function of
the color that takes the prescribed values at the integer colors. Such a
function is far from unique, since any entire function vanishing at the
integer colors may be added to it. Newton interpolation produces a
natural candidate, but only as a formal series, which in general
diverges. This paper is about summing that series and about what the
resulting function tells us about the invariants.

The main results are as follows. For an arbitrary knot at a fixed rank,
we determine the critical level of $q$-Borel summation and show that the
critical Borel transform is meromorphic in the whole plane (Theorem~A).
For Newton coefficients with analytic tails, we describe the Borel
singularities, the Laplace sums and their Stokes jumps through
convergent resolvent expansions, and we show that the jump determines
the poles and residues of the Borel transform (Theorem~B). For an
infinite family of twist knots containing the figure-eight knot, a
positive kernel gives an exact formula for the remainder of the Newton
series. With it we prove Gaussian and Scorer-type limits at the roots of
the Alexander polynomial (Theorem~C) and determine how these limits
change when the rank grows (Theorem~D). Finally, we continue the sums in
$\log q$ and show that the coefficient of the resulting logarithmic
monodromy is a convergent pairing of Borel residues; for the
figure-eight knot its asymptotics are governed by the hyperbolic volume
(Theorem~E).

\subsection{Newton interpolation and its summation}

Let $J_r^{(n)}(K;q)$ be the reduced, zero-framed quantum invariant of a
knot $K\subset S^3$ colored by $\operatorname{Sym}^r(\mathbb C^n)$, that
is, the symmetric colored HOMFLY--PT polynomial specialized at $A=q^n$,
and put $Q=q^2$. We interpolate at the nodes
\[
 x_r^{(n)}=q^{n+2r}+q^{-n-2r},\qquad r=0,1,2,\ldots,
\]
in the Newton basis $B_k^{(n)}(X)=\prod_{j=0}^{k-1}(X-x_j^{(n)})$. The
Newton coefficients $c_{K,n,k}$ are determined by the triangular system
\[
 J_r^{(n)}(K;q)=\sum_{k=0}^r c_{K,n,k}B_k^{(n)}(x_r^{(n)}),
 \qquad r\ge0,
\]
and the formal series
\[
 F_{K,n}(X)=\sum_{k\ge0}c_{K,n,k}B_k^{(n)}(X)
\]
takes the value $J_r^{(n)}(K;q)$ at $X=x_r^{(n)}$, because all terms
with $k>r$ vanish there. By the cyclotomic expansions of
\cite{FangZhou}, the coefficients are Laurent polynomials in $q$, and
$F_{K,n}$ is an element of the Newton completion of \cite{Fang}. We call
$X$ the color variable and the nodes the quantum colors.

At a quantum color $F_{K,n}$ is a finite sum; at any other color it is
an infinite series whose terms usually grow very fast. Indeed, the basis
factors into $q$-Pochhammer symbols, and for every color $X$ that is not
a node
\begin{equation}\label{intro:basis-growth}
 B_k^{(n)}(X)\sim(-1)^kq^{-k^2-(n-1)k}E_n(X),\qquad
 E_n(X)=\prod_{j\ge0}\bigl(1-q^{n+2j}X+q^{2n+4j}\bigr)
\end{equation}
as $k\to\infty$ (Theorem~\ref{thm:bridge}). The entire function $E_n$
has simple zeros exactly at the nodes. The series can therefore
converge at a color that is not a node only if
$c_{K,n,k}q^{-k^2-(n-1)k}\to0$; for the twist knots studied below its
terms are of size $q^{-k^2}$ up to exponential factors.

To sum it, we separate the color from the size of the terms. With an
auxiliary variable $z$, put
\[
 \widehat F_{K,n}(X,z)=\sum_{k\ge0}c_{K,n,k}B_k^{(n)}(X)z^k .
\]
For $0<q<1$ the coefficients of this series are bounded by the
exponential of a quadratic function of $k$
(Corollary~\ref{cor:universal-finite-level}), which is the natural
setting for $q$-Borel--Laplace summation \cite{DVZ,RSZ,Dreyfus}. We use
the moments and the kernel of Kashaev, Mangazeev and
Mari\~no~\cite{KMM}: for $0<p<1$,
\[
 M_k(p)=(p;p)_kp^{-k(k+1)/2},\qquad
 \cB_p\Big(\sum_{k}a_kz^k\Big)(\xi)=\sum_{k}\frac{a_k}{M_k(p)}\,\xi^k,
\]
\[
 \cL_pg(z)=\frac1{-\log p}\int_0^\infty\frac{g(zt)}{(-t;p)_\infty}\,dt .
\]
Since
\[
 \frac1{-\log p}\int_0^\infty\frac{t^k\,dt}{(-t;p)_\infty}=M_k(p)
 \qquad(k\ge0)
\]
(Theorem~\ref{thm:kernel}), $\cL_p\cB_p$ is the identity on polynomials
in $z$. At a node $\widehat F_{K,n}$ is a polynomial in $z$, and its sum
at $z=1$ is the colored invariant. Consequently, whenever the sum is
defined at $z=1$, its value there is a function of the color that
interpolates the colored invariants; in the cases treated below it is
entire. Two interpolations that are entire in the color differ by an
entire multiple of $E_n$ (Theorem~\ref{thm:E}).

Three questions now arise. At which base $p$ should one sum, and what
does the Borel transform look like there? The Borel transform has
poles, and the Laplace sum jumps across the rays through them; what do
these Stokes jumps record? And how does the analytic interpolation
behave as $q\to1$, in particular where the classical limit is singular,
and after continuation to complex values of $q$? We answer the first
question for every knot, the second for Newton coefficients with
analytic tails, and the third for the twist knots of
Section~\ref{subsec:intro-twist}.

\subsection{Critical levels for arbitrary knots}

Fix $0<q<1$ and $n\ge2$, and write $h=-\log Q$. The quadratic growth of
the invariants and of the Newton series is measured by
\begin{align*}
 \delta_{K,n}
 &=\max\left\{0,\limsup_{r\to\infty}
                   \frac{2\log|J_r^{(n)}(K;q)|}{hr^2}\right\},\\
 \beta_{K,n}
 &=\max\left\{0,\limsup_{k\to\infty}
          \frac{2\log|c_{K,n,k}B_k^{(n)}(X)|}{hk^2}\right\},
\end{align*}
where $X$ is any color that is not a node and $\log0=-\infty$; by
\eqref{intro:basis-growth}, $\beta_{K,n}$ does not depend on $X$. If
$\beta_{K,n}>0$, the $Q^s$-Borel transform of $\widehat F_{K,n}$ has
radius of convergence zero for $s<\beta_{K,n}$ and is entire for
$s>\beta_{K,n}$. Thus $Q^{\beta_{K,n}}$ is the only base at which the
Borel transform can have a positive but finite radius of convergence,
and we call $\beta_{K,n}$ the critical level.

\begin{theoremA}
Let $K$ be a knot and $n\ge2$. Then
\[
                  \beta_{K,n}=\max\{\delta_{K,n}-1,0\},
\]
and $\beta_{K,n}$ is rational. The Newton terms obey the bound at the
critical exponent itself: there is $R>0$ such that, for every compact
set $\mathcal K$ of colors,
\[
 \sup_{X\in\mathcal K}\bigl|c_{K,n,k}B_k^{(n)}(X)\bigr|
          \le C_{\mathcal K}R^kQ^{-\beta_{K,n}k^2/2}.
\]
If $\beta_{K,n}>0$, the critical $Q^{\beta_{K,n}}$-Borel transform has a
positive and finite radius of convergence, the same at every color that
is not a node, and it extends meromorphically to the whole plane, with
poles of bounded order on finitely many geometric progressions. At each
pole on its circle of convergence, the principal part is
$E_n(X)/E_n(0)$ times the principal part at $X=0$.
\end{theoremA}

The identity $\beta=\max\{\delta-1,0\}$ holds for arbitrary complex node
values (Theorem~\ref{thm:nodegrowth}) and follows from an explicit
formula for the Newton coefficients in terms of the node values. The
shift by one comes from the basis: on the diagonal,
$B_r^{(n)}(x_r^{(n)})$ has quadratic order two, one more than
$B_r^{(n)}(X)$ at the other colors. For arbitrary node values, however,
a finite quadratic order need not give a Borel transform with positive
radius at the critical base. For knots this is supplied by
$q$-holonomicity \cite{GLL,GLSurvey}: the Newton coefficients again
satisfy a linear $q$-difference equation, and Ramis's growth theorem for
entire solutions of such equations \cite{RamisGrowth} shows that the
critical level is rational and that the bound holds at the critical
exponent itself (Theorem~\ref{thm:critical-knot-borel}). The same
difference equation yields the meromorphic continuation
(Theorem~\ref{thm:allknots-meromorphic}). The factor $E_n(X)$ in the
principal parts vanishes exactly at the nodes, where the Borel transform
is a polynomial (Proposition~\ref{prop:first-singularities-node-factor}).
More generally, the factorization of the basis transports the continued
Borel transform from the color $X=0$ to all colors, with a set of
possible poles that does not depend on the color
(Proposition~\ref{prop:universal-color-transport}).

Reduced invariants are multiplicative under connected sum, so the node
order of a connected power is $\delta_{\#^MK,n}=M\delta_{K,n}$, and
Theorem~A gives $\beta_{\#^MK,n}=\max\{M\delta_{K,n}-1,0\}$
(Corollary~\ref{cor:general-connected}). For instance
$\beta_{\#^M4_1,n}=2M-1$, so every positive odd integer is the critical
level of some knot. Theorem~A concerns the Borel transform only. A
Laplace sum at the critical level requires, in addition, control of the
growth of the Borel transform along the integration ray. We give a
criterion in terms of the difference equation
(Proposition~\ref{prop:holonomic-growth-test}), which needs to be
checked only at the single color $X=0$
(Corollary~\ref{cor:one-color-growth-test}); holonomicity alone does not
provide this control (Remark~\ref{rem:growth}).

\subsection{Spectral summation and Stokes data}

For the knots treated in detail below much more structure is available:
their Newton coefficients have an \emph{analytic tail}. This means that,
from some index on,
\[
 c_k=\sum_{\nu=1}^s\lambda_\nu^kP_\nu(k)f_\nu(Q^k),
\]
where $\lambda_\nu\ne0$, the $P_\nu$ are polynomials and the $f_\nu$ are
holomorphic near zero. The class is natural because of a one-line
computation. If $c_k=\alpha^{-k}b(Q^k)$ with $b(t)=\sum_mb_mt^m$, then
expanding $b(Q^k)$ and summing the geometric series in $k$ gives
\[
 \sum_{k\ge0}\alpha^{-k}b(Q^k)\,\xi^k
   =\sum_{m\ge0}\frac{b_m}{1-Q^m\xi/\alpha}
   =:(\sRes_\alpha b)(\xi).
\]
The Taylor coefficients of the \emph{symbol} $b$ thus become residues:
the \emph{resolvent} $\sRes_\alpha b$ has simple poles at the points
$\alpha Q^{-m}$, with residues $-\alpha Q^{-m}b_m$. The $q$-Pochhammer
factors produced by the basis are themselves analytic functions of
$Q^k$. Hence, after a finite initial part is removed, the Borel
transform of a Newton series with an analytic tail is a finite
combination of resolvents and of their Euler derivatives, which account
for the polynomials $P_\nu$ (Theorem~\ref{thm:analytic-symbol}).

For resolvents, summation and Stokes data are explicit. Let $\W_R$,
$R>1$, be the Banach algebra of power series $b(t)=\sum_mb_mt^m$ with
$\sum_m|b_m|R^m<\infty$, and put $h_p=-\log p$ and
$w_p(t)=1/(-t;p)_\infty$.

\begin{theoremB}
Let $0<p<1$, $\alpha\ne0$, $R>1$ and $b\in\W_R$, and put
$H=\sRes_\alpha b$.
\begin{enumerate}[label=(\roman*),leftmargin=2em]
\item The Laplace sum $\cL_pH(z)$ is holomorphic for $z\ne0$ off the
ray $\alpha\mathbb R_{>0}$. On closed sectors avoiding this ray it has the
asymptotic expansion $\sum_k\alpha^{-k}b(Q^k)M_k(p)z^k$, and the
remainder after $N$ terms is at most $C\,M_N(p)|z/\alpha|^N$.
\item Both lateral limits of $\cL_pH$ exist on the ray. At
$z=e^{i\arg\alpha}r$, the limit from the side of decreasing argument
minus the limit from the side of increasing argument is
\[
 \Delta_b(z)=-\frac{2\pi i}{h_p}\sum_{m\ge0}b_m\,\tau_mw_p(\tau_m),
 \qquad \tau_m=\frac{|\alpha|Q^{-m}}{r}.
\]
\item The jump determines the symbol. If $b_{m_0}$ is the first nonzero
coefficient of $b$, then
$\Delta_b(z)\sim-(2\pi i/h_p)\,b_{m_0}\tau_{m_0}w_p(\tau_{m_0})$ as
$r\downarrow0$, and subtracting this term and repeating recovers every
$b_m$.
\end{enumerate}
\end{theoremB}

Theorem~B collects Theorems~\ref{gen:thm:resolvent},
\ref{gen:thm:sumjump} and~\ref{gen:thm:injective}, and
Proposition~\ref{gen:prop:jets} extends parts (i) and (ii) to Euler
derivatives. Within
the resolvent class, then, the Stokes jump, viewed as a function of the
modulus on the singular ray, loses no information: it determines every
pole and residue of the Borel transform. What the jump cannot see is a
polynomial added to the Borel transform, that is, a change of finitely
many coefficients. The calculus is also multiplicative. If
$g\in\W_R$ and $\T$ is the dilation $(\T f)(\xi)=f(Q\xi)$, then
$g(\T)\sRes_\alpha b=\sRes_\alpha(gb)$, and the same operator $g(\T)$,
acting on functions of $z$, transforms the Laplace sums and the jumps.

For Newton series with analytic tails, Theorem~\ref{thm:analytic-symbol}
turns this into statements about the color. The sum $S(X,z)$ is an
entire function of $X$ and equals the terminating Newton polynomial at
every node; it has lateral limits on every singular ray, locally
uniformly in the color; and its jump is a normally convergent sum of
residues.

Analytic tails are compatible with the divided differences of the
Askey--Wilson calculus \cite{IS}. The $q$-divided difference
$\mathcal D_q$ of \eqref{eq:Dq-color} acts on the Newton basis with node
parameter $A$ by $\mathcal D_qB_k^{(A)}=[k]_qB_{k-1}^{(qA)}$; it lowers
the index and, since $A=q^n$, raises the rank by one. Accordingly,
raising the rank of an analytic tail by $d$ amounts to applying
$z^{-d}\mathcal D_q^d/[d]_q!$ to its sum
(Theorem~\ref{thm:general-rank-raising}). For the figure-eight knot the
Newton coefficients are the $q$-binomial coefficients
\[
 c_{4_1,n,k}=\frac{[k+n-2]_q!}{[k]_q!\,[n-2]_q!},
\]
which are obtained in exactly this way from the constant sequence of
rank two. Hence the figure-eight sums of all ranks are generated by the
sum of rank two:
\begin{equation}\label{intro:rank-ladder}
 S_n(X,z)=\frac{z^{-(n-2)}}{[n-2]_q!}\,\mathcal D_q^{\,n-2}S_2(X,z),
\end{equation}
and the same identity holds for the lateral sums and their jumps
(Theorem~\ref{thm:rank-ladder}). We call \eqref{intro:rank-ladder} the
rank ladder.

\subsection{Twist knots and a positive kernel}\label{subsec:intro-twist}

Let $K_p$ be the twist knots in the convention of
Kawagoe~\cite{Kawagoe}, and put $L_m=\overline{K_{-m}}$ for $m\ge1$;
thus $L_1=4_1$ and $L_2=\overline{6_1}$. The Newton coefficients factor
as
\[
 c_{L_m,n,k}=D_{n,k}\,g_{m,k}^{(n)},\qquad
 D_{n,k}=\frac{[k+n-2]_q!}{[k]_q!\,[n-2]_q!},
\]
where $D_{n,k}$ is the figure-eight coefficient, $g_{1,k}^{(n)}=1$, and
\[
 g_{m,k}^{(n)}=\sum_{j=0}^k Q^{j(j+n-1)}\qb kj\,g_{m-1,j}^{(n)}
 \qquad(m\ge2)
\]
(Proposition~\ref{pa:knot-expansion}). The factor $g_{m,k}^{(n)}$ has a
simple combinatorial meaning. It is a sum over words
$w=(w_1,\ldots,w_k)\in\{1,\ldots,m\}^k$,
\[
 g_{m,k}^{(n)}=\sum_{w}Q^{E_{m,n}(w)},\qquad
 E_{m,n}(w)=\operatorname{inv}(w)
        +\sum_{s=1}^{m-1}N_s(w)\bigl(N_s(w)+n-1\bigr),
\]
where $\operatorname{inv}(w)$ is the number of inversions of $w$ and
$N_s(w)$ is the number of letters of $w$ that are at most $s$
(Lemma~\ref{pa:energy}). In particular $1\le g_{m,k}^{(n)}\le m^k$.

Two analytic properties of these coefficients drive the asymptotic
analysis. First, $g_{m,k}^{(n)}$ is the value at $Q^k$ of an entire
function (Lemma~\ref{pa:interpolation}), so the Newton series of $L_m$
has an analytic tail and the results above apply. Its Borel
transform at base $Q$, which we denote by $H_{m,n,h}(X,\xi)$, has
possible poles only at $\xi=-Q^{n-2-j}$, $j\ge0$, all of them simple,
and its sum at $z=1$,
\[
 S_{m,n,h}(X)=\frac1h\int_0^\infty
     \frac{H_{m,n,h}(X,t)}{(-t;Q)_\infty}\,dt,
 \qquad Q=e^{-h},
\]
is an entire function of the color with
$S_{m,n,h}(x_r^{(n)})=J_r^{(n)}(L_m;q)$ for all $r$
(Theorem~\ref{pa:source}). Second, the coefficients have a positive
spectral representation. Write $g_{m,k}(A;Q)$ for the coefficient
defined by the same recursion with $Q^n$ replaced by $A^2$, so that
$g_{m,k}^{(n)}=g_{m,k}(q^n;Q)$. For independent $0<a,Q<1$,
\[
 (a;Q)_k\,g_{m,k}(\sqrt a;Q)
       =\sum_{\ell\ge0}p_{m,\ell}^{(a)}Q^{k\ell},
 \qquad p_{m,\ell}^{(a)}\ge0,\quad
       \sum_\ell p_{m,\ell}^{(a)}=1
\]
(Proposition~\ref{pa:positive-moments}). The probability distributions
$p_m^{(a)}$ are generated by an explicit Markov kernel, dual to the
$q$-beta-binomial weights of \cite{CMPPush}, which we obtain from a
confluent form of Heine's transformation \cite{GasperRahman,DLMFHeine}.

At $a=Q^n$ this representation, combined with a rank factor whose
weights have bounded total variation, gives an exact formula for the
error made by truncating the Newton series (Theorem~\ref{pa:source}).
For $K\ge1$ and $|X-2|<2\cosh((K+n/2)h)-2$, the difference between
$S_{m,n,h}(X)$ and its Newton polynomial of degree $K-1$ is the $K$th
Newton term multiplied by a linear combination of the averages
\[
 \mathbb E\,\frac1{1+Q^{2-n+j}\,T_K},\qquad j=0,1,2,\ldots,
\]
where $T_K$ is a positive random variable with density proportional to
$t^K/(-t;Q)_\infty$. The coefficients of the combination add up to one
and have explicitly bounded total variation. Each average lies between
zero and one, so the remainder is bounded by the next Newton term times
explicit factors. At rank two and for real colors with
$0\le X-2<2\cosh((K+1)h)-2$, all coefficients are nonnegative and the
normalized remainder is a Stieltjes transform of a probability measure.
The bounds are uniform in $h$, so the formula reduces the analytic
function $S_{m,n,h}$ to a finite Newton sum, with explicit control, in
all the ranges of $h$ and $X$ used below.

The twist factor also acts on Stokes data. The Borel transform of $L_m$
is that of the figure-eight knot with its symbol multiplied by the
entire function interpolating $g_{m,k}^{(n)}$, and the effect is
visible already at the first pole: for $1\le i<j$, the jumps of $L_i$
and $L_j$ differ at small modulus by an explicit nonzero leading term
(Proposition~\ref{merged:twist-separation}). In particular, at every
fixed rank and every color that is not a node, the Stokes jumps
distinguish all the knots $L_m$.

\subsection{Critical limits at Alexander roots}

We now let $q\to1$, writing $Q=e^{-h}$ and $q=e^{-h/2}$ with
$h\to0^+$. For each fixed $k$,
\[
 D_{n,k}\to\binom{k+n-2}{n-2},\qquad
 g_{m,k}^{(n)}\to m^k,\qquad
 B_k^{(n)}(X)\to(X-2)^k,
\]
so, term by term, the Newton series of $L_m$ tends to
\[
 \sum_{k\ge0}\binom{k+n-2}{n-2}\bigl[m(X-2)\bigr]^k
                  =\bigl[1-m(X-2)\bigr]^{-(n-1)}.
\]
In the trace coordinate $X=t+t^{-1}$, the function
$\Delta_m(X)=1-m(X-2)$ is the Alexander polynomial of $L_m$, normalized
to take the value one at $t=1$. At rank two the limit is the inverse
Alexander polynomial that appears in the Melvin--Morton--Rozansky
expansion \cite{GLAsymptotics,GLAnalytic}. The geometric series
converges only for $|m(X-2)|<1$, and its sum blows up at the root
\[
 X_m=2+\frac1m .
\]
What replaces it near $X_m$ is decided by the first decaying factor
that survives the limit, and there are two different mechanisms.

For $m\ge2$ the cutoff comes from the twist coefficient. In the word
model, $g_{m,k}^{(n)}/m^k$ is the average of $e^{-hE_{m,n}(w)}$ over
uniformly random words, and the energy concentrates around its mean
$\alpha_mk^2+O(k)$, where
\[
 \alpha_m=\frac{(m-1)(4m+1)}{12m}.
\]
Hence $g_{m,k}^{(n)}/m^k\approx e^{-\alpha_mhk^2}$, a Gaussian cutoff at
$k$ of order $h^{-1/2}$. For the figure-eight knot the twist coefficient
is identically one, and the first decay comes instead from the Newton
product. At the root $X_1=3$,
\[
 \log B_k^{(n)}(3)=\sum_{j<k}\log\bigl(1-\tau((j+n/2)h)\bigr)
       \approx-\frac{h^2k^3}{3},\qquad \tau(v)=2\cosh v-2,
\]
a cubic cutoff at $k$ of order $h^{-2/3}$. Summing the Newton terms in
the windows $X=X_m+\sqrt h\,s$ and $X=3+h^{2/3}s$ leads to Riemann sums
for the entire functions
\begin{align}
 \mathcal G_{m,n}(s)
   &=\frac1{(n-2)!}\int_0^\infty
        t^{n-2}e^{mst-\alpha_mt^2}\,dt
             \qquad(m\ge2),\label{pa:gaussian-profile}\\
 \mathcal H_n(s)
   &=\frac1{(n-2)!}\int_0^\infty
        t^{n-2}e^{st-t^3/3}\,dt .\label{pa:scorer-profile}
\end{align}

\begin{theoremC}
Fix $m$ and $n$. As $h\to0^+$,
\begin{align*}
 h^{(n-1)/2}S_{m,n,h}(X_m+\sqrt h\,s)
       &=\mathcal G_{m,n}(s)+O(\sqrt h), &&m\ge2,\\
 h^{2(n-1)/3}S_{1,n,h}(3+h^{2/3}s)
       &=\mathcal H_n(s)+O(h^{2/3}).
\end{align*}
The errors are uniform for $s$ in compact subsets of $\mathbb C$, and
the same estimates hold for every fixed number of $s$-derivatives.
\end{theoremC}

To make the heuristics rigorous we control the analytic sum outside the
critical window by the remainder formula, and the twist coefficient by
concentration estimates for the word energy
(Theorems~\ref{pa:gaussian} and~\ref{pa:scorer}). The profiles are
classical special functions. At rank two, $\mathcal G_{m,2}$ is an
error-function integral, and $\mathcal H_2=\pi\operatorname{Hi}$, where
$\operatorname{Hi}(s)=\pi^{-1}\int_0^\infty e^{st-t^3/3}\,dt$ is the
Scorer function, a solution of the inhomogeneous Airy equation
$y''-sy=1/\pi$. The higher-rank Scorer profiles are derivatives of the
rank-two one,
\[
 \mathcal H_n=\frac{\partial_s^{\,n-2}\mathcal H_2}{(n-2)!},
\]
and this is what the rank ladder \eqref{intro:rank-ladder} becomes in
the critical scaling, where the $q$-divided difference turns into
$\partial_s$. We use this to deduce the Scorer limit at all ranks from
rank two. A rotation of $s$ by a cube root of unity relates
$\mathcal H_n$ to derivatives of the Airy function
\eqref{pa:Airy-profile}, and by Theorem~C the same relation holds
asymptotically for the analytic sums at complex colors.

At the quantum colors the theorem gives growth laws for the invariants
themselves. Let $u_m=\operatorname{arccosh}(1+1/(2m))$ and
$h_r=u_m/(r+n/2)$. Then the node $x_r^{(n)}$ is exactly $X_m$, and
Theorem~C at $s=0$ gives, as $r\to\infty$,
\begin{align*}
 J_r^{(n)}(L_m;e^{-h_r/2})
 &\sim\frac{\Gamma((n-1)/2)}
 {2(n-2)!\,\alpha_m^{(n-1)/2}}\,h_r^{-(n-1)/2}
 &&(m\ge2),\\
 J_r^{(n)}(4_1;e^{-h_r/2})
 &\sim\frac{3^{(n-4)/3}\,\Gamma((n-1)/3)}{(n-2)!}\,h_r^{-2(n-1)/3},
\end{align*}
with relative errors $O(h_r^{1/2})$ and $O(h_r^{2/3})$
(Corollary~\ref{pa:critical-nodes}). At rank two these are growth laws
of order $\sqrt N$ and $N^{2/3}$ in the dimension $N=r+1$ of the
representation. For the figure-eight knot at rank two, polynomial
growth at an Alexander root was proved by Hikami and Murakami~\cite{HM}.

On the real axis below the root, the inverse Alexander expression is
again a good approximation, and the analytic sums interpolate between
the two regimes. Put $\eta=1-m(X-2)$, so that the classical limit is
$\eta^{-(n-1)}$. For fixed $m$ and $n$, uniformly in the indicated
ranges,
\begin{align*}
 \eta^{n-1}S_{m,n,h}(X)
 &=1-\alpha_m n(n-1)\frac h{\eta^2}
      +O\Big(\frac h\eta+\frac{h^2}{\eta^4}\Big),
      &&\sqrt h\le\eta\le\tfrac18,\quad m\ge2,\\
 \eta^{n-1}S_{1,n,h}(3-\eta)
 &=1-\frac{(n-1)n(n+1)}3\,\frac{h^2}{\eta^3}
      +O\Big(\frac{h^2}{\eta^2}+\frac{h^4}{\eta^6}\Big),
      &&h^{2/3}\le\eta\le\tfrac18
\end{align*}
(Theorem~\ref{bridge:analytic-abelian-overlap}). The left tails of the
profiles have exactly the same first corrections
(Proposition~\ref{bridge:profile-left-tails}), so the analytic sums
match $\mathcal G_{m,n}(-u)$ and $\mathcal H_n(-u)$, first correction
included, in the simultaneous limit $u\to\infty$, $\sqrt h\,u\to0$
(respectively $h^{2/3}u\to0$)
(Corollary~\ref{bridge:actual-critical-matching}). On the real axis
below the root, the inverse Alexander expression and the critical
profile thus describe the analytic interpolation on overlapping ranges.

\subsection{Growing rank}

The profiles depend on the rank, and it is natural to ask how large the
rank may be before Theorem~C fails. In the Gaussian case the critical
window is centered at $k\approx\sqrt{n/(2\alpha_mh)}$, and the
approximations made at fixed rank break down when $hn^3$ is of order
one. In the Scorer case the window is centered at $k\approx(n/h^2)^{1/3}$,
and the corresponding parameter is $h^{2/3}n^{5/3}$. On these scales a
correction factor appears. Put
\[
 c_m=-\frac{(4m+1)(3m^3+3m^2+8m-2)}
               {180m^2(2\alpha_m)^{3/2}}<0 .
\]

\begin{theoremD}
Let $n=n(h)\to\infty$ as $h\to0^+$, with $m$ fixed. If $m\ge2$ and
$hn^3=O(1)$, then
\[
 \frac{h^{(n-1)/2}S_{m,n,h}(X_m+\sqrt h\,s)}
             {\mathcal G_{m,n}(s)}
       -e^{c_m\sqrt h\,n^{3/2}}\longrightarrow0.
\]
If $m=1$ and $h^{2/3}n^{5/3}=O(1)$, then
\[
 \frac{h^{2(n-1)/3}S_{1,n,h}(3+h^{2/3}s)}
             {\mathcal H_n(s)}
       -e^{-3h^{2/3}n^{5/3}/40}\longrightarrow0.
\]
The convergence is uniform for $s$ in compact subsets of $\mathbb C$,
and on such a set the denominators have no zeros once $h$ is small.
\end{theoremD}

The exponent of the correction collects the terms of order
$\sqrt h\,n^{3/2}$, respectively $h^{2/3}n^{5/3}$, at the center of the
critical window. In the Gaussian case they come from the rank factor
$D_{n,k}$, from the rank-dependent part of the mean word energy, from
the variance of the energy, whose leading coefficient is
\[
 V_m=\frac{(m-1)(m+1)(16m^2+1)}{180m^2},
\]
and from the cubic term of the Newton product; the higher cumulants of
the energy are controlled by a dependency-graph bound
\cite{FMNCumulants}. In the Scorer case they come from the rank factor
and the Newton product alone
(Theorems~\ref{pa:growing-rank-gaussian} and~\ref{pa:growing-rank-scorer}).
Consequently the fixed-rank limits remain valid, with relative error
$o(1)$, as long as $n=o(h^{-1/3})$ in the Gaussian case and
$n=o(h^{-2/5})$ in the Scorer case. These ranges are sharp: for
$n\sim ah^{-1/3}$ and $n\sim ah^{-2/5}$ the limiting ratios at $s=0$ are
$e^{c_ma^{3/2}}$ and $e^{-3a^{5/3}/40}$. In both regimes $nh\to0$, so
that $q^n\to1$.

\subsection{Logarithmic monodromy and geometric periods}

So far $q$ has been real. To continue in the deformation parameter we
use the coordinate $Y=Q^{n/2}X$, which removes the half-integral powers
of $Q$ from the Borel transform. For complex $Q$ with $0<|Q|<1$ the
Borel transform is then meromorphic in the Borel variable, and its
possible poles, all simple, lie at $-Q^{n-2-j}$, $j\ge0$
(Lemma~\ref{pa:complex-continuation}).
The positive real axis is no longer a usable Laplace contour: for
suitable complex $h=-\log Q$ arbitrarily close to zero in any sector
around the positive axis, the kernel $1/(-t;Q)_\infty$ has a pole on it
(Proposition~\ref{pa:fixed-ray}). Instead we integrate over the
logarithmic spiral $\{e^{hx}:x\in\mathbb R\}$, which stays away from the
poles of the kernel and of the Borel transform. This gives a
holomorphic continuation $\widehat S_{m,n}(h,Y)$ of the sums to the
half-plane $\Re h>0$ (Theorem~\ref{pa:spiral-theorem}). The
continuation is not single-valued in $Q=e^{-h}$, but its dependence on
$\log Q$ is very simple.

\begin{theoremE}
There are functions $\mathcal A_{m,n}(Q,Y)$ and $\mathcal B_{m,n}(Q,Y)$,
single-valued and holomorphic for $0<|Q|<1$ and entire in $Y$, such that
\[
 \widehat S_{m,n}(h,Y)
       =\mathcal A_{m,n}(e^{-h},Y)
                       +\frac1h\,\mathcal B_{m,n}(e^{-h},Y).
\]
The coefficient $\mathcal B_{m,n}$ vanishes at the nodes
$Y=Q^{-r}+Q^{n+r}$, $r\ge0$. For real $0<Q<1$,
\[
 \mathcal B_{m,n}(Q,Q^{n/2}X)
               =\sum_{\ell\ge0}Q^\ell\rho_\ell R_\ell(X),
 \qquad
 \rho_\ell=\frac{(-1)^\ell Q^{\ell(\ell-1)/2}}
                   {(Q;Q)_\ell(Q;Q)_\infty},
\]
where $R_\ell(X)$ is the residue of $H_{m,n,h}(X,\xi)$ at
$\xi=-Q^{-\ell}$, and the series converges normally on compact sets of
colors. The Stokes jump at fixed $Q$, as a function of the modulus,
determines all these residues and hence $\mathcal B_{m,n}$.
\end{theoremE}

Since $\mathcal A_{m,n}$ and $\mathcal B_{m,n}$ are single-valued in
$Q$, the factor $1/h$ is the only source of multivaluedness, and the
deck transformation $h\mapsto h+2\pi i$ acts by
\[
 \widehat S_{m,n}(h+2\pi i,Y)-\widehat S_{m,n}(h,Y)
 =-\frac{2\pi i}{h(h+2\pi i)}\,\mathcal B_{m,n}(e^{-h},Y)
\]
(Theorem~\ref{pa:monodromy}). We call $\mathcal B_{m,n}$ the
coefficient of the logarithmic monodromy. The decomposition comes from
the residue theorem applied to the Laplace integrand multiplied by
$\log(-t)$: the poles of the Borel transform at $-Q^{-\ell}$, $\ell\ge0$,
coincide with poles of the kernel, and the resulting double poles
produce the term $\mathcal B_{m,n}/h$. At rank two and for the
figure-eight knot this
decomposition was found by Kashaev, Mangazeev and Mari\~no~\cite{KMM}.
Theorem~E shows that the two discontinuities of the problem, the Stokes
jump in the auxiliary variable $z$ at fixed $Q$ and the logarithmic
monodromy in $h$, are two different linear functionals of the same
sequence of Borel residues
(Proposition~\ref{bridge:residues-stokes-logarithmic-block}).

The coefficient $\mathcal B_{m,n}$ also has a closed contour-integral
representation, in which Euler's expansion of $(w^{-1};Q)_\infty$
assembles the Borel residues (Theorem~\ref{r4:closed-block}). For the
figure-eight knot at the color $X=2$, the integrand is a ratio of
$q$-Pochhammer symbols whose classical phase has the two critical
values $\pm iV$, where
\[
 V=2\operatorname{Im}\operatorname{Li}_2(e^{i\pi/3})
       =\operatorname{Vol}(S^3\setminus4_1)
\]
is the hyperbolic volume of the complement. Writing $d=n-2$, we prove
that in each closed sector $\delta\le\arg h\le\pi/2-\delta$,
$0<\delta<\pi/4$,
\[
 \mathcal B_{1,n}(e^{-h},2e^{-nh/2})
 =\frac{(-1)^d h^{1/2-2d}e^{iV/h-i\pi/4}}
     {d!\,\Gamma(n/2)^2\sqrt{2\pi}\,3^{1/4}}
 \left(1-\frac{i(27d^2+11)}{72\sqrt3}\,h
                                  +O(|h|^2)\right),
\]
and that the formula in the lower sector is obtained by complex
conjugation (Theorem~\ref{r4:sector-period}). At rank two the action and
the first two coefficients agree with the geometric perturbative
periods computed in \cite{KMM}, and the leading factor is the
meridian-normalized adjoint torsion factor $\sqrt{-h/(4\pi\tau_\mu)}$ of
the geometric representation (Appendix~\ref{app:torsion}). At higher
rank the integrand is the rank-two one multiplied by an explicit factor,
and raising the rank by two multiplies it by a Laurent polynomial
(Corollary~\ref{r4:rank-raising}). The geometric action is therefore the
same at every rank; the rank changes only the power of $h$ and the
coefficients. On the positive real axis the behavior at $X=2$ is quite
different: there the remainder formula gives a complete asymptotic
expansion of $S_{m,n,h}(2)$ in powers of $h$
(Proposition~\ref{r4:parabolic-abelian}). Finally, the node divisor
$E_n$, which divides $\mathcal B_{m,n}$, has a simple expression in
terms of the dual nome $e^{-4\pi^2/h}$
(Appendix~\ref{app:modular-divisor}).

\subsection{Relation to earlier work}

Newton interpolation and its divided differences belong to the
Askey--Wilson calculus of Ismail and Stanton~\cite{IS}. The Newton
expansions of colored invariants that are uniform in the rank, and the
completion in which they live, are developed in~\cite{FangZhou,Fang}.
For twist knots, cyclotomic formulas go back to Masbaum~\cite{Masbaum};
we use Kawagoe's formula for the symmetric colored HOMFLY--PT
polynomials~\cite{Kawagoe} together with the differential expansions
of~\cite{IMMM,MorozovTselousov,CLZDoubleTwist}. The finite normalized
weights of the twist family are $q$-beta-binomial
weights~\cite[Appendix~A.2]{CMPPush}, and Heine's
transformation~\cite{GasperRahman,DLMFHeine} produces the infinite dual
kernel needed for the remainder formula.

The summation moments are those of Kashaev, Mangazeev and
Mari\~no~\cite{KMM}, whose analytic quantum invariants are the starting
point of this work; for general $q$-summation and the local theory of
$q$-difference equations see~\cite{DVZ,RSZ,Dreyfus}. The critical level
for arbitrary knots combines the explicit Newton inversion with the
$q$-holonomicity of colored invariants~\cite{GLL,GLSurvey} and Ramis's
growth theorem~\cite{RamisGrowth}.

On the asymptotic side, Hikami and Murakami proved polynomial growth of
the figure-eight invariant at an Alexander root at rank two~\cite{HM},
Murakami and Tran studied the real exponential asymptotics beyond this
point~\cite{MurakamiTran}, and Wong and Au studied fixed even ranks in a
different complex regime~\cite{WongAu}. The inverse-Alexander limit near
the trivial parameter is part of the analytic Melvin--Morton--Rozansky
theory~\cite{GLAsymptotics,GLAnalytic}. Our critical estimates concern
the analytic interpolation rather than the invariants alone, hold for
complex colors, and come with error bounds and with the corrections for
growing rank.

At rank two, the logarithmic decomposition of the figure-eight sum and
the associated geometric series appear in~\cite{KMM}, and the same
$q$-series occurs in the factorization of state integrals by Garoufalidis
and Kashaev~\cite{GKState}. Fantini and Wheeler prove summability and
Stokes relations for nonabelian state integrals, including those of
$4_1$ and $5_2$~\cite{FantiniWheeler}; the broader resurgent setting is
described in~\cite{Witten,GGM,GGMW}. Two different discontinuities
appear in the present paper. The Stokes phenomenon takes place in the
auxiliary variable $z$ at fixed $q$, whereas the logarithmic monodromy
is a deck transformation in $h=-\log Q$. Theorem~E relates the two
through the Borel residues, and the geometric asymptotics above concern
the coefficient $\mathcal B_{m,n}$ of the second.

\subsection{Organization}

Section~\ref{sec:foundations} sets up the summation method and proves
the results on critical levels for arbitrary knots.
Section~\ref{sec:spectral} develops spectral summation, the Stokes jump
and rank raising. Section~\ref{pa:data} treats the twist knots: their
Newton coefficients, the positive kernel, the remainder formula and the
separation of their Stokes data. Section~\ref{a:fixed-critical} proves
the Gaussian and Scorer limits and the matching with the
inverse-Alexander regime, and Section~\ref{pa:growing-rank} lets the
rank grow. Section~\ref{pa:complex} continues the sums in the
deformation parameter and studies the logarithmic monodromy and its
geometric asymptotics. Appendix~\ref{app:boundary} contains the
boundary-value estimates for lateral sums, Appendix~\ref{app:torsion}
the comparison with adjoint torsion, and
Appendix~\ref{app:modular-divisor} the modular form of the node divisor.

\section{Newton interpolation and critical Borel levels}
\label{sec:foundations}\label{sec:continuation}

This section concerns an arbitrary knot at a fixed rank. We first
show that the quadratic growth of a Newton series is determined by the
growth of its values at the quantum colors
(Theorem~\ref{thm:nodegrowth}). For colored knot invariants,
$q$-holonomicity then implies that the critical Borel level is rational
and attained, and that the critical Borel transform continues
meromorphically to the plane (Theorem~\ref{thm:critical-knot-borel}).
The factorization of the Newton basis carries this continuation,
together with its growth estimates, from the single color $X=0$ to all
colors (Proposition~\ref{prop:universal-color-transport}).

We use $\mathbb N_0=\{0,1,2,\ldots\}$. For $0<p<1$ let
\[
 (a;p)_k=\prod_{j=0}^{k-1}(1-ap^j),\qquad
 (a_1,\ldots,a_s;p)_k=\prod_{i=1}^s(a_i;p)_k,
\]
with the empty product equal to one. The subscript $\infty$ denotes
the convergent infinite product. We use the Gaussian binomial
coefficients
\[
 \Gauss{k}{j}=\frac{(Q;Q)_k}{(Q;Q)_j(Q;Q)_{k-j}}
 \quad(0\le j\le k),
\]
which are zero for other $j$, and $[s^j]f$ denotes the coefficient of
$s^j$ in a power series $f$.

First let $q$ be an indeterminate and $R=\mathbb Z[q^{\pm1}]$. The
Newton completion of \cite{Fang} is defined by
\begin{equation}\label{eq:nodes}
 x_j^{(n)}=q^{n+2j}+q^{-n-2j},\quad
 B_k^{(n)}(X)=\prod_{j=0}^{k-1}(X-x_j^{(n)}),\quad
 \mathcal C_n=\varprojlim_k R[X]/(B_k^{(n)}).
\end{equation}
Every element of $\mathcal C_n$ has a unique expansion
\begin{equation}\label{eq:Newton}
 F=\sum_{k\ge0}c_k(q)B_k^{(n)}(X),\qquad c_k(q)\in R.
\end{equation}
Indeed, $B_k^{(n)}$ is monic of degree $k$, so $B_0^{(n)},\ldots,
B_{k-1}^{(n)}$ form a basis of the polynomials of degree less than $k$.
An element of $R[X]/(B_k^{(n)})$ is represented by a unique remainder of
degree less than $k$, and this remainder has unique coordinates in that
basis. Passing from $k$ to $k+1$ keeps the first $k$ coordinates, which
gives \eqref{eq:Newton}. At the node $x_r^{(n)}$ only the terms with
$k\le r$ are nonzero.

For the analytic results we specialize the coefficients and fix
\begin{equation}\label{eq:parameters}
 0<q<1,\qquad Q=q^2,\qquad h_Q=-\log Q>0.
\end{equation}
It is convenient to replace $q^n$ by an independent parameter
$0<A<1$ and to write
\[
 x_j^{(A)}=AQ^j+A^{-1}Q^{-j},\qquad
 B_k^{(A)}(X)=\prod_{j=0}^{k-1}(X-x_j^{(A)}).
\]
The rank-$n$ basis is the case $A=q^n$, and this value of $A$ is
understood whenever no other value is specified. We write
$X(y)=y+y^{-1}$ and
\[
 E_A(X)=(Ay,A/y;Q)_\infty
       =\prod_{j\ge0}(1-AQ^jX+A^2Q^{2j}),\qquad E_n=E_{q^n}.
\]
Expressions such as $E_A$ that are symmetric under $y\mapsto y^{-1}$
depend only on $X$. For a knot $K$, let $J_r^{(n)}(K;q)$ be its reduced, zero-framed
quantum invariant colored by $\operatorname{Sym}^r(\mathbb C^n)$.
Its Newton coefficients are characterized by
\[
 J_r^{(n)}(K;q)=\sum_{k=0}^r
 c_{K,n,k}(q)B_k^{(n)}(x_r^{(n)}).
\]
The cyclotomic expansions of \cite{FangZhou} give the corresponding
Newton element $F_{K,n}$; the analytic arguments below use its
coefficients after specialization. For any Newton series $F$, its
auxiliary series is
\begin{equation}\label{eq:suspension}
 \widehat F_n(X,z)=\sum_{k\ge0}c_k(q)B_k^{(n)}(X)z^k.
\end{equation}
We write $\widehat F_A$ for the same construction with the independent
node parameter $A$. When the sum constructed below is defined at $z=1$,
its value there is a function of the color that takes the values of $F$
at the nodes. We write $\Omega=\mathbb C\setminus(-\infty,0]$.

\subsection{The summation method}

For real $0<p<1$ put
\begin{equation}\label{eq:moment}
 M_k(p)=(p;p)_k p^{-k(k+1)/2},\qquad
 w_p(t)=\frac1{(-t;p)_\infty},\qquad h_p=-\log p.
\end{equation}
We use the transforms of \cite[\S2]{KMM}:
\begin{equation}\label{eq:KMM}
 \cB_p\left(\sum a_kz^k\right)(\xi)
 =\sum_{k\ge0}\frac{a_k}{M_k(p)}\xi^k,\qquad
 \cL_p g(z)=\frac1{h_p}\int_0^\infty g(zt)w_p(t)\,dt,
\end{equation}
and write $\mathcal S_p=\mathcal L_p\mathcal B_p$ whenever the
integral exists.

\begin{definition}\label{def:admissible}
An auxiliary Newton series is \emph{meromorphically $p$-summable on an
open sector} if its Borel transform $\cB_p\widehat F$ converges near
zero, continues to a meromorphic function on $\C$ that is holomorphic
on the sector, and grows at most polynomially on every closed
subsector. Its sum is then $\cL_p\cB_p\widehat F$ for $z$ in the
sector.
\end{definition}

The Laplace integral also converges under the weaker bound
$|g(\xi)|\le C\exp\{\kappa\log^2|\xi|\}$ with
$\kappa<1/(2h_p)$, as the following kernel estimate shows.

\begin{theorem}\label{thm:kernel}\label{pa:moment}
For fixed $0<p<1$,
\begin{align}
 \log(-t;p)_\infty
 &=\frac{(\log t)^2}{2h_p}+\frac{\log t}{2}+O_p(1)
       &&(t\to+\infty),\label{eq:kernel-growth}\\
 \frac1{h_p}\int_0^\infty t^kw_p(t)\,dt
 &=M_k(p) &&(k\ge0).\label{eq:moment-identity}
\end{align}
Let $g$ be holomorphic near zero and along a ray from the origin, and
suppose that $|g(\xi)|\le C\exp\{\kappa(\log|\xi|)^2\}$ on that ray for
$|\xi|\ge1$, with $\kappa<1/(2h_p)$. Then the Laplace integral of $g$
along the ray converges absolutely. In particular polynomial growth
suffices.
\end{theorem}

\begin{proof}
Put $x=\log t>0$ and $J=\lfloor x/h_p\rfloor$, so that
$\log(-t;p)_\infty=\sum_{j\ge0}\log(1+e^{x-h_pj})$. For $j\le J$ write
$\log(1+e^{x-h_pj})=(x-h_pj)+\log(1+e^{-(x-h_pj)})$. The second terms,
and the terms with $j>J$, are bounded by geometric series with ratio
$e^{-h_p}=p$, so together they are $O_p(1)$. The main part is
$\sum_{j=0}^J(x-h_pj)=x^2/(2h_p)+x/2+O_p(1)$, which proves
\eqref{eq:kernel-growth}.

The functional equation $w_p(pt)=(1+t)w_p(t)$ gives
$tw_p(t)=w_p(pt)-w_p(t)$, hence
\[
 \int_0^\infty w_p(t)\,dt
 =\int_0^\infty\frac{w_p(pt)-w_p(t)}t\,dt=h_p.
\]
To evaluate the last integral, integrate over $[\epsilon,R]$ and
substitute $u=pt$ in the first term. This leaves
$\int_{p\epsilon}^\epsilon w_p(u)\,du/u-\int_{pR}^R w_p(u)\,du/u$,
which tends to $w_p(0)\log(1/p)=h_p$ as $\epsilon\to0$ and
$R\to\infty$, because $w_p$ is continuous at zero and decays at
infinity. All integrals converge by \eqref{eq:kernel-growth}. If
$I_k=\int_0^\infty t^kw_p(t)\,dt$, the same functional equation gives
$I_{k+1}=\int_0^\infty t^k(w_p(pt)-w_p(t))\,dt=(p^{-k-1}-1)I_k$, and
$M_{k+1}(p)=(p^{-k-1}-1)M_k(p)$; induction proves
\eqref{eq:moment-identity}. Finally, with $t=e^x$, the bound on $g$
and \eqref{eq:kernel-growth} show that the integrand is at most
$\exp\{(\kappa-1/(2h_p))x^2+O(x)\}$ for large $x$, and the negative
quadratic coefficient makes the integral converge.
\end{proof}

\begin{definition}\label{def:slope}
The quadratic order of a coefficient sequence $(u_k)$ is
\[
 \beta(u)=\max\left\{0,\limsup_{k\to\infty}
           \frac{2\log|u_k|}{h_Qk^2}\right\},
\]
with $\log0=-\infty$; a sequence of finite support has order zero. If
$|u_k|=\exp(\frac12h_Q\beta k^2+O(k))$ with $\beta>0$, the Borel
transform at base $p=Q^\beta$ has positive and finite radius of
convergence, while at every base $Q^{\beta'}$ with $0<\beta'<\beta$ the
radius is zero.
\end{definition}

\subsection{The basis factorization and analytic growth}

The Newton basis factors into $q$-Pochhammer symbols. This gives its
size away from the nodes and a closed form of the Borel transform at
base $Q$.

\begin{theorem}\label{thm:bridge}
For $0<A<1$, choose $y\in\C^\times$ with $X=y+y^{-1}$, and set
$a=Ay$, $b=A/y$. Then
\begin{align}
 B_k^{(A)}(X)
 &=(-1)^k A^{-k}Q^{-k(k-1)/2}(a,b;Q)_k,\label{eq:factor}\\
 \cB_Q\widehat F_A(X,\xi)
 &=\sum_{k\ge0}c_k(q)\frac{(a,b;Q)_k}{(Q;Q)_k}
       \left(-\frac QA\xi\right)^k.\label{eq:bridge}
\end{align}
Both right-hand sides are unchanged when $y$ is replaced by $y^{-1}$.
For $A=q^n$ and $X$ not a node, one has
\begin{equation}\label{eq:basis-asymp}
 B_k^{(n)}(X)\sim
 (-1)^kq^{-k^2-(n-1)k}E_n(X),\qquad
 E_n(X)\ne0.
\end{equation}
\end{theorem}

\begin{proof}
For every $j$,
\[
 X-AQ^j-A^{-1}Q^{-j}
 =-A^{-1}Q^{-j}(1-AQ^jy)(1-AQ^j/y).
\]
Taking the product over $0\le j<k$ gives \eqref{eq:factor}. Dividing by
the moment leaves
\[
 \frac{B_k^{(A)}(X)}{M_k(Q)}
 =(-1)^k(Q/A)^k\frac{(a,b;Q)_k}{(Q;Q)_k},
\]
which proves \eqref{eq:bridge}. At $A=q^n$ the products
$(a,b;Q)_k$ converge to $E_n(X)$, and $(-1)^kA^{-k}Q^{-k(k-1)/2}
=(-1)^kq^{-k^2-(n-1)k}$. The limit $E_n(X)$ vanishes exactly when $a$ or
$b$ equals $Q^{-r}$ for some $r\ge0$, that is, when $X=x_r^{(n)}$.
\end{proof}

Thus the Newton basis has quadratic order one at base $Q$ at every
color that is not a node, and the rank $n$ enters only the linear term
of the exponent.

\subsection{Recovering growth from the node values}
\label{subsec:nodegrowth}
In this subsection $0<A,Q<1$ are fixed and $(J_r)_{r\ge0}$ is an
arbitrary complex sequence, with Newton coefficients determined by
$J_r=\sum_{k=0}^r c_kB_k^{(A)}(x_r^{(A)})$. We show that the quadratic
order of the Newton series is determined by that of the values $J_r$.

\begin{lemma}\label{lem:explicitinverse}
For every $k\ge0$,
\begin{equation}\label{eq:explicitinverse}
 c_k=A^kQ^{k(k+1)/2}\sum_{r=0}^k
 \frac{(-1)^{k-r}J_rQ^{r(r-1)/2}(1-A^2Q^{2r})}
 {(Q;Q)_r(Q;Q)_{k-r}(A^2Q^r;Q)_{k+1}}.
\end{equation}
For every compact color set $\mathcal K$ there is a constant
$C_{\mathcal K}$, independent of $k$, such that
\begin{equation}\label{eq:inversebound}
 \sup_{X\in\mathcal K}|c_kB_k^{(A)}(X)|
 \le C_{\mathcal K}Q^k\sum_{r=0}^k|J_r|Q^{r(r-1)/2}.
\end{equation}
\end{lemma}
\begin{proof}
The coefficient $c_k$ is the leading coefficient of the Lagrange
interpolant of degree $k$, so
$c_k=\sum_{r=0}^kJ_r/\prod_{j\ne r}(x_r-x_j)$. Since
$x_r-x_j=A^{-1}(Q^{-r}-Q^{-j})(1-A^2Q^{r+j})$, the product of the
differences is
\[
 \prod_{\substack{0\le j\le k\\j\ne r}}(x_r-x_j)
 =(-1)^{k-r}A^{-k}Q^{-k(k+1)/2-r(r-1)/2}
 \frac{(Q;Q)_r(Q;Q)_{k-r}(A^2Q^r;Q)_{k+1}}{1-A^2Q^{2r}},
\]
which proves \eqref{eq:explicitinverse}. The three Pochhammer symbols in
the denominator are bounded below by positive constants independent of
$k$ and $r$. By \eqref{eq:factor},
\[
 B_k^{(A)}(X)=(-1)^kA^{-k}Q^{-k(k-1)/2}
 \prod_{j<k}(1-AXQ^j+A^2Q^{2j}),
\]
and the last product is bounded uniformly in $k$ for $X\in\mathcal K$.
Multiplying, the powers of $A$ cancel and the powers of $Q$ leave
$Q^kQ^{r(r-1)/2}$, which proves \eqref{eq:inversebound}.
\end{proof}

With $h=h_Q$ define
\[
 \delta(J)=\max\left\{0,\limsup_{r\to\infty}
             \frac{2\log|J_r|}{hr^2}\right\},\qquad
 \beta(F)=\max\left\{0,\limsup_{k\to\infty}
             \frac{2\log|c_kB_k^{(A)}(X)|}{hk^2}\right\},
\]
where $X$ is any color that is not a node, $\log0=-\infty$, and a
sequence of finite support has order zero. By \eqref{eq:basis-asymp},
$\beta(F)$ does not depend on the choice of $X$.

\begin{theorem}\label{thm:nodegrowth}
For arbitrary complex node values, with infinite orders allowed,
\begin{equation}\label{eq:exactnodegrowth}
                 \beta(F)=\max\{\delta(J)-1,0\}.
\end{equation}
\end{theorem}
\begin{proof}
Suppose $\delta<\infty$. For every $\epsilon>0$ we have
$|J_r|\le C_\epsilon Q^{-(\delta+\epsilon)r^2/2}$, so the $r$-th term on
the right of \eqref{eq:inversebound} is at most
$C_\epsilon Q^{-(\delta+\epsilon-1)r^2/2+O(r)}$. If
$\delta+\epsilon>1$, the largest of these $k+1$ terms is the last one
up to a factor $e^{O(k)}$; otherwise they are bounded by $e^{O(k)}$.
Hence $\beta\le\max\{\delta+\epsilon-1,0\}$, and letting
$\epsilon\downarrow0$ gives $\beta\le\max\{\delta-1,0\}$.

Conversely, suppose $\beta<\infty$. At a fixed color that is not a node,
the product in \eqref{eq:factor} converges to a nonzero limit, so
\[
 |c_k|\le C_\epsilon A^kQ^{k(k-1)/2-(\beta+\epsilon)k^2/2}.
\]
At a node there is the exact formula
\[
 B_k^{(A)}(x_r^{(A)})
 =A^{-k}Q^{-rk}(Q^{r-k+1};Q)_k(A^2Q^r;Q)_k,
 \qquad k\le r,
\]
and both products on the right are at most one. Dividing the logarithm
of the $k$-th term of $J_r=\sum_kc_kB_k^{(A)}(x_r^{(A)})$ by $h$
therefore gives at most
\[
 \tfrac12(\beta+\epsilon-1)k^2+rk+\tfrac12k+O(1).
\]
This bound increases on $0\le k\le r$, so it is at most
$\tfrac12(\beta+\epsilon+1)r^2+O(r)$, and summing over $k\le r$ gives
$\delta\le\beta+1$. If $\delta>1$, this and the first part give
$\beta=\delta-1$; if $\delta\le1$, the first part gives $\beta=0$.
Finally, if $\delta=\infty$, the second part shows that $\beta$ cannot be
finite.
\end{proof}

\subsection{Holonomic continuation at a fixed rank}
\label{subsec:holonomic-continuation}

At a fixed rank the $q$-holonomicity of the colored invariants passes
through the Newton transform, so every Borel transform at a rational
level that converges near the origin continues meromorphically to the
plane (Theorem~\ref{thm:allknots-meromorphic}). Proposition~\ref{prop:holonomic-growth-test} gives a recurrence test
for the growth at infinity.

\begin{proposition}
\label{prop:newton-holonomic}
Fix an integer $n\ge2$.  With $q$ indeterminate, the sequences
$k\mapsto c_{K,n,k}(q)$ and $k\mapsto c_{K,n,k}(q)B_k^{(n)}(X)$ are
$q$-holonomic over $\mathbb C(q)$ and $\mathbb C(q,X)$, respectively.
Let $s=a/b>0$ be rational, with $a,b$ positive integers, put
$\rho=q^{1/b}$ and $p=Q^s=\rho^{2a}$.  Then
\[
 u_k(\rho,X)=\frac{c_{K,n,k}(\rho^b)
                         B_k^{(n)}(X)}{M_k(\rho^{2a})}
\]
is $\rho$-holonomic in $k$.  After specializing to any
$q_0\in(0,1)$ and any $X_0\in\mathbb C$, the numerical sequence still
satisfies a nonzero recurrence with coefficients polynomial in
$\rho_0^k$, where $\rho_0=q_0^{1/b}$.
\end{proposition}

\begin{proof}
The one-row colored HOMFLY sequence is $q$-holonomic by \cite{GLL}.
Passing to the reduced, zero-framed normalization multiplies it by
$q$-hypergeometric factors, which preserves holonomicity.  Its
specialization $A=q^n$ also has a nonzero recurrence.  To see this
without a genericity assumption on the specialization, write the
recurrence with coefficients in $\mathbb C[A^{\pm1},q^{\pm1},M^{\pm1}]$,
where $M$ stands for $q^r$, and clear denominators.  Divide all its
coefficients by the largest common power of $A-q^n$ before making the
substitution.
The resulting coefficient polynomials do not all vanish under the
substitution.  Each colored invariant is regular there as a rational
function of the indeterminate $q$, so the specialized relation holds.

For the Newton transform use Lemma~\ref{lem:explicitinverse}.  At
$A=q^n$ its kernel, including the factors outside the sum, is
\[
 K_n(k,r)=
 \frac{(-1)^{k-r}q^{nk}
       Q^{k(k+1)/2+r(r-1)/2}(1-Q^{n+2r})
       (Q;Q)_{n+r-1}}
      {(Q;Q)_r(Q;Q)_{k-r}(Q;Q)_{n+r+k}},
       \qquad 0\le r\le k.
\]
Thus $c_{K,n,k}=\sum_{r=0}^kK_n(k,r)J_r^{(n)}(K;q)$.
The displayed kernel is a proper $q$-hypergeometric term: its
factorial arguments are integral affine forms in $(k,r)$ and its
monomial exponent is quadratic.  Extended by zero outside
$0\le r\le k$, it is $q$-holonomic; the restriction is given by the
two holonomic step functions imposing $r\ge0$ and $k-r\ge0$.
For each $k$ the sum over $r$ is finite, so closure under products and
finite sums \cite[Theorems~5.2 and~5.3]{GLSurvey} proves the first
assertion.

The basis has the first-order recurrence
\[
 B_{k+1}^{(n)}(X)=
       (X-q^{n+2k}-q^{-n-2k})B_k^{(n)}(X),
\]
so multiplication by it preserves holonomicity, including at its
zeros.  After $q=\rho^b$, every coefficient of a $q$-recurrence is
polynomial in an integral power of $\rho^k$.  Moreover
\[
 \frac{M_k(p)}{M_{k+1}(p)}
       =\frac{p^{k+1}}{1-p^{k+1}}
       =\frac{\rho^{2a(k+1)}}{1-\rho^{2a(k+1)}},
\]
which proves the assertion concerning rational Borel levels.

Finally clear denominators in a recurrence for $u_k$ so that its
coefficients are polynomials in $(\rho,X,M)$, where now $M$ stands for
$\rho^k$ (negative powers are cleared by a monomial factor).  Divide
out the largest common power of
$\rho-\rho_0$, and then specialize $\rho=\rho_0$.  All $u_k$ are
regular at this point: the possible denominators in the Newton formula
and in the moments are nonzero for $0<q_0<1$.  Next divide out the
largest common power of $X-X_0$ and specialize $X=X_0$.  The $u_k$ are
polynomial in $X$, so this second specialization is also legitimate.
The coefficient polynomials in $M$ are now not all zero.  A nonzero
polynomial cannot vanish at all the distinct points $\rho_0^k$; hence
the resulting recurrence is nonzero as a recurrence of numerical
sequences.  Finitely many exceptional initial indices may be omitted.
\end{proof}

\begin{lemma}
\label{lem:holonomic-meromorphy}
Let $0<\rho<1$, and suppose that a sequence $(u_k)_{k\ge0}$ has a
nonzero recurrence with coefficients polynomial in $\rho^k$.
If $H(\xi)=\sum_{k\ge0}u_k\xi^k$ converges near zero, there are
polynomials $D_0,\ldots,D_e,R$, with $D_0\ne0$, such that
\begin{equation}\label{eq:holonomic-germ-equation}
 D_0(\xi)H(\xi)+\sum_{\ell=1}^e
                 D_\ell(\xi)H(\rho^\ell\xi)=R(\xi).
\end{equation}
The function $H$ has a unique meromorphic continuation to $\mathbb C$.
If $\mathcal Z$ is the finite set of nonzero roots of $D_0$, its poles
lie in
\begin{equation}\label{eq:holonomic-pole-spirals}
                 \bigcup_{\alpha\in\mathcal Z}
                    \alpha\rho^{-\mathbb N_0}.
\end{equation}
For each multiplicative orbit $\mathcal O=\alpha\rho^{\mathbb Z}$,
every pole in that orbit has order at most
\[
       \sum_{\alpha'\in\mathcal Z\cap\mathcal O}
                              \operatorname{mult}_{\alpha'}D_0.
\]
In particular all pole orders are at most $\deg D_0$.
\end{lemma}

\begin{proof}
Write the recurrence, valid for $k\ge k_0$, as
\[
 \sum_{j=0}^{d}\sum_{\ell=0}^{e}
                 r_{j\ell}\rho^{\ell k}u_{k+j}=0.
\]
Common powers of $\rho^k$ have been removed, so some $r_{j0}\ne0$.
Multiplying by $\xi^{k+d}$ and summing gives
\eqref{eq:holonomic-germ-equation}, where
\[
 D_\ell(\xi)=\sum_{j=0}^d r_{j\ell}\rho^{-\ell j}\xi^{d-j}.
\]
The omitted initial coefficients give the polynomial $R$.  The
different powers $\xi^{d-j}$ ensure that $D_0$ is not zero.

If $H$ is known in $|\xi|<r$, equation
\eqref{eq:holonomic-germ-equation} defines it meromorphically in
$|\xi|<\rho^{-1}r$, since every argument on the right after solving for
$H(\xi)$ has smaller modulus.  On the overlap it agrees with the
original germ.  Iteration gives the unique continuation.  A zero of
$D_0$ at the origin causes no difficulty: the germ there was assumed
holomorphic and is kept throughout the construction.

A pole can arise only from a nonzero root of $D_0$ or from a pole at
one of the contracted arguments.  This proves
\eqref{eq:holonomic-pole-spirals}.  Along each chain of contracted
arguments the moduli decrease strictly.  A nonzero root of $D_0$ is
therefore encountered at most once, and contributes no more than its
multiplicity to the pole order.  Taking the largest order over the
branches gives the asserted bound, with cancellations only decreasing
the order.
\end{proof}

\begin{lemma}
\label{lem:recurrence-growth-certificate}
Let $0<\rho<1$, and suppose that, for all sufficiently large $r$,
\[
       \sum_{j=0}^d a_j(\rho^r)v_{r+j}=0,
       \qquad d\ge1,\quad a_d\ne0,
\]
where the $a_j$ are polynomials.  Discard zero coefficients, put
$\nu_j=\ord_{t=0}a_j(t)$, and define the nonnegative rational number
\begin{equation}\label{eq:recurrence-growth-slope}
 \tau=\max\left(0,\max_{0\le j<d,\ a_j\ne0}
                    \frac{\nu_d-\nu_j}{d-j}\right).
\end{equation}
An empty inner maximum is omitted.
There exist $C,R>0$ such that
\begin{equation}\label{eq:recurrence-growth-bound}
                       |v_r|\le CR^r\rho^{-\tau r^2/2}.
\end{equation}
If a recurrence has order zero, its solutions are eventually zero and
the same bound holds with $\tau=0$.
\end{lemma}

\begin{proof}
A nonzero polynomial $a_d$ has only finitely many zeros among the
points $\rho^r$.  For all large $r$,
\[
 \left|\frac{a_j(\rho^r)}{a_d(\rho^r)}\right|
                 \le C_0\rho^{(\nu_j-\nu_d)r}.
\]
Let $W_r=R^r\rho^{-\tau r^2/2}$.  Then
\[
 \left|\frac{a_j(\rho^r)}{a_d(\rho^r)}\right|
              \frac{W_{r+j}}{W_{r+d}}
 \le C_0 R^{j-d}
       \rho^{\{\nu_j-\nu_d+\tau(d-j)\}r
                              +\tau(d^2-j^2)/2}
 \le C_0R^{j-d}.
\]
The last inequality uses \eqref{eq:recurrence-growth-slope} and
$r\ge0$.  Choose $R$ so large that the sum of the last bounds over
$j<d$ is at most one, and choose $C$ to control the finitely many
initial values.  Induction using the recurrence proves
\eqref{eq:recurrence-growth-bound}.
\end{proof}

\begin{corollary}
\label{cor:universal-finite-level}
For every knot $K$, fixed rank $n\ge2$ and $q\in(0,1)$, the node order
$\delta(K,n;q)$ and Newton order $\beta(K,n;q)$ are finite.
More explicitly, take a nonzero specialized recurrence
$\sum_j a_j(q^r)J_{r+j}^{(n)}=0$, let $\tau$ be
\eqref{eq:recurrence-growth-slope}, and put
\[
                       s_0=\max\{\tau/2-1,0\}.
\]
Then $\delta(K,n;q)\le\tau/2$, $\beta(K,n;q)\le s_0$, and on every
compact set $\mathcal K\subset\mathbb C$ there exist $C_{\mathcal K}$
and $R_{\mathcal K}$ with
\begin{equation}\label{eq:universal-finite-level}
 \sup_{X\in\mathcal K}|c_{K,n,k}B_k^{(n)}(X)|
             \le C_{\mathcal K}R_{\mathcal K}^k Q^{-s_0k^2/2}.
\end{equation}
If $s_0>0$, the $Q^{s_0}$-Borel transform has positive radius of
convergence, locally uniformly in $X$, and has the meromorphic
continuation of Theorem~\ref{thm:allknots-meromorphic}.  If $s_0=0$,
the untransformed auxiliary series converges near $z=0$.
The level $s_0$ depends on the chosen recurrence and is an upper bound
for the critical level $\beta(K,n;q)$ studied in
Theorem~\ref{thm:critical-knot-borel}.
\end{corollary}

\begin{proof}
The first part of Proposition~\ref{prop:newton-holonomic}, including
its specialization argument, gives the recurrence for the node values.
Apply Lemma~\ref{lem:recurrence-growth-certificate} with $\rho=q$.
Since $Q=q^2$, its bound reads
\[
                 |J_r^{(n)}|\le CR^rQ^{-\tau r^2/4},
\]
which gives $\delta\le\tau/2$.  Inserting this estimate in
\eqref{eq:inversebound} gives
\[
 \sup_{X\in\mathcal K}|c_kB_k^{(n)}(X)|
 \le C_{\mathcal K}Q^k
       \sum_{r=0}^k R^rQ^{-(\tau/2-1)r^2/2-r/2}.
\]
If $\tau/2>1$, bound the quadratic part by its value at $r=k$;
if $\tau/2\le1$, discard its nonnegative exponent.  In both cases
the remaining exponential factors and $k+1$ are absorbed in
$R_{\mathcal K}^k$.  This proves
\eqref{eq:universal-finite-level}.  For $s_0>0$, division by
$M_k(Q^{s_0})$ leaves an exponential bound because
$(Q^{s_0};Q^{s_0})_k$ is bounded away from zero.  The remaining
conclusions follow from Lemma~\ref{lem:holonomic-meromorphy} and the
root test.
\end{proof}

\begin{theorem}
\label{thm:allknots-meromorphic}
Let $K$ be any knot, $n\ge2$, $0<q<1$ and $s\in\mathbb Q_{>0}$.
For every fixed $X\in\mathbb C$, if
\[
 H_{K,n,s}(X,\xi)=
       \sum_{k\ge0}\frac{c_{K,n,k}(q)B_k^{(n)}(X)}{M_k(Q^s)}\xi^k
\]
has positive radius of convergence, it continues meromorphically to
the whole $\xi$-plane.  Its poles have uniformly bounded orders and
belong to finitely many geometric progressions, determined by a
polynomial difference equation as in
Lemma~\ref{lem:holonomic-meromorphy}.  This applies in particular at a
critical rational level whenever the critical transform has positive
radius of convergence.  If $s>\beta(K,n;q)$, the transform is entire.
\end{theorem}

\begin{proof}
Apply Proposition~\ref{prop:newton-holonomic} and
Lemma~\ref{lem:holonomic-meromorphy}.  The final assertion follows
directly from the definition of quadratic order: dividing by
$M_k(Q^s)$ leaves Gaussian decay in $k$ when $s>\beta$.
\end{proof}

Theorem~\ref{thm:allknots-meromorphic} is pointwise in the rank and the
color. The poles it allows are those of an annihilating equation, and
some of them may be absent; at a node, for instance, the transform is a
polynomial. Proposition~\ref{prop:universal-color-transport} below gives
a set of possible poles that does not depend on the color.

\subsection{The critical Borel level of a knot}
\label{subsec:critical-knot-level}

For arbitrary complex node values, a finite quadratic order $\beta$ does
not imply that the Borel transform at base $Q^\beta$ has a positive
radius of convergence: the terms may exceed $Q^{-\beta k^2/2}$ by a
factor that grows faster than any exponential. For colored invariants
this cannot happen. They are $q$-holonomic, and Ramis's growth theorem
for entire solutions of linear $q$-difference equations shows that the
critical level is rational and attained.

\begin{lemma}
\label{lem:exact-holonomic-growth}
Let $0<\rho<1$ and let $(u_k)_{k\ge0}$ satisfy a nonzero recurrence
with coefficients polynomial in $\rho^k$.  If the sequence is not
eventually zero, there exists $\lambda\in\mathbb Q$ such that
\begin{equation}\label{eq:exact-holonomic-growth}
 \limsup_{k\to\infty}
       \frac{2\log|u_k|}{(-\log\rho)k^2}=\lambda,
 \qquad |u_k|\le CR^k\rho^{-\lambda k^2/2}
\end{equation}
for some $C,R>0$.  In particular the bound holds at the exact
quadratic exponent, and not merely at every larger exponent.
\end{lemma}

\begin{proof}
We use the following coefficient form of
\cite[Proposition~2.1 and Theorem~4.8]{RamisGrowth}.  If an entire
function $G(\xi)=\sum_{k\ge0}b_k\xi^k$ is not a polynomial and
satisfies a nonzero linear $\rho$-difference equation with polynomial
coefficients and polynomial inhomogeneous term, then there exists
$\gamma\in\mathbb Q_{>0}$ such that
\begin{equation}\label{eq:ramis-coefficient-form}
                    |b_k|\le CR^k\rho^{\gamma k^2/2},
\end{equation}
and no estimate of this form holds with any $\gamma'>\gamma$.
A polynomial inhomogeneous term can first be removed by applying
$\prod_{j=0}^{D}(T_\rho-\rho^j)$, where
$(T_\rho f)(\xi)=f(\rho\xi)$ and $D$ is its degree.
The resulting homogeneous difference operator is nonzero, since the
coefficient of its highest shift is a product of nonzero coefficients.
Recast the equation with dilation base
$\bar q=\rho^{-1}>1$ and put $h_\rho=-\log\rho$.
In explicit growth notation, Ramis's Theorem~4.8 gives a rational
$\kappa>0$ and constants $C,a$ such that, for $|\xi|$ large,
\[
 |G(\xi)|\le C\exp\left\{
       \frac{\kappa}{2h_\rho}\log^2|\xi|+a\log|\xi|\right\},
\]
with no bound of this form for a smaller positive $\kappa$.
Cauchy's estimate, optimized at
$\log r=h_\rho(k-a)/\kappa$, gives
$|b_k|\le C'R^k\rho^{k^2/(2\kappa)}$.
Conversely a bound with any $\gamma'>1/\kappa$ would, by summing
the Gaussian coefficient bound (equivalently by Proposition~2.1),
give a quadratic growth coefficient $1/\gamma'<\kappa$.
Thus \eqref{eq:ramis-coefficient-form} holds optimally with
$\gamma=1/\kappa$. The term $a\log|\xi|$ in Ramis's bound contributes
only the exponential factor $R^k$.

Lemma~\ref{lem:recurrence-growth-certificate} supplies a finite
quadratic upper bound for $u_k$.  Choose a positive integer $m$
sufficiently large that
\[
              G(\xi)=\sum_{k\ge0}u_k\rho^{mk^2}\xi^k
\]
is entire.  Multiplication of the sequence by $\rho^{mk^2}$ preserves
holonomicity, since its consecutive-term ratio is
$\rho^{m(2k+1)}$.  The generating-function calculation in
Lemma~\ref{lem:holonomic-meromorphy} gives the required equation for
$G$.  The function $G$ is not a polynomial because $u_k$ is not
eventually zero.

Apply \eqref{eq:ramis-coefficient-form} to
$b_k=u_k\rho^{mk^2}$.  It gives
\[
                    |u_k|\le CR^k\rho^{-(2m-\gamma)k^2/2}.
\]
Moreover
$\limsup 2\log|b_k|/(( -\log\rho)k^2)=-\gamma$:
a strict inequality would imply a bound with some
$\gamma'>\gamma$, contrary to optimality.  Thus
\eqref{eq:exact-holonomic-growth} holds with
$\lambda=2m-\gamma\in\mathbb Q$.  The value of $\lambda$ is fixed by
the displayed limsup and is therefore independent of the auxiliary
integer $m$.
\end{proof}

\begin{theorem}
\label{thm:critical-knot-borel}
Let $K$ be a knot, $n\ge2$ and $0<q<1$, and let
$\beta=\beta(K,n;q)$ be the quadratic order at a color that is not a
node.  Then $\beta\in\mathbb Q_{\ge0}$, and there is a constant
$R>0$ such that, for every compact color set $\mathcal K$,
\begin{equation}\label{eq:critical-uniform-bound}
 \sup_{X\in\mathcal K}|c_{K,n,k}(q)B_k^{(n)}(X)|
             \le C_{\mathcal K}R^k Q^{-\beta k^2/2}.
\end{equation}
Consequently:
\begin{enumerate}[label=(\roman*),leftmargin=2em]
\item If $\beta=0$, the auxiliary Newton series converges normally
      for $X$ in compact sets and $z$ sufficiently small.
\item If $\beta>0$, its critical $Q^\beta$-Borel transform has a
      positive finite radius $R_{K,n,q}$, independent of every color
      that is not a node.  It continues meromorphically to the plane,
      with poles of bounded order on finitely many geometric
      progressions.  At a node the transform is a polynomial.
\end{enumerate}
\end{theorem}

\begin{proof}
The color $X=0$ is not a node, since all nodes are positive and larger
than two.  Put $u_k=c_{K,n,k}(q)B_k^{(n)}(0)$.  By
Proposition~\ref{prop:newton-holonomic}, this is a $q$-holonomic
numerical sequence after specialization.  If $u_k$ is eventually
zero, the Newton coefficients are eventually zero and all assertions
in the zero-order case are immediate.  Otherwise apply
Lemma~\ref{lem:exact-holonomic-growth} with $\rho=q$.
Since $h_Q=2(-\log q)$, its exponent is related to the order used in
this paper by
\[
                   \beta=\frac{\max\{\lambda,0\}}2.
\]
This proves rationality and the critical bound at $X=0$.

Write $A=q^n$ and
$P_k(A,X)=\prod_{j<k}(1-AXQ^j+A^2Q^{2j})$.
The factorization of the basis gives
\[
       \frac{B_k^{(n)}(X)}{B_k^{(n)}(0)}
               =\frac{P_k(A,X)}{P_k(A,0)}.
\]
The numerator is bounded uniformly in $k$ on every compact color
set, while $P_k(A,0)\ge1$.  Hence the bound at $X=0$ gives
\eqref{eq:critical-uniform-bound}.  If $\beta=0$, normal convergence
near $z=0$ follows.  If $\beta>0$, dividing by $M_k(Q^\beta)$ leaves
an exponential bound, so the critical Borel transform has positive
radius, locally uniformly in $X$.

Let $H_0$ be this transform at $X=0$.  It cannot be entire.  Otherwise,
writing $\beta=a/b$, it would be a nonpolynomial entire function
satisfying a $\rho$-difference equation with $\rho=q^{1/b}$, by
Proposition~\ref{prop:newton-holonomic}, and formula
\eqref{eq:ramis-coefficient-form} would give
\[
 |[\xi^k]H_0|\le C R_1^k\rho^{\gamma k^2/2}
                    \qquad(\gamma>0).
\]
Multiplication by $M_k(Q^\beta)$ would then bound $u_k$ at a quadratic
order strictly less than $\beta$, a contradiction.  Thus the radius
$R_{K,n,q}$ is finite.  Theorem~\ref{thm:allknots-meromorphic} gives
meromorphic continuation.

For every color that is not a node,
$P_k(A,X)/P_k(A,0)\to E_A(X)/E_A(0)\ne0$.
The root test therefore gives exactly the same critical radius as at
$X=0$.  At a node the Newton expansion terminates.
\end{proof}

Theorem~\ref{thm:critical-knot-borel} concerns the Borel transform only.
The critical Laplace integral converges once the Borel transform grows
slowly enough on the integration ray, for instance under condition
\eqref{eq:holonomic-laplace-test} of
Proposition~\ref{prop:holonomic-growth-test}, and by
Corollary~\ref{cor:one-color-growth-test} it suffices to check such a
bound at the single color $X=0$.

\begin{proposition}
\label{prop:first-singularities-node-factor}
Assume $\beta(K,n;q)>0$, and write $H(X,\xi)$ for the critical Borel
transform and $R_*=R_{K,n,q}$.  Then the germ
\begin{equation}\label{eq:first-shell-comparison}
 H(X,\xi)-\frac{E_n(X)}{E_n(0)}H(0,\xi)
\end{equation}
is holomorphic on $\mathbb C_X\times\{|\xi|<R_*/Q\}$.
In particular, if $\alpha$ is a pole of $H(0,\xi)$ with
$|\alpha|=R_*$, then
\begin{equation}\label{eq:first-pole-principal-part}
 \operatorname{PP}_{\xi=\alpha}H(X,\xi)
       =\frac{E_n(X)}{E_n(0)}
                     \operatorname{PP}_{\xi=\alpha}H(0,\xi).
\end{equation}
There is at least one such pole.  Each of these nearest poles is
present, with the same order, at every color that is not a node, and
disappears exactly at the nodes.
\end{proposition}

\begin{proof}
The infinite-product tail gives, uniformly for $X$ in compact sets,
\[
 \frac{P_k(A,X)}{P_k(A,0)}
          -\frac{E_A(X)}{E_A(0)}=O(Q^k).
\]
Indeed $P_k(A,X)-E_A(X)=O(Q^k)$ uniformly on such sets, and the same
holds at $X=0$, where both the finite products and their limit are
bounded away from zero.  Only $P_k(A,0)$ and $E_A(0)$ occur in
denominators, so the estimate also holds at the nodes.

The $k$th coefficient of \eqref{eq:first-shell-comparison} is therefore
$O_{\mathcal K}(Q^k[\xi^k]H(0,\xi))$.
Its series converges normally for $|\xi|<R_*/Q$ by the root test,
which proves the first assertion.  Since $H(0,\xi)$ is meromorphic and
has finite Taylor radius $R_*$, a singularity on $|\xi|=R_*$ is a
pole.  Subtracting the holomorphic comparison term gives
\eqref{eq:first-pole-principal-part}.  Finally $E_n$ vanishes exactly
at the nodes, and a nonzero scalar multiple of a principal part has
the same pole order.
\end{proof}

\subsection{Color transport from one reference color}
\label{subsec:general-color-transport}

The factorization of the Newton basis transports the continued Borel
transform from the color $X=0$ to every color, with a set of possible
poles that does not depend on the color. Write
$P_k(A,X)=\prod_{j<k}(1-AXQ^j+A^2Q^{2j})$ as above.

\begin{proposition}
\label{prop:universal-color-transport}
Fix $A=q^n$, and a Borel level for which
\[
 H(X,\xi)=\sum_{k\ge0}b_k
              \frac{P_k(A,X)}{P_k(A,0)}\xi^k
\]
converges near $\xi=0$.  Suppose that $H_0(\xi)=H(0,\xi)$ has a
meromorphic continuation to $\mathbb C$, with pole set $\Sigma_0$.
Then $H$ has a meromorphic continuation, entire in the color away from
a locally finite common possible pole set
\[
             \Sigma_* = \bigcup_{m\ge0}Q^{-m}\Sigma_0.
\]
More precisely it is jointly holomorphic on
$\mathbb C_X\times(\mathbb C_\xi\setminus\Sigma_*)$.
If $\Sigma_0$ is contained in finitely many progressions
$\alpha\rho^{-\mathbb N_0}$ and $Q=\rho^d$ for a positive integer $d$,
the same is true of $\Sigma_*$.  A uniform bound on the pole orders of $H_0$ is also a
bound for the pole orders of $H(X,\cdot)$ at every color.

There are polynomials $v_m(X)$ of degree at most $m$, depending only
on $A,Q$, such that every principal part coefficient has the form
\begin{equation}\label{eq:general-color-principal-parts}
 [ (\xi-\alpha)^{-j}]H(X,\xi)
  =\frac{E_A(X)}{E_A(0)}
    \sum_{\substack{m\ge0\\Q^m\alpha\in\Sigma_0}}
       v_m(X)Q^{-mj}
       [ (\xi-Q^m\alpha)^{-j}]H_0(\xi).
\end{equation}
The sum is finite.  Consequently each coefficient of every principal
part is $E_A(X)$ times a polynomial in $X$. Away from the nearest
poles, this polynomial may also vanish at colors that are not nodes.
\end{proposition}

\begin{proof}
For small $t$, set
\[
 V_X(t)=\prod_{j\ge0}
       \frac{1+A^2t^2Q^{2j}}
            {1-AXtQ^j+A^2t^2Q^{2j}}
       =\sum_{m\ge0}v_m(X)t^m.
\]
This is analytic near $t=0$, uniformly for $X$ in compact sets.  Its
functional equation is
\[
 (1-AXt+A^2t^2)V_X(t)=(1+A^2t^2)V_X(Qt),
\]
so $v_0=1$, $v_{-1}=0$, and, for $m\ge1$,
\[
 (1-Q^m)v_m
     =AXv_{m-1}-A^2(1-Q^{m-2})v_{m-2}.
\]
This proves the polynomial assertion.  On a given compact color set
$\mathcal K$, choose $r>0$ so that $V_X(t)$ is analytic for
$|t|\le r$, uniformly in $X\in\mathcal K$, and then choose $N\ge1$
with $Q^N<r$.  For every $k\ge N$, the product identity gives
\[
 \frac{P_k(A,X)}{P_k(A,0)}
       =\frac{E_A(X)}{E_A(0)}V_X(Q^k).
\]
This identity is valid at nodes too: the denominator in $V_X(Q^k)$
does not vanish for the chosen $k$ and compact set.

Let $H_{<N}(X,\xi)$ denote the first $N$ Taylor terms of $H$ and put
$H_{0,<N}=H_{<N}(0,\cdot)$.  Interchanging normally convergent series
near zero gives the exact continuation formula
\begin{equation}\label{eq:general-color-transport}
\begin{split}
 H(X,\xi)={}&H_{<N}(X,\xi)\\
 &+\frac{E_A(X)}{E_A(0)}
       \sum_{m\ge0}v_m(X)
       \bigl[H_0(Q^m\xi)-H_{0,<N}(Q^m\xi)\bigr].
\end{split}
\end{equation}
Indeed Cauchy's bound gives $|v_m(X)|\le C_{\mathcal K}r^{-m}$,
while the bracket is $O(Q^{mN})$, uniformly on bounded $\xi$ sets for
large $m$.  Thus the sum converges normally on compact subsets of
$\mathcal K\times(\mathbb C\setminus\Sigma_*)$.
The set $\Sigma_0$ is bounded away from zero, so only finitely many
$m$ contribute poles in any bounded disk.  This proves local
finiteness, joint continuation, and the assertion about pole orders.
The continuations agree as $\mathcal K,N$ vary by uniqueness from the
original germ.  Taking principal parts in the finite set of singular
terms yields \eqref{eq:general-color-principal-parts}; the factor
$Q^{-mj}$ comes from the change of variable $\xi\mapsto Q^m\xi$.
\end{proof}

\begin{corollary}
\label{cor:one-color-growth-test}
In Proposition~\ref{prop:universal-color-transport}, assume that
$\Sigma_0$ lies on finitely many progressions $\alpha\rho^{-\mathbb N_0}$,
that $Q=\rho^d$ with $d\ge1$ integral, and that the pole orders are
bounded.
Suppose on a sector, or outside fixed relative disks about $\Sigma_*$,
that
\[
 |H_0(\xi)|\le C(1+|\xi|)^a
             \exp\{\kappa\log^2(2+|\xi|)\},\qquad \kappa\ge0.
\]
Then the same quadratic growth exponent $\kappa$ is valid for
$H(X,\xi)$, uniformly for $X$ in any compact set, after changing the
constant and the polynomial factor.  In particular a subcritical
bound $\kappa<1/(2h_p)$ at $X=0$ gives ordinary $p$-Laplace sums on
the corresponding pole-free sector at every color, locally uniformly
in $X$. If the bound holds in a wedge about a singular ray outside
fixed relative pole disks, it also gives both lateral sums and the
Stokes jump on that ray.
\end{corollary}

\begin{proof}
The domain under consideration is preserved by multiplication by
$Q^m$.  The Taylor remainder
$G_N=H_0-H_{0,<N}$ obeys, on this domain,
\[
 |G_N(\xi)|\le C_N|\xi|^N(1+|\xi|)^{a_N}
               \exp\{\kappa\log^2(2+|\xi|)\}
\]
for some $a_N\ge0$: near zero use the Taylor remainder, and away from
zero absorb the finite polynomial and the factor $|\xi|^{-N}$.
In \eqref{eq:general-color-transport}, the $m$th remainder is
therefore bounded by the displayed right side times $Q^{mN}$.
The resulting geometric series converges because $Q^N<r$.
The subcritical bound gives the Laplace integral on each such
sector. In a wedge about a singular ray it also controls the fixed
relative detours and their residue contributions, as in
Proposition~\ref{prop:holonomic-growth-test}.
\end{proof}

For a knot at its critical rational level, the hypotheses on the poles
of $H_0$ and on their orders follow from
Theorem~\ref{thm:critical-knot-borel}; the growth bound has to be
established separately.

\subsection{Growth on Laplace contours}

\begin{proposition}
\label{prop:holonomic-growth-test}
In Lemma~\ref{lem:holonomic-meromorphy}, omit the zero $D_\ell$ and put
\[
 h_\rho=-\log\rho,\qquad
 \mu=\max\left(0,\max_{\ell\ge1}
              \frac{\deg D_\ell-\deg D_0}{\ell}\right),
\]
where the inner maximum is omitted if there are no positive shifts.
For every $\kappa>\mu/(2h_\rho)$ there is a constant $C_\kappa$ such
that
\begin{equation}\label{eq:holonomic-global-growth}
 |H(\xi)|\le C_\kappa
      \exp\{\kappa\log^2(2+|\xi|)\}
\end{equation}
outside sufficiently small disks of fixed relative radius around the
possible poles in \eqref{eq:holonomic-pole-spirals}.
If $\mu=0$, one may instead take a polynomial bound there.
Consequently the $p$-Laplace integral converges on every pole-free
sector whenever
\begin{equation}\label{eq:holonomic-laplace-test}
                 \mu<\frac{h_\rho}{h_p},\qquad h_p=-\log p.
\end{equation}
Under this strict inequality the lateral sums also exist.
Their difference is the convergent sum of the residues of
$H(\xi)w_p(\xi/z)$, multiplied by $2\pi i/(h_pz)$, with the
orientation specified in Lemma~\ref{merged:subcritical-laterals}.
\end{proposition}

\begin{proof}
After solving \eqref{eq:holonomic-germ-equation} for $H$, the rational
coefficient of $H(\rho^\ell\xi)$ is bounded for large $|\xi|$ by
$C|\xi|^{d_\ell}$, where $d_\ell=\deg D_\ell-\deg D_0$.
For $W(r)=\exp(\kappa\log^2r)$,
\[
 r^{d_\ell}\frac{W(\rho^\ell r)}{W(r)}
   =\exp(\kappa\ell^2h_\rho^2)
      r^{d_\ell-2\kappa\ell h_\rho}\longrightarrow0.
\]
The sum of these ratios is less than $1/2$ for large $r$.  The
inhomogeneous term is rational and has at most polynomial growth, so
induction on annuli gives \eqref{eq:holonomic-global-growth}.  The
induction stays outside the excluded disks: if $\rho^\ell\xi$ is in
the disk about a possible pole $\alpha\rho^{-m}$, then $\xi$ is in
the disk with the same relative radius about
$\alpha\rho^{-m-\ell}$.  A finite collection of orbits has positive
relative separation after coincident poles are identified, so these
disks can be chosen disjoint.  The bounded starting annulus supplies
the induction constant.

If $\mu=0$, all $d_\ell\le0$.  Replacing $W(r)$ by $(1+r)^L$ and
taking $L$ sufficiently large makes the finite sum of coefficient
bounds times $\rho^{\ell L}$ smaller than $1/2$; take $L$ also larger
than the degree of the inhomogeneous term.  This proves polynomial
growth.  If there are no positive shifts, $H$ is rational and both
conclusions are immediate.

Under \eqref{eq:holonomic-laplace-test}, choose $\kappa$ strictly
between $\mu/(2h_\rho)$ and $1/(2h_p)$.  The Gaussian decay in
$\log t$ of $w_p(t)$ then dominates the growth of $H(zt)$.
Lemma~\ref{merged:subcritical-laterals} gives the lateral limits and
the residue formula under the same bound.
\end{proof}

\begin{remark}\label{rem:growth}
The growth condition in \eqref{eq:holonomic-laplace-test} cannot be
dropped. The function
\[
 H(\xi)=\sum_{k\ge0}Q^{k^2}\xi^k,
 \qquad H(\xi)=1+Q\xi H(Q^2\xi),
\]
is entire and $q$-holonomic, but its $p$-Laplace integral with
$p=Q^3$ diverges for every positive $z$. Indeed, all terms are
nonnegative on the positive axis, and Tonelli's theorem gives the
integral as the divergent series
$\sum_{k\ge0}Q^{k^2}M_k(Q^3)z^k$.
\end{remark}

\subsection{Connected sums}
\label{subsec:connected}
For knots we use reduced, zero-framed quantum invariants colored by an
irreducible representation. Cutting a knot open gives a
$(1,1)$-tangle, whose Reshetikhin--Turaev operator is an equivariant
endomorphism of the irreducible representation, hence a scalar by
Schur's lemma, and this scalar is the reduced invariant. Composing the
tangles of $K$ and $L$ gives a tangle for $K\#L$ and multiplies the
scalars \cite{RT}. Hence
\begin{equation}\label{eq:connected-sum-values}
 J_r^{(n)}(K\mathbin{\#}L;q)
       =J_r^{(n)}(K;q)J_r^{(n)}(L;q).
\end{equation}
A Newton series is determined by its values at the nodes: if all
values vanish, evaluating successively at $x_0,x_1,\ldots$ and using
$B_k(x_k)\ne0$ shows that every coefficient vanishes. The product
$F_{K,n}F_{L,n}$ has the values on the right of
\eqref{eq:connected-sum-values}, so
\begin{equation}\label{eq:connected-sum-Newton}
                 F_{K\# L,n}=F_{K,n}F_{L,n}.
\end{equation}
The same argument works over the integral domain
\(\mathbb Z[q^{\pm1}]\) before specialization.

Write $\delta_K$ for the quadratic order of the reduced, zero-framed
invariants of $K$ at the nodes, at fixed rank and $q$, and $\beta_K$
for the quadratic order of its auxiliary series at a color that is not
a node.

\begin{corollary}\label{cor:general-connected}
For every knot $K$, every $m\ge1$, and every fixed rank,
\begin{equation}\label{eq:general-connected}
 \beta_{\#^mK}=\max\{m\delta_K-1,0\}.
\end{equation}
If $\beta_K>0$, this reads $\beta_{\#^mK}=m\beta_K+m-1$. For knots
$K_1,\ldots,K_m$ one always has
\[
 \beta_{K_1\#\cdots\#K_m}\le
          \max\{\delta_{K_1}+\cdots+\delta_{K_m}-1,0\},
\]
with equality when the limits
$\lim_{r\to\infty}2\log|J_r(K_i)|/(h_Qr^2)$ exist and are finite and
nonnegative.
\end{corollary}
\begin{proof}
Reduced invariants multiply under connected sum. For a power,
$\log|J_r(K)^m|=m\log|J_r(K)|$, also where $J_r(K)=0$, so the node order
of $\#^mK$ is exactly $m\delta_K$, and Theorem~\ref{thm:nodegrowth}
applies. For different knots, the $\limsup$ of a sum is at most the sum
of the $\limsup$s, and it is equal to it when the limits exist. Without
this hypothesis the maximal growth of the factors may occur along
different subsequences, and the bound need not be attained.
\end{proof}

\section{Spectral summation and Stokes data}
\label{sec:spectral}\label{sec:stokes}

This section treats Borel transforms that are convergent sums of
simple poles along a geometric progression. For them the Borel
singularities, the Laplace sums and the Stokes jumps are explicit, and
the jump determines the transform (Theorem~\ref{gen:thm:injective}).
After finitely many initial terms are removed, the Borel transform of
every Newton series whose coefficients have an analytic tail is a finite
combination of such sums and of their Euler derivatives
(Theorem~\ref{thm:analytic-symbol}). At
the end of the section we show that this calculus is compatible with
the divided differences of the Newton basis, which change the rank.

\subsection{The resolvent}\label{subsec:resolvent}
For $R>1$ let
\[
 \W_R=\Big\{b(t)=\sum_{m\ge0}b_mt^m:\
 \|b\|_R=\sum_{m\ge0}|b_m|R^m<\infty\Big\}.
\]
With pointwise multiplication this is a unital Banach algebra, and its
elements are holomorphic on the disk $|t|<R$. For $\alpha\ne0$ and
$b\in\W_R$ put
\begin{equation}\label{gen:eq:resolvent}
 (\sRes_\alpha b)(\xi)=\sum_{m\ge0}\frac{b_m}{1-Q^m\xi/\alpha},
 \qquad \alpha_m=\alpha Q^{-m}.
\end{equation}
We call $b$ the \emph{spectral symbol} of $\sRes_\alpha b$ and $\alpha$
its \emph{pole seed}. The operator $\T f(\xi)=f(Q\xi)$ shifts the index
$m$ in \eqref{gen:eq:resolvent} by one, and this gives the
multiplication rule \eqref{gen:eq:spectralproduct}.

\begin{theorem}\label{gen:thm:resolvent}
Let $b,g\in\W_R$.
\begin{enumerate}[label=(\roman*),leftmargin=2em]
\item $\sRes_\alpha b$ is meromorphic on $\C$. Its poles are simple and
lie among the points $\alpha_m$, and $\alpha_m$ is a pole exactly when
$b_m\ne0$. Moreover
\[
 \Res_{\xi=\alpha_m}\sRes_\alpha b=-\alpha_m b_m,
 \qquad
 [\xi^k]\sRes_\alpha b=\alpha^{-k}b(Q^k).
\]
\item On each closed sector at positive angular distance $\theta$ from
the ray $\alpha\mathbb R_{>0}$, $|\sRes_\alpha b|\le C_\theta\|b\|_R$,
uniformly in $|\xi|$.
\item The series $g(\T)=\sum_{\ell\ge0}g_\ell\T^\ell$ converges normally
on such functions, and
\begin{equation}\label{gen:eq:spectralproduct}
 g(\T)\sRes_\alpha b=\sRes_\alpha(gb).
\end{equation}
If $R_m$ is the residue of $\sRes_\alpha b$ at $\alpha_m$, the residues
of $g(\T)\sRes_\alpha b$ are
\begin{equation}\label{gen:eq:resconv}
 R'_m=\sum_{\ell=0}^m g_\ell Q^{-\ell}R_{m-\ell}.
\end{equation}
\end{enumerate}
\end{theorem}
\begin{proof}
On a compact set that avoids the points $\alpha_m$, the denominators in
\eqref{gen:eq:resolvent} are bounded away from zero and tend uniformly
to one as $m\to\infty$. Since $\sum|b_m|<\infty$, the series converges
normally there, and only the term with index $m$ is singular at
$\alpha_m$, where it has residue $-\alpha_mb_m$. For $|\xi|<|\alpha|$,
expanding each term in a geometric series and summing over $m$ first
gives the Taylor coefficients. If $w=Q^m\xi/\alpha$ satisfies
$|\arg w|\ge\theta$, then $|1-w|\ge c_\theta(1+|w|)$ with $c_\theta>0$,
which gives (ii). Finally, $\T^\ell$ replaces the term with index $m$
by the term with index $m+\ell$, so the coefficient of
$(1-Q^j\xi/\alpha)^{-1}$ in $g(\T)\sRes_\alpha b$ is
$\sum_{\ell\le j}g_\ell b_{j-\ell}$, the $j$-th Taylor coefficient of
$gb$. Absolute convergence in $\W_R$ justifies the rearrangement and
proves \eqref{gen:eq:spectralproduct}. Since $R_m=-\alpha_mb_m$ and
$\alpha_m/\alpha_{m-\ell}=Q^{-\ell}$, this gives \eqref{gen:eq:resconv}.
\end{proof}

\subsection{Laplace sums and lateral limits}
Fix $0<p<1$, independently of $Q$, and let $\cL_p$ be the Laplace
transform of \eqref{eq:KMM}. By Theorem~\ref{thm:kernel} the normalized
moments of $w_p$ are $M_k(p)$, and $w_p$ decays faster than every power
of $t$.

\begin{theorem}\label{gen:thm:sumjump}
Let $H=\sRes_\alpha b$ with $b\in\W_R$.
\begin{enumerate}[label=(\roman*),leftmargin=2em]
\item The integral $\cL_pH(z)$ converges for $z\ne0$ off the ray
$\alpha\mathbb R_{>0}$ and is holomorphic there. On each closed sector
at angular distance $\theta>0$ from this ray, and for every $N\ge0$,
\begin{equation}\label{gen:eq:remainder}
 \left|\cL_pH(z)-\sum_{k=0}^{N-1}\alpha^{-k}b(Q^k)M_k(p)z^k\right|
 \le C_\theta\|b\|_R M_N(p)|z/\alpha|^N.
\end{equation}
\item $\cL_pH$ has boundary values from both sides of the ray, locally
uniformly for $z=e^{i\arg\alpha}r$ with $r>0$. Let $\Delta_b(z)$ be the
boundary value from the side of decreasing argument minus the boundary
value from the side of increasing argument; on the negative axis this
is the upper value minus the lower value. With
$\tau_m=|\alpha|Q^{-m}/r$,
\begin{align}
 \Delta_b(z)
 &=\frac{2\pi i}{h_pz}\sum_{m\ge0}
 \Res_{\xi=\alpha_m}H(\xi)\,w_p(\alpha_m/z) \label{gen:eq:jumpres}\\
 &=-\frac{2\pi i}{h_p}\sum_{m\ge0}b_m\tau_m w_p(\tau_m),
 \label{gen:eq:jumpcoeff}
\end{align}
and the series converge normally for $r$ in compact subsets of
$(0,\infty)$.
\item Let $\T$ act on functions of $z$ by $z\mapsto Qz$. For every
$g\in\W_R$, $\cL_p\sRes_\alpha(gb)=g(\T)\cL_p\sRes_\alpha b$, and the
same holds for both lateral sums and for the jump:
$\Delta_{gb}=g(\T)\Delta_b$.
\end{enumerate}
\end{theorem}

In general $\cL_pH$ does not extend holomorphically to $z=0$; by (i),
its asymptotic expansion there is $\sum_k\alpha^{-k}b(Q^k)M_k(p)z^k$.

\begin{proof}
By Theorem~\ref{gen:thm:resolvent}(ii), $|H(zt)|\le C_\theta\|b\|_R$ for
$z$ in the closed sector and $t>0$, so the integral converges and is
holomorphic. For the remainder we use
\[
 \frac1{1-w}-\sum_{k<N}w^k=\frac{w^N}{1-w},\qquad w=Q^m\xi/\alpha,
\]
whose modulus is at most $C_\theta|\xi/\alpha|^NQ^{mN}$ in the sector.
Summing over $m$ bounds the Taylor remainder of order $N$ of $H$ by
$C_\theta\|b\|_R|\xi/\alpha|^N$, and integration against $w_p$ gives
\eqref{gen:eq:remainder}, because the $N$-th normalized moment is
$M_N(p)$.

For the lateral limits, the term of index $m$ in $\cL_pH(z)$ is
\[
 \frac{b_m}{h_p}\int_0^\infty\frac{w_p(t)}{1-t/\tau}\,dt,
 \qquad \tau=\frac{\alpha_m}{z}.
\]
As $z$ tends to the point $e^{i\arg\alpha}r$ of the ray, $\tau$ tends to
$\tau_m>0$, from above when $z$ comes from the side of decreasing
argument and from below otherwise. Lemma~\ref{lem:lateral-tail} with $K=0$ shows that
these integrals are bounded uniformly in $m$ and near the ray, and that
their boundary values from above and from below differ by
$-2\pi i\tau_mw_p(\tau_m)$. For $r$ in a compact subset of $(0,\infty)$
the $\tau_m$ are bounded below by a positive constant, so
$\sum|b_m|<\infty$ allows us to pass to the boundary values term by
term. This gives \eqref{gen:eq:jumpcoeff}, and \eqref{gen:eq:jumpres}
follows from $\Res_{\alpha_m}H=-\alpha_mb_m$ and $\alpha_m/z=\tau_m$.

The bounds in (i) and (ii) are linear in $\|b\|_R$, so $\cL_p$, the
lateral limits and $\Delta$ are continuous linear maps on $\W_R$. Since
$\cL_p(\T H)=\T\cL_pH$ and $\T$ preserves both sides of the ray,
applying these maps to \eqref{gen:eq:spectralproduct} proves (iii).
\end{proof}

\begin{remark}\label{gen:rem:tails}
Theorem~\ref{gen:thm:sumjump} extends to functions
$P(\xi)+(\xi/\alpha)^N\sRes_\alpha b$ with a polynomial $P$ and a fixed
$N\ge0$. The factor $\xi^N$ adds only polynomial growth, the kernel has
moments of all orders, and in the lateral argument $t^Nw_p(t)$ replaces
$w_p(t)$, which is the case $K=N$ of Lemma~\ref{lem:lateral-tail}. The
polynomial changes the Laplace sum but has no jump. Since
$\T^\ell(\xi^Nf)=Q^{N\ell}\xi^N\T^\ell f$, the operator $g(\T)$ replaces
the symbol $b(t)$ of the tail by $g(Q^Nt)b(t)$. This form arises when
the coefficients are $\alpha^{-k}g(Q^k)$ for a function $g$ that is
holomorphic only near $t=0$: the coefficients with $k\ge N$ are those of
$(\xi/\alpha)^N\sRes_\alpha b_N$ with $b_N(t)=g(Q^Nt)$, and $b_N$ is
holomorphic on a disk of radius greater than one once $N$ is large.
\end{remark}

Polynomial factors in the index $k$ correspond to Euler derivatives.

\begin{proposition}\label{gen:prop:jets}
Let $\mathsf E=\xi\partial_\xi$, let $b_0,\ldots,b_d\in\W_R$, and put
\[
 H=\sum_{j=0}^d\mathsf E^j\sRes_\alpha b_j.
\]
Then $[\xi^k]H=\alpha^{-k}\sum_{j=0}^dk^jb_j(Q^k)$. The function $H$ has
poles of order at most $d+1$ at the points $\alpha_m$, it is bounded on
closed sectors away from the ray $\alpha\mathbb R_{>0}$, and its Laplace
sums and lateral sums exist as in Theorem~\ref{gen:thm:sumjump}. The
operator $g(\T)$ multiplies each $b_j$ by $g$, and
\[
 \cL_p(\mathsf E H)=(z\partial_z)\cL_pH,\qquad
 \Delta_{\mathsf E H}=(z\partial_z)\Delta_H.
\]
If $h_{m,j}$ is the coefficient of $(\xi-\alpha_m)^{-j}$ in the principal
part of $H$ at $\alpha_m$, then
\[
 \Delta_H(z)=\frac{2\pi i}{h_pz}
 \sum_{m\ge0}\sum_{j=1}^{d+1}
 \frac{h_{m,j}}{(j-1)!\,z^{j-1}}w_p^{(j-1)}(\alpha_m/z),
\]
with normal convergence for $r$ in compact subsets of $(0,\infty)$.
\end{proposition}
\begin{proof}
With $w=Q^m\xi/\alpha$ we have $\mathsf E=w\,d/dw$, so
$\mathsf E^j(1-w)^{-1}$ is a linear combination of
$(1-w)^{-1},\ldots,(1-w)^{-j-1}$ with coefficients independent of $m$.
These functions are bounded on closed sectors away from $w\in(0,\infty)$,
so the bounds of Theorem~\ref{gen:thm:resolvent} hold for $H$, and
differentiation under the integral sign gives the formula for
$\cL_p(\mathsf EH)$. After a rotation of $z$ that moves the ray
$\alpha\mathbb R_{>0}$ to the negative axis, each $\cL_p\sRes_\alpha b_j$
is a sum of the kind treated in Proposition~\ref{prop:lateral-jets}.
That proposition gives the boundary values of all fixed Euler
derivatives, locally uniformly in $r$, and the identity
$\Delta_{\mathsf EH}=(z\partial_z)\Delta_H$, because $z\partial_z=r\partial_r$
on the ray. The jump of $\int_0^\infty f(t)(t-\tau)^{-j}\,dt$ across
$\tau>0$ is $2\pi if^{(j-1)}(\tau)/(j-1)!$, and applied to the principal
parts this gives the formula for $\Delta_H$. Finally, $\mathsf E$
commutes with $\T$, so $g(\T)$ acts on each $b_j$ separately.
\end{proof}

\subsection{The jump determines the symbol}\label{subsec:injective}

\begin{theorem}\label{gen:thm:injective}
For fixed $p$, $Q$ and $\alpha$, the map $b\mapsto\Delta_b$ is injective
on $\W_R$. More precisely, if $m_0$ is the least index with
$b_{m_0}\ne0$, then
\begin{equation}\label{gen:eq:firstpoleasym}
 \Delta_b(e^{i\arg\alpha}r)
 \sim-\frac{2\pi i}{h_p}b_{m_0}\tau_{m_0}w_p(\tau_{m_0})
 \qquad(r\downarrow0),
\end{equation}
so the coefficients $b_m$ can be read off from $\Delta_b$ one at a time.
If $g,b\in\W_R$ are nonzero, then $\Delta_{gb}\ne0$. If $g$ is a unit
of $\W_R$, then $g(\T)$ is invertible, with inverse $g^{-1}(\T)$, on
Borel transforms, Laplace sums and jumps.
\end{theorem}
\begin{proof}
Put $f(\tau)=\tau w_p(\tau)$. Then
\[
 \frac{d\log f(\tau)}{d\log\tau}
 =1-\sum_{j\ge0}\frac{\tau p^j}{1+\tau p^j},
\]
which tends to $-\infty$ as $\tau\to\infty$. Hence $f$ is eventually
decreasing and $f(c\tau)/f(\tau)\to0$ for every fixed $c>1$. As
$r\downarrow0$ we have $\tau_{m_0}\to\infty$ and
$\tau_m=Q^{-(m-m_0)}\tau_{m_0}$. After division of
\eqref{gen:eq:jumpcoeff} by $f(\tau_{m_0})$, each term with $m>m_0$
therefore tends to zero and is eventually bounded by $|b_m|$, and
dominated convergence proves \eqref{gen:eq:firstpoleasym}, whose right
side is not zero. Subtracting the term with index $m_0$ and repeating
the argument recovers the remaining coefficients. Elements of $\W_R$ are
holomorphic on a disk, so $gb\ne0$ when $g$ and $b$ are nonzero. The
last assertion follows from Theorem~\ref{gen:thm:sumjump}(iii).
\end{proof}

Injectivity holds within the class of transforms $\sRes_\alpha b$:
adding a polynomial to the Borel transform changes the Laplace sum but
not the jump.

\subsection{Coefficients with analytic tails}
\label{subsec:analytic-tails}

The Newton coefficients of the twist knots in
Section~\ref{pa:data} have the following form.

\begin{definition}
At a fixed real \(0<q<1\), a Newton coefficient sequence has an
\emph{analytic tail} if there are finitely many numbers
\(\lambda_\nu\ne0\), polynomials \(P_\nu\), functions \(f_\nu(s)\)
holomorphic near \(s=0\), and an integer \(K\), such that
\[
 c_k(q)=\sum_\nu\lambda_\nu^kP_\nu(k)f_\nu(Q^k)
 \qquad(k\ge K).
\]
The integer $K$ may always be increased, so we may assume that all
points $Q^k$, $k\ge K$, lie in the disks where the $f_\nu$ are
holomorphic.
\end{definition}

\begin{lemma}\label{lem:remainder}
Suppose $g$ is holomorphic near zero and on an open sector, and is
bounded by a polynomial on each closed subsector. Write
$g(\xi)=\sum b_k\xi^k$ near zero. On every closed subsector there are
constants $C,D,N_0$ such that for $N\ge N_0$ and all $\xi$ in the
subsector
\[
 \left|g(\xi)-\sum_{k<N}b_k\xi^k\right|\le CD^N|\xi|^N.
\]
Consequently $\cL_pg$ has the asymptotic expansion
$\sum b_kM_k(p)z^k$: for $z$ in the same closed subsector and
$N\ge N_0$,
\begin{equation}\label{eq:general-remainder}
 \left|\cL_pg(z)-\sum_{k<N}b_kM_k(p)z^k\right|
 \le CD^NM_N(p)|z|^N.
\end{equation}
\end{lemma}

\begin{proof}
Choose $0<r<R$ such that $g$ is holomorphic on $|\xi|\le R$. Cauchy's
estimates give $|b_k|\le C_1R^{-k}$, and for $|\xi|\le r$ the Taylor
remainder is at most $C_2R^{-N}|\xi|^N$. For $|\xi|\ge r$ in the
subsector we bound $g$ and the first $N$ Taylor terms separately and
divide by $|\xi|^N$. If $|g(\xi)|\le C_3(1+|\xi|)^e$, then for $N\ge e$
\[
 \frac{|g(\xi)|}{|\xi|^N}\le C_3(1+r)^er^{-N},
 \qquad
 \sum_{k<N}\frac{|b_k\xi^k|}{|\xi|^N}
 \le C_1r^{-N}\sum_{k<N}(r/R)^k\le\frac{C_1r^{-N}}{1-r/R}.
\]
Hence the first inequality holds with $D=1/r$. For $z$ in the
subsector, the points $zt$, $t>0$, lie in it as well. Applying the
inequality with $\xi=zt$, integrating against $w_p(t)\,dt/h_p$, and
using \eqref{eq:moment-identity} gives \eqref{eq:general-remainder}.
\end{proof}

\begin{theorem}\label{thm:analytic-symbol}
Suppose the coefficients of a Newton element have an analytic tail.
Then the \(Q\)-Borel transform of its auxiliary series continues
meromorphically to the whole Borel plane, with poles only at
\[
 \xi=-\frac{A}{Q\lambda_\nu}Q^{-m},\qquad m\ge0.
\]
A pole coming from the $\nu$th term has order at most
\(1+\deg P_\nu\); where the progressions of several terms meet, the
order is at most the largest of their bounds. The Borel transform grows
at most polynomially on closed sectors that avoid the finitely many
rays through these poles, locally uniformly in \(X\). Its Laplace sum
$S$ is an entire function of $X$, it equals the terminating Newton
polynomial at every node, and both lateral sums exist on each singular
ray, locally uniformly in the color and in positive modulus. Their
difference is the normally convergent residue sum
\[
 \Delta(X,z)=\frac{2\pi i}{h_Qz}
 \sum_{\alpha\in z\mathbb R_{>0}}
 \operatorname{Res}_{\xi=\alpha}
       \bigl(H(X,\xi)w_Q(\xi/z)\bigr),
\]
where $H=\mathcal B_Q\widehat F$ and the jump is taken from decreasing
argument minus increasing argument. On each closed sector avoiding the
singular rays and each compact set of colors, it satisfies
\[
 \left|S(X,z)-\sum_{k<N}c_kB_k^{(A)}(X)z^k\right|
 \le CD^N M_N(Q)|z|^N\qquad(N\ge N_0).
\]
\end{theorem}

\begin{proof}
By \eqref{eq:bridge}, the Borel transform is
$\sum_kc_kp_kv^k$ with
\[
 p_k=\frac{(a,b;Q)_k}{(Q;Q)_k},
 \qquad v=-\frac QA\xi.
\]
Fix a compact set of colors and choose $K$ sufficiently large that
$|aQ^K|,|bQ^K|<1$ uniformly. Cancellation of infinite products gives
$p_{K+\ell}=p_KC_K\Psi_K(Q^\ell)$, where
\[
 C_K=\frac{(aQ^K,bQ^K;Q)_\infty}{(Q^{K+1};Q)_\infty}
\]
and
\[
 \Psi_K(s)=\frac{(Q^{K+1}s;Q)_\infty}
 {(aQ^Ks,bQ^Ks;Q)_\infty}
 =\frac{(Q^{K+1}s;Q)_\infty}
 {\prod_{j\ge0}(1-AQ^{K+j}sX+A^2Q^{2K+2j}s^2)}.
\]
The second expression shows that $\Psi_K$, $p_K$ and $C_K$ are
holomorphic functions of $X$; no choice of $y$ is involved. Choose $K$
so large that $\Psi_K(s)$ and all $f_\nu(Q^Ks)$ are holomorphic on a
common disk of radius greater than one, for all colors in the compact
set. Then
\[
 \Psi_K(s)f_\nu(Q^Ks)=\sum_{m\ge0}e_{\nu,m}s^m,
 \qquad \sum_m|e_{\nu,m}|<\infty,
\]
uniformly on the compact set. Write $\theta_v=v\,d/dv$, so that
$P(\theta_v)v^k=P(k)v^k$. For $|v|$ small, the part of the Borel
transform coming from the indices $k\ge K$ of the $\nu$th term is
\[
 p_KC_K
 P_\nu(\theta_v)
 \left[
  (\lambda_\nu v)^K
  \sum_{m\ge0}\frac{e_{\nu,m}}{1-\lambda_\nu Q^m v}
 \right],
\]
as one sees by expanding $(1-\lambda_\nu Q^mv)^{-1}$ in a geometric
series and interchanging the two absolutely convergent sums. The series
in brackets converges normally away from the points
$v=\lambda_\nu^{-1}Q^{-m}$ and has simple poles there; applying the
differential operator $P_\nu(\theta_v)$ raises the pole order by at
most $\deg P_\nu$. With $v=-(Q/A)\xi$ this gives the stated poles. On a
closed sector that avoids the ray through these poles,
$|1-\lambda_\nu Q^mv|\ge c(1+|\lambda_\nu Q^mv|)$, and each
$\theta_v$-derivative of $(1-\lambda_\nu Q^mv)^{-1}$ is bounded by a
constant independent of $m$. Hence $P_\nu(\theta_v)$ applied to the
bracket grows at most polynomially, and so does the Borel transform. The
remaining finitely many indices $k<K$ contribute a polynomial in $v$.

By Theorem~\ref{thm:kernel} the Laplace integral converges absolutely
and locally uniformly in $X$, so $S$ is entire in $X$. At a node,
$p_k=0$ for all $k$ beyond the index of the node, the Borel transform
is a polynomial, and the moment identity \eqref{eq:moment-identity}
reproduces the terminating Newton polynomial term by term.
The displayed tail representation is a finite sum of polynomial
prefactors times Euler derivatives of resolvents with absolutely
summable, geometrically decaying symbols.
Theorem~\ref{gen:thm:sumjump}, Remark~\ref{gen:rem:tails}, and
Proposition~\ref{gen:prop:jets} therefore give the lateral limits and
the residue formula; a rotation of $z$ treats each singular ray.
Finally, Lemma~\ref{lem:remainder}, applied on each closed subsector
with constants that are uniform on the compact set of colors, gives
the remainder bound.
\end{proof}

Changing finitely many coefficients adds the same polynomial in $z$ to
the formal series and, by the moment identity
\eqref{eq:moment-identity}, to its sum.

\subsection{The node divisor}

\begin{theorem}\label{thm:E}
The product
\begin{equation}\label{eq:E-product}
 E_n(X)=\prod_{j\ge0}(1-AQ^jX+A^2Q^{2j})
\end{equation}
is an entire function with simple zeros exactly at the nodes
$x_j^{(n)}$, and
\begin{equation}\label{eq:E-tail}
 E_n(X)=-A(X-x_0^{(n)})E_{n+2}(X).
\end{equation}
If two entire functions agree at all nodes, their difference is
$E_n(X)H(X)$ for a unique entire function $H$. Thus the entire
functions that vanish at all nodes form the principal ideal
$E_n\cO(\C)$.
\end{theorem}

\begin{proof}
On $|X|\le R$ the $j$th factor differs from one by at most
$AQ^jR+A^2Q^{2j}$, which is summable in $j$, so the product converges
normally and is entire. The $j$th factor equals
$-AQ^j(X-x_j^{(n)})$, so it has the single zero $x_j^{(n)}$, and the
nodes are distinct because they increase strictly with $j$. Removing the
factor with $j=0$ and noting that the remaining product is $E_{n+2}$
proves \eqref{eq:E-tail}. If an entire function vanishes at all nodes,
its quotient by $E_n$ has removable singularities at the simple zeros of
$E_n$ and no other singularities.
\end{proof}

In particular, if two sums in different sectors are entire in the
color and both equal the terminating Newton polynomials at the nodes,
their difference lies in $E_n\mathcal O(\mathbb C)$.

\subsection{Polynomial operations under summation}

Summation commutes with multiplication by polynomials in the auxiliary
variable and with positive dilations. Both facts are used for rank
raising.

\begin{lemma}\label{lem:poly-action}
Suppose $g=\cB_p f$ is holomorphic near zero and on an open sector,
with polynomial growth on every closed subsector. Then, for $z$ in the
sector,
\begin{equation}\label{eq:poly-action}
 \cS_p(P(z)f(z))=P(z)\cS_pf(z),\qquad P\in\C[z].
\end{equation}
Summation also commutes with every dilation $z\mapsto cz$ by a
positive real number $c$, in particular with $z\mapsto p^\ell z$ for
$\ell\in\Z$.
\end{lemma}

\begin{proof}
It suffices to treat $P(z)=z$. Comparing coefficients, using
$M_{k+1}(p)=M_k(p)(1-p^{k+1})/p^{k+1}$, gives
\[
 \cB_p(zf)(\xi)=\xi\sum_{r\ge1}p^rg(p^r\xi).
\]
The series converges normally. If $|g(\xi)|\le C(1+|\xi|^K)$ on a
closed subsector, then
\[
 \left|\xi\sum_{r\ge1}p^rg(p^r\xi)\right|
 \le C|\xi|\left(\frac p{1-p}+
                   \frac{p^{K+1}|\xi|^K}{1-p^{K+1}}\right),
\]
so $\cB_p(zf)$ again grows polynomially. In the Laplace integral of the
$r$th term substitute $u=p^rt$. The kernel factors that appear
telescope, because the functional equation of $w_p$ gives
\[
 up^{-r}w_p(up^{-r})
 =w_p(up^{-(r-1)})-w_p(up^{-r}),
\]
and their sum over $r\ge1$ is $w_p(u)$. These factors are positive and
$g$ grows polynomially, so the interchange of sum and integral is
justified by absolute convergence. The result is $z\cL_pg(z)$, which
proves \eqref{eq:poly-action} for $P(z)=z$, and induction gives the
general case. For $c>0$ the Borel transform of $f(cz)$ is $g(c\xi)$,
and the substitution is visible directly in the integral defining
$\cL_p$.
\end{proof}

\subsection{Divided differences and the change of node parameter}

For an entire function of the color put
\begin{equation}\label{eq:Dq-color}
 (\mathcal D_qf)(y+y^{-1})=
 \frac{f(qy+q^{-1}y^{-1})-f(q^{-1}y+qy^{-1})}
      {(q-q^{-1})(y-y^{-1})}.
\end{equation}
Here $[k]_q=(q^k-q^{-k})/(q-q^{-1})$,
$[k]_q!=\prod_{j=1}^k[j]_q$, and $[0]_q!=1$.

\begin{lemma}\label{lem:divided-newton-basis}
The apparent singularities at $y=\pm1$ are removable,
$\mathcal D_q$ is a continuous endomorphism of $\mathcal O(\mathbb C_X)$,
and
\begin{equation}\label{eq:Dq-Newton}
 \mathcal D_qB_k^{(A)}=[k]_qB_{k-1}^{(qA)}\quad(k\ge1),
 \qquad\mathcal D_qB_0^{(A)}=0.
\end{equation}
\end{lemma}
\begin{proof}
Both the numerator and the denominator of \eqref{eq:Dq-color} change
sign under $y\mapsto y^{-1}$, so the quotient is invariant. At $y=\pm1$
the two arguments $X(qy)$ and $X(y/q)$ coincide, so the numerator
vanishes there and the singularities are removable. A holomorphic
function of $y\in\C^\times$ that is invariant under $y\mapsto y^{-1}$ is
a holomorphic function of $X=y+y^{-1}$, so $\mathcal D_qf$ is entire.
Continuity on $\mathcal O(\mathbb C)$ follows from Cauchy's estimates
near $y=\pm1$ and from the fact that $X(q^{\pm1}y)$ stays in a compact
set when $X$ does.

For \eqref{eq:Dq-Newton} first take $A$ generic. Evaluate
$\mathcal D_qB_k^{(A)}$ at $y=AqQ^j$ with $0\le j\le k-2$. Then
$X(qy)=x_{j+1}^{(A)}$ and $X(y/q)=x_j^{(A)}$ are both zeros of
$B_k^{(A)}$, so $\mathcal D_qB_k^{(A)}$ vanishes at
$X(AqQ^j)=x_j^{(Aq)}$, $0\le j\le k-2$, which are the zeros of
$B_{k-1}^{(Aq)}$. Moreover $\mathcal D_q$ maps a monic polynomial of
degree $k$ to a polynomial of degree $k-1$ with leading coefficient
$[k]_q$, as the leading terms $y^k$ show. This proves
\eqref{eq:Dq-Newton} for generic $A$, and since both sides are
polynomial in $A^{\pm1}$, it holds for all $A\ne0$.
\end{proof}

\subsection{Rank raising for analytic tails}

The figure-eight sums of different ranks are related by the rank ladder
of Theorem~\ref{thm:rank-ladder}. It is a special case of the following
identity, valid for every analytic tail. For a node parameter
\(0<A<1\) write
\[
 \widehat F_{A,c}(X,z)=\sum_{k\ge0}c_k B_k^{(A)}(X)z^k,
 \qquad S_{A,c}=\mathcal L_Q\mathcal B_Q\widehat F_{A,c}.
\]

\begin{theorem}\label{thm:general-rank-raising}
Let \(c\) have an analytic tail, let \(d\) be a nonnegative integer, and set
\begin{equation}\label{eq:raised-coefficients}
 c^{\langle d\rangle}_k
   =\frac{[k+d]_q!}{[k]_q![d]_q!}\,c_{k+d}.
\end{equation}
Then \(c^{\langle d\rangle}\) has an analytic tail with the same possible
singular directions, and on every common sector of summation
\begin{equation}\label{eq:general-rank-raising}
 S_{Aq^d,c^{\langle d\rangle}}(X,z)
       =\frac{z^{-d}}{[d]_q!}\mathcal D_q^d S_{A,c}(X,z).
\end{equation}
If all $\lambda_\nu$ are positive, the same identity holds for the
lateral sums on the negative axis, and the jumps satisfy
\begin{equation}\label{eq:rank-jump-transport}
 \Delta_{Aq^d,c^{\langle d\rangle}}(X,r)
       =\frac{(-r)^{-d}}{[d]_q!}\mathcal D_q^d
                       \Delta_{A,c}(X,r).
\end{equation}
\end{theorem}

For a single tail \(c_k=\lambda^k f(Q^k)\) the raised tail is
\begin{equation}\label{eq:raised-symbol}
 c^{\langle d\rangle}_k
 =(\lambda q^{-d})^k f_d(Q^k),\qquad
 f_d(s)=\lambda^d\frac{(Qs;Q)_d}{(Q;Q)_d}f(Q^ds),
\end{equation}
and a polynomial factor \(P(k)\) is replaced by \(P(k+d)\).

\begin{proof}
By \eqref{eq:Dq-Newton}, applied coefficientwise,
\[
 \mathcal D_q^d\widehat F_{A,c}
 =[d]_q!z^d\widehat F_{Aq^d,c^{\langle d\rangle}}.
\]
Since $[k+d]_q!/([k]_q![d]_q!)=q^{-dk}(Q^{k+1};Q)_d/(Q;Q)_d$,
substituting the tail of \(c\) gives \eqref{eq:raised-symbol}. The new
tail has $\lambda q^{-d}$ in place of $\lambda$ and node parameter
$Aq^d$, so by Theorem~\ref{thm:analytic-symbol} its poles lie at
$-(Aq^d/Q)(\lambda q^{-d})^{-1}Q^{-m}=-(A/Q)\lambda^{-1}Q^{d-m}$, on the
same rays as before. By the continuity in Lemma~\ref{lem:divided-newton-basis},
$\mathcal D_q^d$ commutes with the continued Borel transform and with
the Laplace integral, and Lemma~\ref{lem:poly-action} handles the factor
$z^d$; this proves \eqref{eq:general-rank-raising}. The lateral limits
of Lemma~\ref{lem:lateral-tail} are locally uniform in $X$, so the same
finite combinations of color values can be taken on the negative axis,
where $z^{-d}=(-r)^{-d}$ on both sides. This gives
\eqref{eq:rank-jump-transport}. When the tail contains polynomial
factors $P(k)$, the boundary values are provided by
Appendix~\ref{app:jets}.
\end{proof}

Write $\mathsf R_d$ for the operation $c\mapsto c^{\langle d\rangle}$,
combined with the change $A\mapsto Aq^d$ of the node parameter.
Canceling the intermediate factorials in
\eqref{eq:raised-coefficients} gives the composition law
\begin{equation}\label{eq:raising-composition}
 \mathsf R_e\mathsf R_d
 =\frac{[d+e]_q!}{[d]_q![e]_q!}\mathsf R_{d+e}.
\end{equation}

\subsection{The figure-eight rank ladder}

The figure-eight Newton coefficients are
\begin{equation}\label{eq:D}
 D_{n,k}
 =q^{-(n-2)k}\frac{(Q^{n-1};Q)_k}{(Q;Q)_k}
 =\frac{[k+n-2]_q!}{[k]_q![n-2]_q!}.
\end{equation}
This is the case $m=1$ of Proposition~\ref{pa:knot-expansion}. At rank
two the coefficients are all one. Write
$S_n(X,z)=\mathcal S_Q\widehat F_{4_1,n}(X,z)$. The coefficients
$D_{n,k}$ have an analytic tail with the single base
$\lambda=q^{-(n-2)}>0$, so by Theorem~\ref{thm:analytic-symbol} the only
singular direction is the negative axis: the sums are defined for
$z\in\Omega$ and have upper and lower boundary values on the negative
axis.

\begin{theorem}\label{thm:rank-ladder}
With $d=n-2$, the figure-eight sums satisfy
\begin{equation}\label{eq:rank-ladder}
 S_n(X,z)=\frac{z^{-d}}{[d]_q!}\mathcal D_q^dS_2(X,z),
                         \qquad z\in\Omega.
\end{equation}
The same identity holds for both lateral sums and for their jumps.
\end{theorem}
\begin{proof}
The coefficient of the rank-two figure-eight Newton series is one.
The raised sequence in \eqref{eq:raised-coefficients} is
$[k+d]_q!/([k]_q![d]_q!)=D_{n,k}$. Apply
Theorem~\ref{thm:general-rank-raising}, including its statement on
lateral sums.
\end{proof}

The figure-eight knot is special here because its rank-$n$ coefficients
are exactly the raised rank-two coefficients. For $m\ge2$ the twist
knots $L_m$ of Section~\ref{pa:data} have no ladder of this form: their
coefficients contain the rank-dependent factor $g_{m,k}^{(n)}$, and
raising the rank-two coefficients does not produce the rank-$n$ ones.

\section{Twist knots and positive spectral kernels}
\label{pa:data}

We use the product notation of Section~\ref{sec:foundations},
now also for complex $Q$ with $0<|Q|<1$. The parameters used in the asymptotic
theorems are
\begin{equation}\label{pa:parameters}
\begin{gathered}
 h>0,\qquad Q=e^{-h},\qquad q=e^{-h/2},\qquad
 n\in\mathbb Z_{\ge2},\\
 A=q^n=Q^{n/2},\qquad
 a=A^2=Q^n,\qquad \lambda_n=Q^{2-n}.
\end{gathered}
\end{equation}
Unless otherwise stated, the rank $n$ and the twist parameter $m$ are fixed
as $h$ tends to zero; Section~\ref{pa:growing-rank} allows the rank to grow. The symmetric color is the one-row partition $(r)$,
or equivalently $\operatorname{Sym}^r\mathbb C^n$ after specialization.

Let $K_p$ be the twist knots in the convention of Kawagoe~\cite{Kawagoe},
and set $L_m=\overline{K_{-m}}$ for $m\ge1$. In this convention
$K_1$ is the left-handed trefoil and $L_1=4_1$.
We write $J_r^{(n)}(K;q)$ for the reduced, zero-framed symmetric
colored HOMFLY--PT invariant specialized at $A=q^n$; the unknot has
value one. At rank two, $J_r^{(2)}(K;q)$ is the standard reduced
colored Jones polynomial of dimension $r+1$, in variable $Q=q^2$.

Define the finite coefficients for an independent $A$ by
\begin{equation}\label{pa:g-recursion}
 g_{1,k}(A;Q)=1,\qquad
 g_{m,k}(A;Q)=\sum_{j=0}^k
       A^{2j}Q^{j(j-1)}\genfrac{[}{]}{0pt}{}{k}{j}_Q g_{m-1,j}(A;Q),
\end{equation}
and put $g_{m,k}^{(n)}=g_{m,k}(Q^{n/2};Q)$. Thus
\begin{equation}\label{pa:rank-recursion}
 g_{m,k}^{(n)}=\sum_{j=0}^k
       Q^{j(j+n-1)}\genfrac{[}{]}{0pt}{}{k}{j}_Q g_{m-1,j}^{(n)}.
\end{equation}
For reference we repeat the basis and the rank factor:
\begin{align}
 x_r^{(n)}&=q^{n+2r}+q^{-n-2r},\qquad
 B_k^{(n)}(X;h)=\prod_{j=0}^{k-1}(X-x_j^{(n)}),
 \label{pa:basis}\\
 D_{n,k}&=q^{-(n-2)k}\frac{(Q^{n-1};Q)_k}{(Q;Q)_k}.
 \label{pa:rank-factor}
\end{align}

\begin{proposition}\label{pa:knot-expansion}
For $m\ge1$, $n\ge2$, and $r\ge0$, the invariant just specified has
the finite expansion
\begin{equation}\label{pa:finite-knot-formula}
 J_r^{(n)}(L_m;q)
 =\sum_{k=0}^r D_{n,k}g_{m,k}^{(n)}B_k^{(n)}(x_r^{(n)};h).
\end{equation}
The equality is an identity of Laurent polynomials, with the
specialization $q=e^{-h/2}$ understood on the right when $h>0$.
\end{proposition}

\begin{proof}
Write $\{z\}=z-z^{-1}$ and
$[j]_q=(q^j-q^{-j})/(q-q^{-1})$.
Kawagoe's skein formula \cite[Theorem~5.1 and Corollary~5.2]{Kawagoe}
for the reduced zero-framed invariant is
\begin{equation}\label{pa:Kawagoe}
 H_r(K_p;q,A)=\sum_{k=0}^r F_k^{(p)}(A,q)
       \frac{[r]_q!}{[k]_q![r-k]_q!}
       \prod_{j=0}^{k-1}\{Aq^{r+j}\}\{Aq^{j-1}\}.
\end{equation}
Its finite spectral expression is
\[
 F_k^{(p)}=\sum_{j=0}^k V_{kj}\Lambda_j^p,\qquad
 \Lambda_j=A^{2j}q^{2j(j-1)},\qquad F_k^{(-1)}=1,
\]
where
\[
 V_{kj}=(-1)^jA^kq^{k(k-1)/2}
       \frac{[k]_q!}{[j]_q![k-j]_q!}
       \frac{\{Aq^{2j-1}\}}
       {\prod_{b=j-1}^{k+j-1}\{Aq^b\}}.
\]
We now convert this formula into the recursion \eqref{pa:g-recursion}.
Let $a=A^2$ and define
\[
 \mathsf B_{ki}=(-1)^{k-i}a^k
       Q^{\binom{k}{2}+\binom{i}{2}}\genfrac{[}{]}{0pt}{}{k}{i}_Q.
\]
The conversion of symmetric to Gaussian binomials rewrites the
spectral matrix as
\[
 V_{kj}=(-1)^{k+j}a^kQ^{\binom{k}{2}+\binom{j}{2}}
 \genfrac{[}{]}{0pt}{}{k}{j}_Q\frac{1-aQ^{2j-1}}{(aQ^{j-1};Q)_{k+1}}.
\]
For $d\ge0$ the finite identity
\begin{equation}\label{pa:finite-sum-identity}
 S_d(a):=\sum_{t=0}^d \genfrac{[}{]}{0pt}{}{d}{t}_Q
       a^tQ^{t(t-1)}(aQ^t;Q)_{d-t}=1
\end{equation}
follows from the Gaussian Pascal identity: $S_0=1$ and
$S_{d+1}(a)=(1-aQ^d)S_d(a)+aQ^dS_d(aQ)$.
For $t=i-j$, cancellation of the displayed products gives
\[
 \frac{(\mathsf BV)_{kj}}{\Lambda_jV_{kj}}
 =\sum_{t=0}^{k-j}\genfrac{[}{]}{0pt}{}{k-j}{t}_Q
       (aQ^{2j})^tQ^{t(t-1)}
       (aQ^{2j+t};Q)_{k-j-t}=1.
\]
Thus $F^{(p)}=\mathsf B F^{(p-1)}$.
Finite Gaussian inversion gives
\[
 (\mathsf B^{-1})_{ij}
 =A^{-2j}Q^{-j(j-1)}\genfrac{[}{]}{0pt}{}{i}{j}_{Q^{-1}}.
\]
Indeed its product with $\mathsf B$ reduces off the diagonal to
the finite sum $(1;Q)_{i-j}=0$, and has diagonal entries one.
These identities are first proved for generic parameters and then
extend by clearing denominators.

Mirroring replaces $(A,q)$ by $(A^{-1},q^{-1})$.
The symmetric quantum binomial in \eqref{pa:Kawagoe} is unchanged,
and the $2k$ braces contribute an even number of minus signs.
Under this substitution $\mathsf B^{-1}$ becomes the matrix
$P_{ij}=A^{2j}Q^{j(j-1)}\genfrac{[}{]}{0pt}{}{i}{j}_Q$.
Since $F^{(-1)}=\mathbf1$, the differential coefficients of $L_m$
are $P^{m-1}\mathbf1$, which is \eqref{pa:g-recursion}.
This evolution is also recorded in \cite[(8)--(10)]{MorozovTselousov}.

Finally, at $A=q^n$,
\[
 x_r^{(n)}-x_j^{(n)}
 =\{q^{r-j}\}\{q^{n+r+j}\},\qquad
 D_{n,k}=
 \frac{\prod_{j=0}^{k-1}\{q^{n-1+j}\}}
      {\prod_{j=1}^k\{q^j\}}.
\]
Each summand of \eqref{pa:Kawagoe} is therefore exactly the
corresponding summand of \eqref{pa:finite-knot-formula}.
\end{proof}

\begin{lemma}\label{pa:interpolation}
For fixed $m$, the values $g_{m,k}(A;Q)$ are bounded uniformly in
$k$ on compact subsets of $\{0<|Q|<1\}\times\mathbb C_A$.
They have the entire interpolation in $w$
\begin{equation}\label{pa:interpolating-symbol}
 G_1(A,w;Q)=1,\qquad
 G_m(A,w;Q)=\sum_{j\ge0}
 \frac{A^{2j}Q^{j(j-1)}(w;Q^{-1})_j}{(Q;Q)_j}
                     g_{m-1,j}(A;Q),
\end{equation}
with $G_m(A,Q^k;Q)=g_{m,k}(A;Q)$.
The series is normally convergent jointly on compact subsets of
$0<|Q|<1$, $A,w\in\mathbb C$.
\end{lemma}

\begin{proof}
On each such compact set the Gaussian binomials are bounded
uniformly in both indices, since their numerator products are
bounded above and their denominator products are bounded away from
zero. Induction in \eqref{pa:g-recursion}, using the convergent
series $\sum_j R^{2j}|Q|^{j(j-1)}$, proves the first assertion.
On $|w|\le T$,
$|(w;Q^{-1})_j|\le(1+T)^j|Q|^{-j(j-1)/2}$.
The remaining quadratic power proves normal convergence of
\eqref{pa:interpolating-symbol}. At $w=Q^k$ the terms with $j>k$
vanish, while $(Q^k;Q^{-1})_j/(Q;Q)_j=\genfrac{[}{]}{0pt}{}{k}{j}_Q$ for $j\le k$.
\end{proof}

If $X=y+y^{-1}$, direct factorization gives
\[
 B_k^{(n)}(X;h)=(-1)^kA^{-k}Q^{-k(k-1)/2}(Ay,A/y;Q)_k.
\]
Consequently the $Q$-Borel transform of the Newton series is
\begin{equation}\label{pa:Borel}
 H_{m,n,h}(X,\xi)
 =\sum_{k\ge0}g_{m,k}^{(n)}u_{n,k}(X)(-\lambda_n\xi)^k,
 \qquad
 u_{n,k}(X)=\frac{(Q^{n-1},Ay,A/y;Q)_k}{(Q;Q)_k^2}.
\end{equation}
It converges normally for $|\lambda_n\xi|<1$ and locally uniformly
in $X$. Symmetric products in $y,y^{-1}$ are interpreted as
functions of $X$; for instance
$(Ay,A/y;Q)_k=\prod_{j<k}(1-AQ^jX+A^2Q^{2j})$.
Theorem~\ref{thm:analytic-symbol} applies to these coefficients, and
the sum of the Newton series at $z=1$ is the entire function of the
color
\begin{equation}\label{pa:sum}
 S_{m,n,h}(X)=\frac1h\int_0^\infty
             \frac{H_{m,n,h}(X,t)}{(-t;Q)_\infty}\,dt.
\end{equation}

\subsection{The positive kernel}\label{a:kernel-remainder}
\label{pa:positive-kernel-section}

In this subsection $Q,a\in(0,1)$ are independent, and we write
$c_{m,k}^{(a)}=g_{m,k}(\sqrt a;Q)$.
Set
\[
 P^{(a)}_{k,j}=a^jQ^{j(j-1)}\genfrac{[}{]}{0pt}{}{k}{j}_Q,\qquad
 \mathsf K^{(a)}_{k,j}=
       \frac{(a;Q)_k}{(a;Q)_j}P^{(a)}_{k,j}.
\]
These finite weights are a specialization of the $q$-beta-binomial
weights \cite[Appendix~A.2]{CMPPush}. The coefficient recurrence also
requires an infinite dual kernel. We now construct it and show that it
is nonnegative, with row sums one.
For integers $r\ge0$ define
\begin{equation}\label{pa:C-positive}
 C_r=(a;Q)_\infty
 \sum_{j\ge0}\frac{a^jQ^{j(j-1)+rj}}{(Q;Q)_j(a;Q)_j}.
\end{equation}

\begin{lemma}\label{pa:duality}
One has $0<C_r\le C_0=1$ and
\begin{equation}\label{pa:C-alternating}
 C_r=\sum_{j\ge0}\frac{(Q^r;Q)_j}{(Q;Q)_j}
                 (-a)^jQ^{j(j-1)/2}.
\end{equation}
The nonnegative kernel
\begin{equation}\label{pa:dual-kernel}
 \widehat{\mathsf K}^{(a)}_{r,l}
 =\frac{a^l(Q^r;Q)_l}{(Q;Q)_l}C_{r+l}
\end{equation}
has row sums one and satisfies
\begin{equation}\label{pa:duality-identity}
 \sum_{j=0}^k\mathsf K^{(a)}_{k,j}Q^{rj}
 =\sum_{l\ge0}\widehat{\mathsf K}^{(a)}_{r,l}Q^{kl}.
\end{equation}
The finite kernel $\mathsf K^{(a)}$ also has row sums one.
\end{lemma}

\begin{proof}
In \eqref{pa:C-positive}, replace
$(a;Q)_\infty/(a;Q)_j$ by $(aQ^j;Q)_\infty$ and use Euler's
product expansion. Put $N=j+l$ in the resulting double series.
The quadratic exponent becomes
$\binom N2+\binom j2+rj$, and the finite $q$-binomial theorem
gives \eqref{pa:C-alternating}. The rearrangement is absolute:
$a<1$, both denominator products are bounded away from zero,
and all quadratic exponents before regrouping are nonnegative.
The alternating formula gives $C_0=1$, and the positive formula
gives $0<C_r\le C_0$.

We use the following confluent form of Heine's symmetry
\cite{GasperRahman,DLMFHeine}, with convention
\[
 {}_1\phi_1(B;C;Q,z)
 =\sum_{j\ge0}\frac{(B;Q)_j}{(C;Q)_j(Q;Q)_j}
                       (-z)^jQ^{\binom j2}:
\]
\begin{equation}\label{pa:Heine}
 (C;Q)_\infty\,{}_1\phi_1(B;C;Q,z)
 =(z;Q)_\infty\,{}_1\phi_1(Bz/C;z;Q,C).
\end{equation}
It follows from \cite[Eq.~17.6.7]{DLMFHeine} by substituting
$(a,b,c,z)=(B,b,C,z/b)$ and letting $b\to\infty$.
A direct proof is as follows. Expand
$(CQ^j;Q)_\infty$ on the left by Euler's identity. After
collecting total degree $N$ in $C,z$, the remaining polynomial is
\[
 \sum_{j=0}^N\genfrac{[}{]}{0pt}{}{N}{j}_Q(B;Q)_jz^jC^{N-j}.
\]
Expanding $(B;Q)_j$ and using
$\genfrac{[}{]}{0pt}{}{N}{j}_Q\genfrac{[}{]}{0pt}{}{j}{r}_Q=\genfrac{[}{]}{0pt}{}{N}{r}_Q\genfrac{[}{]}{0pt}{}{N-r}{j-r}_Q$
rewrites it as
\[
 \sum_{r=0}^N\genfrac{[}{]}{0pt}{}{N}{r}_Q(-Bz)^rQ^{\binom r2}
 \sum_{l=0}^{N-r}\genfrac{[}{]}{0pt}{}{N-r}{l}_Qz^lC^{N-r-l}.
\]
The inner sum is symmetric in $C,z$, proving
\eqref{pa:Heine}. The quadratic powers justify normal convergence
of the series, and removable parameter values follow by continuation.

For $|at|<1$ define
\[
 K_r(t)=\frac{(a;Q)_\infty}{(at;Q)_\infty}
 \sum_{j\ge0}
 \frac{a^jQ^{j(j-1)+rj}(t;Q^{-1})_j}{(Q;Q)_j(a;Q)_j}.
\]
At $t=Q^k$ this is the left side of
\eqref{pa:duality-identity}. Apply \eqref{pa:Heine} with
$(B,C,z)=(t^{-1},a,atQ^r)$. Since $r$ is an integer, it gives
\[
 K_r(t)=\sum_{j\ge0}
 \frac{(Q^r;Q)_j(-a)^jQ^{\binom j2}}
      {(Q;Q)_j(at;Q)_{r+j}}.
\]
The substitution is initially made for $t\ne0$ and extends
holomorphically to zero. Expand
$1/(at;Q)_{r+j}=\sum_{l\ge0}(Q^{r+j};Q)_l(at)^l/(Q;Q)_l$.
The identity
$(Q^r;Q)_j(Q^{r+j};Q)_l=(Q^r;Q)_l(Q^{r+l};Q)_j$
and \eqref{pa:C-alternating} show that
$[t^l]K_r(t)=\widehat{\mathsf K}^{(a)}_{r,l}$.
These rearrangements are absolute on compact subsets of $|at|<1$.
At $t=1$, only $j=0$ remains in the first formula for $K_r$,
so its nonnegative coefficients sum to one. At $r=0$ the dual
row is concentrated at zero, which proves the final assertion.
\end{proof}

\begin{proposition}\label{pa:positive-moments}
There are probability distributions $p^{(a)}_{m,l}$ on
$\mathbb Z_{\ge0}$ such that
\begin{equation}\label{pa:positive-moment-identity}
 (a;Q)_kc_{m,k}^{(a)}=\sum_{l\ge0}p^{(a)}_{m,l}Q^{kl}.
\end{equation}
They are specified by
\[
 p^{(a)}_{1,l}=(a;Q)_\infty\frac{a^l}{(Q;Q)_l},\qquad
 p^{(a)}_{m,l}=\sum_{r\ge0}p^{(a)}_{m-1,r}
                              \widehat{\mathsf K}^{(a)}_{r,l}.
\]
In particular $1\le c_{m,k}^{(a)}\le(a;Q)_k^{-1}$ and
\begin{equation}\label{pa:positive-symbol}
 \frac{G_m(\sqrt a,t;Q)}{(at;Q)_\infty}
 =\frac1{(a;Q)_\infty}\sum_{l\ge0}p^{(a)}_{m,l}t^l,
 \qquad |at|<1.
\end{equation}
\end{proposition}

\begin{proof}
The $q$-binomial theorem proves the moment identity at $m=1$.
Multiply the coefficient recursion by $(a;Q)_k$ and apply
\eqref{pa:duality-identity}. Tonelli's theorem proves the induction,
and the row sums give normalization. The lower coefficient bound
comes from the $j=0$ term of the recursion; the upper bound follows
from the moment identity.
Since $C_r\le1$, one has
$p^{(a)}_{m,l}\le a^l/(Q;Q)_l$. Thus the series in
\eqref{pa:positive-symbol} is holomorphic for $|at|<1$.
Its values at $t=Q^k$ agree with the left side by the moment identity.
The identity theorem proves the displayed interpolation formula.
\end{proof}

\subsection{The remainder formula}
\label{pa:tail-section}

Return to \eqref{pa:parameters}. Let $p_{m,n,i}=p^{(Q^n)}_{m,i}$
and, for an integer $K\ge1$, define the tilted probability
\begin{equation}\label{pa:tilted-p}
 \pi_{m,n,K,i}=
 \frac{p_{m,n,i}Q^{Ki}}{(Q^n;Q)_K g_{m,K}^{(n)}}.
\end{equation}
Write $D=X-2$ and $\tau(v)=2\cosh v-2$.
The ratio of $u_{n,k}$ from \eqref{pa:Borel} to $(Q^n;Q)_k$
has the tail generating function
\begin{equation}\label{pa:psi-original}
 \psi_{n,K,D}(v)=
 \frac{(Q^{K+1}v;Q)_\infty^2}
 {(1-Q^{K+n-1}v)
  (Q^{K+n/2}yv,Q^{K+n/2}v/y;Q)_\infty}.
\end{equation}
We factor it as $F_{n,K,h}(v)\psi^*_{n,K,D}(v)$, where
\begin{align}
 \psi^*_{n,K,D}(v)
 &=\frac{(Q^{K+n/2}v,Q^{K+n/2+1}v;Q)_\infty}
       {(Q^{K+n/2}yv,Q^{K+n/2}v/y;Q)_\infty},
 \label{pa:positive-psi}\\
 F_{n,K,h}(v)
 &=\frac{(Q^{K+1}v;Q)_\infty^2}
 {(1-Q^{K+n-1}v)
  (Q^{K+n/2}v,Q^{K+n/2+1}v;Q)_\infty}.
 \label{pa:rank-correction}
\end{align}
Define the normalized coefficients by
\begin{equation}\label{pa:mu-nu}
 \frac{\psi^*_{n,K,D}(v)}{\psi^*_{n,K,D}(1)}
   =\sum_{l\ge0}\mu_l(D)v^l,\qquad
 \frac{F_{n,K,h}(v)}{F_{n,K,h}(1)}
   =\sum_{r\ge0}\nu_rv^r.
\end{equation}

\begin{lemma}\label{pa:weight-bounds}
If $|D|<\tau((K+n/2)h)$, both coefficient sequences in
\eqref{pa:mu-nu} are absolutely summable and have sum one.
For real $0\le D<\tau((K+n/2)h)$, the weights $\mu_l(D)$
are nonnegative. If $Kh\ge U>0$, then
\begin{equation}\label{pa:nu-bound}
 \sum_r|\nu_r|\le(1-e^{-U})^{-2(n+1)}.
\end{equation}
If $D_0\ge0$ and
$|D-D_0|\le\eta<\tau((K+n/2)h)-D_0$, then
\begin{equation}\label{pa:mu-bound}
 \sum_l|\mu_l(D)|\le
 \prod_{j\ge0}
 \frac{\tau((K+n/2+j)h)-D_0+\eta}
      {\tau((K+n/2+j)h)-D_0-\eta}.
\end{equation}
In particular, if $Kh\ge U$, $\tau(U)>D_0$, and
$\eta\le(\tau(U)-D_0)/2$, then
\begin{equation}\label{pa:mu-uniform-bound}
 \sum_l|\mu_l(D)|\le
 \exp\left\{
 \frac{4\eta}{h(1-D_0/\tau(U))(e^U-1)}\right\}.
\end{equation}
At rank two $F_{2,K,h}=1$.
\end{lemma}

\begin{proof}
Each factor of \eqref{pa:positive-psi}, with
$c=Q^{K+n/2+j}$, is
\[
 \frac{(1-cv)(1-Qcv)}{(1-cv)^2-Dcv}
 =\frac{1-Qcv}{1-cv}
       \left(1-\frac{Dcv}{(1-cv)^2}\right)^{-1}.
\]
Its Taylor coefficients are jointly nonnegative in $D,v$.
The stated domain ensures convergence on a disk in $v$ larger
than the unit disk. Majorizing at $D$ by the coefficients at
$|D|$ and canceling the factors independent of $D$ in the
normalization proves \eqref{pa:mu-bound}.

Let $\|\sum a_rv^r\|_{\mathcal W}=\sum|a_r|$.
The coefficient of $v^r$ in $\log F_{n,K,h}$ is
\[
 \frac{Q^{Kr}}r\left[
 x^{n-1}+\frac{x^{n/2}+x^{n/2+1}-2x}{1-x}\right],\qquad x=Q^r.
\]
The bracket has absolute value at most $n+1$, by
$(1-x^b)/(1-x)\le\max(1,b)$ for $b\ge0$, $0\le x<1$.
Thus $\|\log F\|_{\mathcal W}\le-(n+1)\log(1-e^{-U})$.
Submultiplicativity bounds the norms of $F$ and $F^{-1}$ by
$(1-e^{-U})^{-(n+1)}$, and normalization at one proves
\eqref{pa:nu-bound}.

For \eqref{pa:mu-uniform-bound}, the logarithm of each factor
in \eqref{pa:mu-bound} is at most
$4\eta/(\tau(v)-D_0)$, and
$\tau(v)-D_0\ge(1-D_0/\tau(U))\tau(v)$.
The arguments $(K+n/2+j)h$, $j\ge0$, in \eqref{pa:mu-bound} are all at
least $(K+1)h$.
Decreasing-integrand comparison and
\[
 \int_U^\infty\frac{dv}{2\cosh v-2}=\frac1{e^U-1}
\]
give the result. All normalizations in \eqref{pa:mu-nu} have
value one at $v=1$, which proves the mass assertions.
\end{proof}

For each integer $K\ge1$, let $T_K$ be the positive random variable
with density
\begin{equation}\label{pa:T-probability}
 \frac{t^K}{hM_K(Q)(-t;Q)_\infty}\,dt.
\end{equation}
The moment formula \eqref{eq:moment-identity} shows that this density
has total mass one.

\begin{theorem}\label{pa:source}
For fixed $m,n,h$, the germ \eqref{pa:Borel} extends
meromorphically to the $\xi$-plane, with possible simple poles only at
$\xi=-Q^{n-2-j}$, $j\ge0$. The integral \eqref{pa:sum} converges
locally uniformly in $X\in\mathbb C$ and defines an entire function.
It interpolates the symmetric colored invariants:
\begin{equation}\label{pa:node-interpolation}
 S_{m,n,h}(x_r^{(n)})=J_r^{(n)}(L_m;q).
\end{equation}
For $K\ge1$ and $|X-2|<\tau((K+n/2)h)$, the exact remainder is
\begin{align}
 S_{m,n,h}(X)-\sum_{k<K}D_{n,k}g_{m,k}^{(n)}B_k^{(n)}(X;h)
 ={}&D_{n,K}g_{m,K}^{(n)}B_K^{(n)}(X;h)\nonumber\\
 &\times\sum_{i,r,l\ge0}\pi_{m,n,K,i}\nu_r\mu_l(X-2)
       \mathbb E\frac1{1+\lambda_nQ^{i+r+l}T_K}.
 \label{pa:tail}
\end{align}
In particular, its absolute value is at most
\begin{equation}\label{pa:tail-bound}
 D_{n,K}g_{m,K}^{(n)}|B_K^{(n)}(X;h)|
       \left(\sum_r|\nu_r|\right)
       \left(\sum_l|\mu_l(X-2)|\right).
\end{equation}
\end{theorem}

\begin{proof}
Put $\widetilde u_{n,k}=u_{n,k}/(Q^n;Q)_k$.
Cancellation of infinite products in \eqref{pa:psi-original}
gives the ratio of $\widetilde u_{n,K+j}$ to
$\widetilde u_{n,K}$ as the $j$th $Q$-moment of the convolution
of $\nu$ and $\mu$. Combining it with
\eqref{pa:positive-moment-identity} yields
\[
 \frac{g_{m,K+j}^{(n)}u_{n,K+j}(X)}
      {g_{m,K}^{(n)}u_{n,K}(X)}
 =\sum_{i,r,l\ge0}\pi_{m,n,K,i}\nu_r\mu_l(X-2)Q^{(i+r+l)j}.
\]
First take $X$ away from the finitely many zeros of $u_{n,K}$.
Summing geometric series in $|\lambda_n\xi|<1$ gives the
Borel tail as
\[
 g_{m,K}^{(n)}u_{n,K}(X)(-\lambda_n\xi)^K
 \sum_{i,r,l\ge0}
 \frac{\pi_{m,n,K,i}\nu_r\mu_l(X-2)}
      {1+\lambda_nQ^{i+r+l}\xi}.
\]
The weights are absolutely summable. On each compact set avoiding
the displayed pole sequence the resolvents are uniformly bounded,
which proves meromorphic continuation of the original germ.
On the positive ray their absolute values are at most one.
Theorem~\ref{thm:kernel} therefore permits termwise integration,
and
\[
 (-\lambda_n)^KM_K(Q)u_{n,K}(X)
       =D_{n,K}B_K^{(n)}(X;h)
\]
proves \eqref{pa:tail}.

For every compact set of $X$, a single sufficiently large $K$
satisfies the domain condition. All products and normalized
weights are holomorphic there with locally bounded total variation.
The formula consequently extends across the zeros excluded in the
first step and proves local normal convergence of the integral and
entireness in $X$. At a node the Borel series terminates; applying
\eqref{eq:moment-identity} term by term gives
\eqref{pa:node-interpolation} and Proposition~\ref{pa:knot-expansion}.
\end{proof}

At rank two and real $0\le X-2<\tau((K+1)h)$, all the weights in
\eqref{pa:tail} are nonnegative, and the normalized remainder is a
Stieltjes transform of a probability measure. At higher rank the
weights $\nu_r$ may be negative, but their total variation is bounded by
\eqref{pa:nu-bound}.

\subsection{Stokes separation at a fixed rank}

The coefficient of \(\xi^k\) in \eqref{pa:Borel} is the corresponding
figure-eight coefficient multiplied by \(G_m(A,Q^k;Q)\), where \(G_m\)
is the entire function of \eqref{pa:interpolating-symbol}. After a
finite initial part is removed, this is the symbol multiplication of
Theorem~\ref{gen:thm:resolvent}. Thus \(G_m\) acts as a spectral
multiplier that turns the Borel transform, the Laplace sums and the
jumps of the figure-eight knot into those of \(L_m\).

Write \(S_{m,n,h}(X,z)=\mathcal L_QH_{m,n,h}(X,\cdot)(z)\);
then \(S_{m,n,h}(X,1)=S_{m,n,h}(X)\). Put
\[
 \Delta_{m,n,h}(X,r)
   =S_{m,n,h}(X,-r+i0)-S_{m,n,h}(X,-r-i0),
 \qquad r>0,
\]
and, for \(X=y+y^{-1}\), define
\[
 \Lambda_n(X)=
       \frac{(Q^{n-1},Ay,A/y;Q)_\infty}{(Q;Q)_\infty^2},
 \qquad \rho_n=Q^{n-2},\qquad A=Q^{n/2}.
\]
The zeros of \(\Lambda_n\) are exactly the interpolation nodes.

\begin{proposition}\label{merged:twist-separation}
Fix \(h>0\), \(n\ge2\), and integers \(1\le i<j\). Let
\[
 g_{ij}(A)=G_j(A,0;Q)-G_i(A,0;Q).
\]
Then
\begin{equation}\label{merged:positive-residue-gap}
 g_{ij}(A)\ge \frac1{1-Q}\sum_{a=i}^{j-1}A^{2a}>0.
\end{equation}
At every non-node color \(X\),
\begin{align}
 \Delta_{j,n,h}(X,r)-\Delta_{i,n,h}(X,r)
 ={}&-\frac{2\pi i\Lambda_n(X)}h\,
       \frac{\rho_n}{r}w_Q(\rho_n/r)\nonumber\\
 &\qquad\times\{g_{ij}(A)+O(r/\rho_n)\},
       \qquad r\downarrow0. \label{merged:twist-jump-separation}
\end{align}
The estimate is uniform on compact sets of non-node colors.
In particular, these jumps are pairwise distinct functions of
the modulus, and each pair is separated at every sufficiently
small positive modulus.
\end{proposition}

\begin{proof}
Write \(d_m(A,t)=G_m(A,t;Q)-G_{m-1}(A,t;Q)\).
The recurrence has nonnegative weights for \(A,Q>0\).
For \(m=2\), its positive-index terms give
\[
 d_2(A,Q^k)\ge A^2\frac{1-Q^k}{1-Q}.
\]
For \(m\ge3\), subtracting two recurrences gives
\[
 d_m(A,Q^k)=\sum_{a=0}^k
 A^{2a}Q^{a(a-1)}\Gauss{k}{a}\,d_{m-1}(A,Q^a).
\]
Inductively all terms are nonnegative. The same recurrence at
\(k=1\) gives
\(G_m(A,Q;Q)=1+A^2+\cdots+A^{2m-2}\).
Keeping the \(a=1\) term proves
\[
 d_m(A,Q^k)\ge A^{2m-2}\frac{1-Q^k}{1-Q}.
\]
Let \(k\to\infty\) and telescope over \(m\) to obtain
\eqref{merged:positive-residue-gap}.

Near \(t=0\), define
\[
 \psi_n(t)=
 \frac{(Qt;Q)_\infty^2}
      {(Q^{n-1}t,Ayt,At/y;Q)_\infty},
 \qquad
 \psi_n(t)\{G_j(A,t;Q)-G_i(A,t;Q)\}
       =\sum_{\ell\ge0}a_\ell(X)t^\ell .
\]
Thus \(a_0(X)=g_{ij}(A)\). Product cancellation in
\eqref{pa:Borel}, followed by geometric summation of a sufficiently
far tail, gives
\[
 \operatorname{Res}_{\xi=-\rho_nQ^{-\ell}}
       (H_{j,n,h}-H_{i,n,h})
    =\Lambda_n(X)a_\ell(X)\rho_nQ^{-\ell}.
\]
On every compact color set, the germ of \(\psi_n\) has a common
positive radius and its coefficients satisfy
\(|a_\ell(X)|\le C D^\ell\).
The lateral residue formula of
Theorem~\ref{gen:thm:sumjump}, with its finite-tail version, now
gives the difference of jumps as
\[
 -\frac{2\pi i\Lambda_n(X)}h
     \sum_{\ell\ge0}a_\ell(X)\tau Q^{-\ell}
                         w_Q(\tau Q^{-\ell}),
 \qquad \tau=\rho_n/r.
\]
The kernel functional equation gives
\[
 \frac{\tau Q^{-\ell}w_Q(\tau Q^{-\ell})}
      {\tau w_Q(\tau)}
 =\frac{Q^{-\ell}}
        {\prod_{a=1}^{\ell}(1+\tau Q^{-a})}
 \le \tau^{-\ell}Q^{\ell(\ell-1)/2}.
\]
For \(\tau>D\), the sum of the terms with \(\ell\ge1\) after
this normalization is \(O(\tau^{-1})\), uniformly on the compact
set. The remaining term is \(g_{ij}(A)>0\), which proves
\eqref{merged:twist-jump-separation}. The factor
\(\Lambda_n(X)\) is bounded away from zero on a compact set
without nodes.
\end{proof}

\section{Critical limits at Alexander roots}\label{a:fixed-critical}

In the trace coordinate $X=t+t^{-1}$, the Alexander polynomial of
$L_m$, normalized to take the value one at $t=1$, is
\(\Delta_m(X)=1-m(X-2)\). Its zero is \(X_m=2+1/m\).
For each fixed Newton index, the classical limits of the coefficient
and basis are
\[
 D_{n,k}\longrightarrow\binom{k+n-2}{n-2},\qquad
 g_{m,k}^{(n)}\longrightarrow m^k,\qquad
 B_k^{(n)}(X;h)\longrightarrow(X-2)^k.
\]
For \(|m(X-2)|<1\) the corresponding series converges:
\[
 \sum_{k\ge0}\binom{k+n-2}{n-2}[m(X-2)]^k
             =[1-m(X-2)]^{-(n-1)}.
\]
At rank two this is the inverse Alexander polynomial of the
Melvin--Morton--Rozansky expansion \cite{GLAsymptotics,GLAnalytic}.
Proposition~\ref{pa:subcritical-nodes} proves the limit, uniformly, at
the quantum colors below the root, and
Theorem~\ref{bridge:analytic-abelian-overlap} proves it for the analytic
interpolation up to a distance of order $\sqrt h$ (respectively
$h^{2/3}$) from the root. At the root itself the geometric series
diverges, and the sum is controlled by the first decaying factor that
survives the limit: the Gaussian decay of the coefficient
$g_{m,k}^{(n)}/m^k$ when $m\ge2$ (Section~\ref{pa:gaussian-section}),
and the cubic decay of the Newton product when $m=1$
(Section~\ref{pa:scorer-section}).

\subsection{Word estimates}
\label{pa:coefficient-section}

For a word $w=(w_1,\ldots,w_k)\in\{1,\ldots,m\}^k$, let
\[
 N_s(w)=\#\{j:w_j\le s\},\qquad
 \operatorname{inv}(w)=\#\{i<j:w_i>w_j\},
\]
and set
\begin{equation}\label{pa:word-energy}
 E_{m,n}(w)=\operatorname{inv}(w)
             +\sum_{s=1}^{m-1}N_s(w)(N_s(w)+n-1).
\end{equation}
Put
\begin{equation}\label{pa:alpha-beta}
 \alpha_m=\frac{(m-1)(4m+1)}{12m},\qquad
 \beta_{m,n}=\frac{(m-1)(8m-1)}{12m}
                     +\frac{(n-2)(m-1)}2.
\end{equation}

\begin{lemma}\label{pa:energy}
If a word is chosen uniformly from $\{1,\ldots,m\}^k$, then
\begin{equation}\label{pa:energy-expectation}
 \frac{g_{m,k}^{(n)}}{m^k}=\mathbb E e^{-hE_{m,n}},\qquad
 \mathbb E E_{m,n}=\alpha_m k^2+\beta_{m,n}k.
\end{equation}
With $b_{m,n}=n(m-1)+1$, the centered energy obeys
\begin{equation}\label{pa:energy-mgf}
 \mathbb E e^{t(E_{m,n}-\mathbb E E_{m,n})}
 \le \exp\left(\frac{b_{m,n}^2}{8}t^2k^3\right)
 \qquad(t\in\mathbb R).
\end{equation}
For $m\ge2$, let $p_m=(m-1)/m$. For every $h>0$ and $k\ge0$,
\begin{equation}\label{pa:global-coefficients}
 \frac{g_{m,k}^{(n)}}{m^k}
 \le e^{-hp_m^2k^2/4}+e^{-p_m^2k/2}.
\end{equation}
\end{lemma}

\begin{proof}
Remove all occurrences of the letter $m$ and let the remaining
subword have length $j$. The cross inversions have generating
polynomial $\genfrac{[}{]}{0pt}{}{k}{j}_Q$, and the additional energy is
$j(j+n-1)$. Induction gives exactly the recursion
\eqref{pa:rank-recursion}, proving the first identity.
The random variable $N_s$ is binomial with success probability
$s/m$. Therefore
\[
 \mathbb E\operatorname{inv}(w)=\frac{m-1}{4m}(k^2-k),\qquad
 \mathbb E[N_s(N_s+n-1)]
 =\frac{s^2}{m^2}k^2+
       \left(\frac{ns}{m}-\frac{s^2}{m^2}\right)k.
\]
Summing gives the displayed mean.

Changing one letter changes the inversion count by at most $k-1$
and each quadratic count term by at most $2k+n-2$.
For $k\ge1$, the total change is at most $b_{m,n}k$.
Reveal the independent letters one at a time. The conditional range
of each increment of the associated Doob martingale has length at
most $b_{m,n}k$.
If a centered random variable has range of length $b$, its
logarithmic moment generating function has second derivative at
most $b^2/4$: this derivative is the variance under exponential
tilting, whose support has the same range. Integrating twice from
zero bounds the logarithmic moment generating function by
$b^2t^2/8$. Applying this conditional bound to the martingale
increments proves \eqref{pa:energy-mgf}.

For the last estimate, use $g_{m-1,j}^{(n)}\le(m-1)^j$,
$\genfrac{[}{]}{0pt}{}{k}{j}_Q\le\binom{k}{j}$, and $n\ge2$ in the recurrence.
The result is bounded by
$\mathbb E e^{-hJ(J+1)}$ for $J\sim\operatorname{Bin}(k,p_m)$.
On $J\ge p_mk/2$ it is at most $e^{-hp_m^2k^2/4}$.
The same bounded-range moment estimate applied to the sum of
Bernoulli variables gives
$\mathbb P(J-p_mk\le-p_mk/2)\le e^{-p_m^2k/2}$.
This proves \eqref{pa:global-coefficients}.
\end{proof}

\begin{lemma}\label{pa:weighted-coefficient-error}
For each fixed $m\ge2$ and $n\ge2$, there are constants
$a_0,C>0$ such that, with $\varepsilon=\sqrt h$ and
$t=\varepsilon k$, whenever $hk\le a_0$,
\begin{align}
 \frac{g_{m,k}^{(n)}}{m^k}&\le e^{-\alpha_mt^2/2},
 \label{pa:weighted-coefficient-bound}\\
 \left|\frac{g_{m,k}^{(n)}}{m^k}-e^{-\alpha_mt^2}\right|
 &\le C\varepsilon(t+t^3)e^{-\alpha_mt^2/2}.
 \label{pa:weighted-coefficient-rate}
\end{align}
\end{lemma}

\begin{proof}
Let $B=b_{m,n}^2/8$ and choose $a_0\le\alpha_m/(2B)$.
Jensen's inequality and \eqref{pa:energy-mgf} give
\[
 1\le\mathbb E e^{-h(E_{m,n}-\mathbb E E_{m,n})}
           \le e^{Bh^2k^3}.
\]
Since $\beta_{m,n}\ge0$ and
$Bh^2k^3\le\alpha_mhk^2/2$, the first estimate follows.
The difference in the second estimate is at most
\[
 e^{-\alpha_mhk^2}
       \left(\beta_{m,n}hk+e^{Bh^2k^3}-1\right)
 \le C(hk+h^2k^3)e^{-\alpha_mhk^2/2}.
\]
Use $hk=\varepsilon t$ and $h^2k^3=\varepsilon t^3$.
\end{proof}

The rank factor obeys the same estimate at both critical scales.

\begin{lemma}\label{pa:rank-factor-estimates}
Let $d=n-2$. For all $h>0$,
\begin{equation}\label{pa:rank-factor-products}
 D_{n,k}
 =\prod_{r=1}^d\frac{\sinh(h(k+r)/2)}{\sinh(hr/2)}
 \le e^{hdk/2}\binom{k+d}{d}.
\end{equation}
Fix $U>0$, let $0<h\le1$, and take
$\rho\in\{\sqrt h,h^{2/3}\}$. For $hk\le U$ and $t=\rho k$,
\begin{align}
 \rho^dD_{n,k}&\le C_{n,U}(1+t)^d,
 \label{pa:D-scaled-bound}\\
 \left|\rho^dD_{n,k}-\frac{t^d}{d!}\right|
 &\le C_{n,U}\rho(1+t)^{d+2}.
 \label{pa:D-scaled-rate}
\end{align}
\end{lemma}

\begin{proof}
Cancel the Pochhammer products to obtain
$D_{n,k}=e^{hdk/2}\prod_{r=1}^d(1-Q^{k+r})/(1-Q^r)$.
This proves the equality in \eqref{pa:rank-factor-products};
the inequality follows from the monotonicity of
$(1-e^{-x})/x$. For the rate, write
\[
 D_{n,k}=\binom{k+d}{d}
   \prod_{r=1}^d\frac{f(h(k+r)/2)}{f(hr/2)},\qquad
 f(x)=\frac{\sinh x}{x},\quad f(0)=1.
\]
On the bounded argument interval in question, $f(x)=1+O(x^2)$
and $f$ is bounded away from zero. Thus the product differs from
one by at most $C_{n,U}h^2(k+d)^2$.
Also $\rho^d\binom{k+d}{d}=\prod_{r=1}^d(t+\rho r)/d!$.
Its difference from $t^d/d!$ is bounded by
$C_n\rho(1+t)^{\max(d-1,0)}$.
Since $(h/\rho)^2\le\rho$ at the two specified scales,
these estimates prove \eqref{pa:D-scaled-rate}.
The case $d=0$ is exact.
\end{proof}

\subsection{The Gaussian limit}
\label{pa:gaussian-section}

For $m\ge2$ put
\begin{equation}\label{pa:gaussian-constants}
 X_m=2+\frac1m,\qquad
 U_m=\operatorname{arccosh}\left(1+\frac3{4m}\right),\qquad
 A_m=\frac{p_m^2U_m^2}{4}.
\end{equation}
Thus $\tau(U_m)=3/(2m)$ and $U_m<1$.

\begin{lemma}\label{pa:gaussian-tail}
Fix $m,n$ and $S>0$. With $K=\lceil U_m/h\rceil$, uniformly
for complex $|s|\le S$,
\begin{equation}\label{pa:gaussian-tail-rate}
 S_{m,n,h}(X_m+\sqrt h\,s)
 -\sum_{k<K}D_{n,k}g_{m,k}^{(n)}
                  B_k^{(n)}(X_m+\sqrt h\,s;h)
 =O_{m,n,S}(e^{-A_m/(2h)}).
\end{equation}
For any fixed sufficiently small $a_0>0$, the part of this finite
sum with $a_0/h\le k<K$ is also exponentially small in $1/h$.
\end{lemma}

\begin{proof}
For small $h$ all the unperturbed normalized factors through $K$
satisfy
$|1-m\tau((j+n/2)h)|\le1$, because the largest argument tends
to $U_m$ and $\tau(U_m)<2/m$. Hence, with $\varepsilon=\sqrt h$,
\begin{equation}\label{pa:gaussian-basis-bound}
 m^k|B_k^{(n)}(X_m+\varepsilon s;h)|
 \le e^{mS\varepsilon k}\qquad(k\le K).
\end{equation}
By \eqref{pa:global-coefficients},
$g_{m,K}^{(n)}/m^K\le2e^{-A_m/h}$.
The rank factor is $O_{m,n}(h^{-(n-2)})$.
Apply \eqref{pa:mu-uniform-bound} with
$D_0=1/m$, $U=U_m$, and $\eta=S\sqrt h$; its exponent is
at most $12S/((e^{U_m}-1)\sqrt h)$.
Together with \eqref{pa:nu-bound} and \eqref{pa:tail-bound},
this bounds the remainder by a fixed power of $h^{-1}$ times
$\exp(-A_m/h+O_{m,n,S}(h^{-1/2}))$, proving
\eqref{pa:gaussian-tail-rate}.

For the finite later terms use
\eqref{pa:global-coefficients}, \eqref{pa:gaussian-basis-bound},
and \eqref{pa:rank-factor-products}. On $k\ge a_0/h$, each
negative coefficient exponent dominates $mS\sqrt h\,k$.
There are $O(h^{-1})$ terms, with a fixed polynomial rank factor,
so their sum is $O(e^{-c/h})$ for some $c>0$.
\end{proof}

\begin{theorem}\label{pa:gaussian}
Fix integers $m\ge2$ and $n\ge2$. Let $\mathcal G_{m,n}$ be defined by
\eqref{pa:gaussian-profile}.
For every compact set $C\subset\mathbb C$ and every integer $j\ge0$,
\begin{equation}\label{pa:gaussian-rate}
 \sup_{s\in C}\left|
 \frac{d^j}{ds^j}\left[
 h^{(n-1)/2}S_{m,n,h}(X_m+\sqrt h\,s)
              -\mathcal G_{m,n}(s)\right]\right|
 =O_{m,n,C,j}(\sqrt h).
\end{equation}
\end{theorem}

\begin{proof}
Let $d=n-2$, $\varepsilon=\sqrt h$, and $t=\varepsilon k$.
Choose $a_0$ small enough for Lemma~\ref{pa:weighted-coefficient-error}
and for $1-m\tau((j+n/2)h)>0$ when $k\le a_0/h$ and $j<k$.
Shrinking $h$ achieves the latter condition uniformly in these indices.
For $|s|\le S$, set
$\mathcal B_k=m^kB_k^{(n)}(X_m+\varepsilon s;h)$.
Telescoping its product against $(1+m\varepsilon s)^k$ gives
\[
 |\mathcal B_k-(1+m\varepsilon s)^k|
 \le m\sum_{j<k}\tau((j+n/2)h)
                      (1+mS\varepsilon)^{k-1}
 \le C\varepsilon(1+t^3)e^{mSt}.
\]
Here $\sum_{j<k}\tau((j+n/2)h)\le C_nh^2(k+1)^3$.
The elementary inequality
\begin{equation}\label{pa:exponential-difference}
 |e^u-e^v|\le |u-v|\max(e^{\Re u},e^{\Re v})
\end{equation}
follows by integrating the exponential on the segment from $v$ to $u$.
Apply it with $u=k\log(1+m\varepsilon s)$ and $v=mst$.
The logarithmic difference is $O_S(k\varepsilon^2)$, and both
real parts are at most $mSt$. Thus
\begin{equation}\label{pa:gaussian-basis-rate}
 |\mathcal B_k-e^{mst}|
 \le C\varepsilon(1+t^3)e^{mSt}
 \qquad(k\le a_0/h).
\end{equation}

Combine \eqref{pa:weighted-coefficient-bound},
\eqref{pa:weighted-coefficient-rate},
\eqref{pa:D-scaled-rate}, and \eqref{pa:gaussian-basis-rate}.
For a constant $C$, the resulting estimate is
\begin{align}
 &\left|\varepsilon^dD_{n,k}g_{m,k}^{(n)}
       B_k^{(n)}(X_m+\varepsilon s;h)
  -\frac{t^d}{d!}e^{mst-\alpha_mt^2}\right|\nonumber\\
 &\hspace{25mm}\le
 C\varepsilon(1+t)^{d+3}e^{mSt-\alpha_mt^2/2}.
 \label{pa:gaussian-summand-rate}
\end{align}
The majorant is integrable on $t\ge0$, uniformly for $|s|\le S$.
Its sampled sums multiplied by $\varepsilon$ are bounded uniformly
for small $\varepsilon$. Summing \eqref{pa:gaussian-summand-rate}
therefore costs $O(\varepsilon)$ after the outer Riemann-sum factor
$\varepsilon$ is included.

For the limiting function
$f_s(t)=t^de^{mst-\alpha_mt^2}/d!$, both $f_s$ and its first
$t$ derivative are integrable uniformly on $|s|\le S$.
The elementary bound
\[
 \left|\varepsilon\sum_{k\ge0}f_s(\varepsilon k)
                   -\int_0^\infty f_s(t)\,dt\right|
 \le\varepsilon\int_0^\infty|f_s'(t)|\,dt
\]
gives the same rate for the Riemann sum. The limiting tail beyond
$t=a_0/\varepsilon$ is exponentially small. Lemma~\ref{pa:gaussian-tail}
controls the finite later indices and the analytic remainder.
This proves \eqref{pa:gaussian-rate} at $j=0$.
Entireness of the color function and Cauchy's formula on a slightly
larger compact set give the assertion for every fixed $j$.
\end{proof}

Let $\operatorname{erfc}(z)=1-2\pi^{-1/2}\int_0^z e^{-t^2}\,dt$.
The rank-two profile and its rank derivatives are
\begin{align}
 \mathcal G_{m,2}(s)
 &=\frac{\sqrt\pi}{2\sqrt{\alpha_m}}
   \exp\left(\frac{m^2s^2}{4\alpha_m}\right)
   \operatorname{erfc}\left(-\frac{ms}{2\sqrt{\alpha_m}}\right),
 \label{pa:erfc-profile}\\
 \mathcal G_{m,n}(s)
 &=\frac1{(n-2)!m^{n-2}}
                  \frac{d^{n-2}}{ds^{n-2}}\mathcal G_{m,2}(s).
 \label{pa:gaussian-rank-profile}
\end{align}
Integration by parts also gives
$\mathcal G_{m,2}'-m^2s\mathcal G_{m,2}/(2\alpha_m)
=m/(2\alpha_m)$.

\subsection{The Scorer limit}
\label{pa:scorer-section}

For $m=1$ the twist coefficient is identically one. The critical
scale is now $\delta=h^{2/3}$. Fix a rank $n\ge2$ and set
\begin{equation}\label{pa:scorer-cutoffs}
 L=\lfloor1/(2h)\rfloor,\qquad K=\lceil1/h\rceil.
\end{equation}
The inequalities $\tau(1/2)<1<\tau(1)<2$ imply, for sufficiently
small $h$ depending on $n$,
\begin{equation}\label{pa:scorer-factor-ranges}
 \tau((j+n/2)h)<1\quad(j<L),\qquad
 \tau((j+n/2)h)<2\quad(j<K).
\end{equation}
These inequalities are elementary: the Taylor series of $\tau$ has
positive coefficients, so $\tau(1)\ge13/12$; $(2j)!\ge2^j$ gives
$\cosh(1/2)\le8/7$; and $\cosh1<2$.

\begin{lemma}\label{pa:scorer-tail}
For every fixed $n\ge2$ and $S>0$, uniformly for complex $|s|\le S$,
\begin{align}
 \sum_{L\le k<K}
  \left|D_{n,k}B_k^{(n)}(3+h^{2/3}s;h)\right|
 &=O_{n,S}(e^{-1/(384h)}),
 \label{pa:scorer-finite-tail}\\
 S_{1,n,h}(3+h^{2/3}s)
  -\sum_{k<K}D_{n,k}B_k^{(n)}(3+h^{2/3}s;h)
 &=O_{n,S}(e^{-1/(384h)}).
 \label{pa:scorer-actual-tail}
\end{align}
\end{lemma}

\begin{proof}
For $k\le L$, the unperturbed factors are positive. If $t=\delta k$,
then $\log(1+x)\le x$ and $\tau(v)\ge v^2$ give
\begin{align}
 |B_k^{(n)}(3+\delta s;h)|
 &\le \exp\left(S\delta k-
        \sum_{j<k}\tau((j+n/2)h)\right)\nonumber\\
 &\le e^{St-t^3/3}.
 \label{pa:scorer-majorant}
\end{align}
The last step uses $n/2\ge1$, $\sum_{j=1}^k j^2\ge k^3/3$,
and $\delta^3=h^2$.
When $h\le1/4$, one has $L\ge1/(4h)$, so the first $L$ factors
give a loss at least $h^2L^3/3\ge1/(192h)$. The remaining factors
up to $K$ have absolute value at most $1+S\delta$, by
\eqref{pa:scorer-factor-ranges}. Consequently
\begin{equation}\label{pa:scorer-prefix-bound}
 |B_k^{(n)}(3+\delta s;h)|
 \le\exp\left(-\frac1{192h}+2S h^{-1/3}\right)
 \qquad(L\le k\le K).
\end{equation}
Here $K\le2/h$ for small $h$.
The rank factor is at most $C_n(k+1)^{n-2}$ in this range,
by \eqref{pa:rank-factor-products}. Summing proves
\eqref{pa:scorer-finite-tail}, since the polynomial factor and
$e^{2S h^{-1/3}}$ are absorbed by the exponential margin.

For the analytic remainder, apply \eqref{pa:mu-uniform-bound} with
$D_0=1$, $U=1$, and $\eta=S\delta$.
Since $1-1/\tau(1)\ge1/13$, this yields
\[
 \sum_l|\mu_l(1+\delta s)|
 \le\exp\left(\frac{52S}{(e-1)h^{1/3}}\right).
\]
The fixed-rank bound \eqref{pa:nu-bound},
\eqref{pa:scorer-prefix-bound} at $k=K$, and
\eqref{pa:tail-bound} now give \eqref{pa:scorer-actual-tail}.
\end{proof}

We prove the Scorer limit directly at rank two and deduce the higher
ranks from the rank ladder of Theorem~\ref{thm:rank-ladder}; this also
explains why the higher-rank profiles are normalized derivatives of
$\mathcal H_2$. The deduction needs the following estimate for the
divided difference at the shrinking color scale.

\begin{lemma}\label{bridge:scaled-divided-difference}
Put $q=e^{-h/2}$ and $\rho=h^{2/3}$.  For a function $F$ holomorphic
near a compact subset of the $s$-plane, define
\[
 (\mathscr D_hF)(s)
  =\rho\left.\mathcal D_q
       \left[F\left(\frac{X-3}{\rho}\right)\right]
                                                  \right|_{X=3+\rho s}.
\]
For every $R>0$ there are constants $C_R,h_R>0$ such that
\begin{align}
 \|\mathscr D_hF-F'\|_{|s|\leq R}
     &\leq C_R\rho\|F\|_{|s|\leq R+1},
       \label{bridge:scaled-dd-error}\\
 \|\mathscr D_hF\|_{|s|\leq R}
     &\leq C_R\|F\|_{|s|\leq R+1},
       \label{bridge:scaled-dd-bound}
\end{align}
whenever $0<h<h_R$ and $F$ is holomorphic on a neighborhood of
$\{|s|\leq R+1\}$.
\end{lemma}

\begin{proof}
Write $X=3+\rho s$ and choose the holomorphic square root
$v=\sqrt{X^2-4}$ that is positive at $X=3$.  The two color
arguments in \eqref{eq:Dq-color} have center
$X\cosh(h/2)$ and half-distance $v\sinh(h/2)$.  Thus
\[
 (\mathscr D_hF)(s)
   =\frac{F(s+a_h(s)+b_h(s))-F(s+a_h(s)-b_h(s))}
                  {2b_h(s)},
\]
where
\[
 a_h(s)=\frac{X(\cosh(h/2)-1)}{\rho}=O_R(h^{4/3}),
 \qquad
 b_h(s)=\frac{v\sinh(h/2)}{\rho}=O_R(h^{1/3}).
\]
For small $h$, all arguments lie in $|s|<R+1/2$ when $|s|\leq R$.
The quotient is
\[
 \frac12\int_{-1}^1 F'(s+a_h(s)+t b_h(s))\,dt.
\]
Taylor expansion of $F'$ through its linear term, with Cauchy
bounds for its next derivative, gives
\[
 |\mathscr D_hF(s)-F'(s)|
   \leq C_R\bigl(|a_h(s)|+|a_h(s)|^2+|b_h(s)|^2\bigr)
                      \|F\|_{|s|\leq R+1}.
\]
The term linear in $t b_h$ integrates to zero.  Since
$|a_h|=O_R(\rho^2)$ and $|b_h|^2=O_R(\rho)$, this proves
\eqref{bridge:scaled-dd-error}.  The integral formula and
Cauchy's derivative estimate prove \eqref{bridge:scaled-dd-bound}.
\end{proof}

\begin{proposition}\label{bridge:scorer-from-rank-two}
Suppose that the figure-eight sum at rank two satisfies
\begin{equation}\label{bridge:rank-two-scorer-input}
 h^{2/3}S_{1,2,h}(3+h^{2/3}s)
       =\mathcal H_2(s)+O(h^{2/3})
\end{equation}
locally uniformly for complex $s$.  Then for every fixed $n\geq2$,
\begin{equation}\label{bridge:all-rank-scorer-output}
 h^{2(n-1)/3}S_{1,n,h}(3+h^{2/3}s)
       =\mathcal H_n(s)+O_n(h^{2/3})
\end{equation}
locally uniformly for complex $s$.  The same estimate holds after
any fixed number of $s$-derivatives.
\end{proposition}

\begin{proof}
Let $d=n-2$, $\rho=h^{2/3}$, and
$F_h(s)=\rho S_{1,2,h}(3+\rho s)$.  At $z=1$, the exact identity
\eqref{eq:rank-ladder} becomes
\begin{equation}\label{bridge:scaled-rank-ladder}
 \rho^{d+1}S_{1,n,h}(3+\rho s)
    =\frac1{[d]_q!}\mathscr D_h^dF_h(s).
\end{equation}
The hypothesis gives $F_h=\mathcal H_2+O(\rho)$ on every compact
set, and in particular gives locally uniform bounds for $F_h$.
Apply Lemma~\ref{bridge:scaled-divided-difference} successively on
$d+1$ nested disks.  Its operator bound controls the propagated
$O(\rho)$ error, while its difference estimate compares each
application of $\mathscr D_h$ with one derivative.  Induction gives
\[
 \mathscr D_h^dF_h
     =\mathcal H_2^{(d)}+O_d(\rho)
\]
on every fixed compact set.  Since
$[d]_q!=d!\,(1+O_d(h^2))$ and
$\mathcal H_2^{(d)}/d!=\mathcal H_n$, substitution into
\eqref{bridge:scaled-rank-ladder} proves
\eqref{bridge:all-rank-scorer-output}.  Cauchy's formula on a
larger compact set gives the derivative estimates.
\end{proof}

\begin{theorem}\label{pa:scorer}
Fix an integer $n\ge2$. Let $\mathcal H_n$ be defined by
\eqref{pa:scorer-profile}.
For every compact set $C\subset\mathbb C$ and every integer $j\ge0$,
\begin{equation}\label{pa:scorer-rate}
 \sup_{s\in C}\left|
 \frac{d^j}{ds^j}\left[
 h^{2(n-1)/3}S_{1,n,h}(3+h^{2/3}s)-\mathcal H_n(s)
 \right]\right|=O_{n,C,j}(h^{2/3}).
\end{equation}
\end{theorem}

\begin{proof}
By Proposition~\ref{bridge:scorer-from-rank-two}, it suffices to
prove the locally uniform estimate at rank two. In this rank
\(D_{2,k}=1\). Put \(\delta=h^{2/3}\), \(t=\delta k\), and
\(\phi_s(t)=st-t^3/3\), and fix \(|s|\le S\).
For \(k\le L\), all factors of \(B_k^{(2)}(3+\delta s;h)\)
lie in a fixed right half-plane for sufficiently small \(h\).
The Taylor estimate \(\log(1+z)=z+O(|z|^2)\) is uniform for
\(z=\delta s-\tau((j+1)h)\), \(j<L\). Moreover,
\begin{align*}
 \sum_{j<k}\tau((j+1)h)
 &=\frac{h^2k^3}{3}
      +O\bigl(h^2(k+1)^2+h^4(k+2)^5\bigr),\\
 \sum_{j<k}|\delta s-\tau((j+1)h)|^2
 &\le C_S\bigl(k\delta^2+h^4(k+2)^5\bigr).
\end{align*}
Summing the logarithms gives
\begin{equation}\label{pa:scorer-log-rate}
 \left|\log B_k^{(2)}(3+\delta s;h)-\phi_s(t)\right|
       \le C_S\delta(1+t^5),\qquad k\le L.
\end{equation}
The real parts of both logarithms are at most \(St-t^3/3\),
by \eqref{pa:scorer-majorant}. Hence
\eqref{pa:exponential-difference} implies
\begin{equation}\label{pa:scorer-basis-rate}
 \left|B_k^{(2)}(3+\delta s;h)-e^{\phi_s(t)}\right|
 \le C_S\delta(1+t^5)e^{St-t^3/3}.
\end{equation}
After multiplication by the outer Riemann-sum factor \(\delta\),
the summed error is \(O_S(\delta)\). The integrable majorant
and its first derivative have uniformly bounded integrals, so
the Riemann sum of \(e^{\phi_s(t)}\) differs from its integral
by \(O_S(\delta)\). Its tail beyond \(\delta L\) is exponentially
small in \(1/h\). Lemma~\ref{pa:scorer-tail}, at \(n=2\),
controls the remaining finite indices and the analytic remainder.
This proves \eqref{bridge:rank-two-scorer-input}. The exact rank
ladder and Proposition~\ref{bridge:scorer-from-rank-two} now give
\eqref{pa:scorer-rate} at every fixed rank, with all stated
derivative estimates.
\end{proof}

The rank-two function is
$\mathcal H(s)=\mathcal H_2(s)=\pi\operatorname{Hi}(s)$,
where $\operatorname{Hi}$ is the Scorer function. The higher-rank
profiles satisfy
\begin{equation}\label{pa:scorer-rank-profile}
 \mathcal H_n(s)=\frac{\mathcal H^{(n-2)}(s)}{(n-2)!},
 \qquad \mathcal H''-s\mathcal H=1.
\end{equation}
The differential equation follows by integration by parts.
Let $\operatorname{Ai}$ be the entire solution of $f''=sf$ with
\[
 \operatorname{Ai}(0)=\frac1{3^{2/3}\Gamma(2/3)},\qquad
 \operatorname{Ai}'(0)=-\frac1{3^{1/3}\Gamma(1/3)}.
\]
If $\omega=e^{2\pi i/3}$, then
\begin{equation}\label{pa:Airy-profile}
 \omega^{n-1}\mathcal H_n(\omega s)
 -\omega^{2(n-1)}\mathcal H_n(\omega^2s)
 =\frac{2\pi i}{(n-2)!}\operatorname{Ai}^{(n-2)}(s).
\end{equation}
For $n=2$, the left side solves the homogeneous Airy equation.
The Gamma integrals
$\mathcal H(0)=3^{-2/3}\Gamma(1/3)$ and
$\mathcal H'(0)=3^{-1/3}\Gamma(2/3)$, together with the
reflection formula, identify its initial data with those of
$2\pi i\operatorname{Ai}(s)$. Differentiation proves
\eqref{pa:Airy-profile} in every rank. By Theorem~\ref{pa:scorer}, \eqref{pa:Airy-profile} continues to hold,
with an error $O(h^{2/3})$ on compact sets, when $\mathcal H_n(s)$ is
replaced by $h^{2(n-1)/3}S_{1,n,h}(3+h^{2/3}s)$.

\subsection{Integer colors}
\label{pa:node-section}

Put
\begin{equation}\label{pa:node-scale}
 u_m=\operatorname{arccosh}\left(1+\frac1{2m}\right),
 \qquad h_r=\frac{u_m}{r+n/2}.
\end{equation}
The node $x_r^{(n)}$ at $h=h_r$ is exactly $2+1/m$.

\begin{corollary}\label{pa:critical-nodes}
Fix $m\ge2$ and $n\ge2$. As $r\to\infty$,
\begin{equation}\label{pa:Gaussian-node-constant}
 J_r^{(n)}(L_m;e^{-h_r/2})
 =\frac{\Gamma((n-1)/2)}
 {2(n-2)!\alpha_m^{(n-1)/2}}
 h_r^{-(n-1)/2}\bigl(1+O_{m,n}(h_r^{1/2})\bigr).
\end{equation}
For the figure-eight knot, with $u_1=\operatorname{arccosh}(3/2)$,
\begin{equation}\label{pa:Scorer-node-constant}
 J_r^{(n)}(4_1;e^{-h_r/2})
 =\frac{3^{(n-4)/3}\Gamma((n-1)/3)}{(n-2)!}
 h_r^{-2(n-1)/3}\bigl(1+O_n(h_r^{2/3})\bigr).
\end{equation}
\end{corollary}

\begin{proof}
Apply the exact node identity \eqref{pa:node-interpolation} and
Theorems~\ref{pa:gaussian} and \ref{pa:scorer} at $s=0$.
The elementary Gamma integrals give the two positive constants,
so the absolute profile errors are also the stated relative errors.
\end{proof}

At rank two, $N=r+1$ is the representation dimension and the
standard Jones variable is $Q=q^2=e^{-u_m/N}$. In particular,
\eqref{pa:Gaussian-node-constant} becomes
\[
 J_{N-1}^{(2)}(L_m;e^{-u_m/(2N)})
 =\frac{\sqrt\pi}{2\sqrt{\alpha_m u_m}}\sqrt N
       \bigl(1+O_m(N^{-1/2})\bigr),\qquad m\ge2.
\]

\begin{proposition}\label{pa:subcritical-nodes}
Fix $m\ge1$, $n\ge2$, and $0\le\theta_0<1$.
Uniformly for $\theta\in[0,\theta_0]$,
\begin{equation}\label{pa:subcritical-limit}
 J_r^{(n)}\left(L_m;
  \exp\left[-\frac{\theta u_m}{2(r+n/2)}\right]\right)
 \longrightarrow
 \left[1-m\bigl(2\cosh(\theta u_m)-2\bigr)\right]^{-(n-1)}
 \qquad(r\to\infty).
\end{equation}
The expression at $\theta=0$ is understood by Laurent-polynomial
evaluation at $q=1$.
\end{proposition}

\begin{proof}
Let $h=\theta u_m/(r+n/2)$ and
$D_\theta=2\cosh(\theta u_m)-2$.
At the node, every factor of $B_k^{(n)}(x_r^{(n)};h)$ for
$k\le r$ is nonnegative and at most $D_\theta$.
The word representation gives $g_{m,k}^{(n)}\le m^k$, and
\eqref{pa:rank-factor-products} gives
\[
 0\le D_{n,k}g_{m,k}^{(n)}B_k^{(n)}(x_r^{(n)};h)
 \le C_{m,n,\theta_0}\binom{k+n-2}{n-2}\rho^k,
 \qquad \rho=mD_{\theta_0}<1.
\]
For each fixed $k$, the summand tends uniformly in $\theta$ to
$\binom{k+n-2}{n-2}(mD_\theta)^k$.
The displayed summable majorant permits a uniform passage to the
infinite sum. The binomial generating function evaluates it as
$(1-mD_\theta)^{-(n-1)}$.
At $q=1$ all positive-degree Newton terms vanish, so the value is
one, in agreement with the same formula.
\end{proof}

\subsection{Matching with the subcritical regime}

As $s\to-\infty$ the critical profiles decay like powers of $|s|$, and
their leading terms reproduce the inverse-Alexander expression. We now
show that the analytic interpolation joins the two regimes, with
matching first corrections.

\begin{proposition}\label{bridge:profile-left-tails}
For fixed $m\geq2$ and $n\geq2$, as $u\to+\infty$,
\begin{equation}\label{bridge:gaussian-profile-left-tail}
 \mathcal G_{m,n}(-u)
  =\frac1{(mu)^{n-1}}
     \left(1-\frac{\alpha_m n(n-1)}{m^2u^2}
                         +O_{m,n}(u^{-4})\right).
\end{equation}
For fixed $n\geq2$,
\begin{equation}\label{bridge:scorer-profile-left-tail}
 \mathcal H_n(-u)
  =\frac1{u^{n-1}}
     \left(1-\frac{(n-1)n(n+1)}{3u^3}
                         +O_n(u^{-6})\right).
\end{equation}
The bounds for the remainders hold uniformly for $u\geq1$.
\end{proposition}

\begin{proof}
In the Gaussian integral set $v=mut$.  The remaining factor is
$\exp(-\alpha_m v^2/(m^2u^2))$.  For $x\geq0$,
\[
 0\leq e^{-x}-1+x\leq\frac{x^2}{2}.
\]
Integrating this inequality against
$v^{n-2}e^{-v}/(n-2)!$ gives
\eqref{bridge:gaussian-profile-left-tail}, because
$\Gamma(n+1)/\Gamma(n-1)=n(n-1)$.
For the Scorer integral, the substitution $v=ut$ leaves the
factor $\exp(-v^3/(3u^3))$.  The same inequality, with its sixth
moment controlling the remainder, gives
\eqref{bridge:scorer-profile-left-tail}.
\end{proof}

\begin{theorem}\label{bridge:analytic-abelian-overlap}
Fix $m\geq1$ and $n\geq2$, and set
\[
 \eta_0=\frac18,\qquad
 X=2+\frac{1-\eta}{m}.
\]
There are constants $C,h_0>0$, depending only on $m,n$, with
the following properties.

If $m\geq2$, $0<h<h_0$, and $\sqrt h\leq\eta\leq\eta_0$, then
\begin{equation}\label{bridge:gaussian-abelian-overlap}
 \eta^{n-1}S_{m,n,h}(X)
   =1-\alpha_m n(n-1)\frac{h}{\eta^2}
       +O_{m,n}\left(\frac h\eta+\frac{h^2}{\eta^4}\right).
\end{equation}
If $m=1$, $0<h<h_0$, and $h^{2/3}\leq\eta\leq\eta_0$, then
\begin{equation}\label{bridge:scorer-abelian-overlap}
 \eta^{n-1}S_{1,n,h}(3-\eta)
   =1-\frac{(n-1)n(n+1)}3\frac{h^2}{\eta^3}
       +O_n\left(\frac{h^2}{\eta^2}
                           +\frac{h^4}{\eta^6}\right).
\end{equation}
The constants in both error terms are uniform throughout the
stated ranges.
\end{theorem}

\begin{proof}
Put $d=n-2$, $r=1-\eta$, and
$C_{n,k}=\binom{k+d}{d}$.  Choose $U_m>0$ by
\[
 \tau(U_m)=\frac3{2m},\qquad
 \tau(v)=2\cosh v-2,
\]
and set $K=\lceil U_m/h\rceil$.  Decrease $h_0$ so that, for
$h<h_0$,
\begin{equation}\label{bridge:overlap-cutoff}
 \tau((K+n/2)h)>\frac1m,
 \qquad
 m\tau((j+n/2)h)\leq\frac74\quad(0\leq j<K).
\end{equation}
This is possible because the largest argument in the second
inequality tends to $U_m$.  Since $r\geq7/8$, these inequalities
imply
\begin{equation}\label{bridge:overlap-basis-majorant}
 V_{h,k}:=m^kB_k^{(n)}(X;h)
       =\prod_{j<k}\bigl(r-m\tau((j+n/2)h)\bigr),
 \qquad |V_{h,k}|\leq r^k\quad(k\leq K).
\end{equation}
The word representation \eqref{pa:energy-expectation} gives
$0<g_{m,k}^{(n)}/m^k\leq1$.  The rank-factor formula
\eqref{pa:rank-factor-products} also gives
\begin{equation}\label{bridge:overlap-rank-bounds}
 D_{n,k}\leq C C_{n,k},\qquad
 |D_{n,k}-C_{n,k}|
       \leq C h^2(k+1)^2 C_{n,k}\qquad(k\leq K).
\end{equation}
For the second estimate, write the ratio to $C_{n,k}$ as
\[
 \frac{D_{n,k}}{C_{n,k}}
    =\prod_{a=1}^{d}
       \frac{f(h(k+a)/2)}{f(ha/2)},\qquad
 f(v)=\frac{\sinh v}{v},\quad f(0)=1.
\]
All arguments lie in a fixed bounded interval, $f$ is bounded
away from zero there, and $f(v)=1+O(v^2)$.

We first bound the analytic remainder.  Here
$D=X-2=r/m$ is real and nonnegative, and
\eqref{bridge:overlap-cutoff} ensures
$D<\tau((K+n/2)h)$.  The weights $\pi_{m,n,K,i}$ and
$\mu_l(D)$ in \eqref{pa:tail} are probability weights.  The
signed rank correction has uniformly bounded total variation by
\eqref{pa:nu-bound}, since $Kh\geq U_m$.  The resolvent in
\eqref{pa:tail} has absolute value at most one.  Therefore
\begin{equation}\label{bridge:overlap-actual-remainder}
 \left|S_{m,n,h}(X)-
       \sum_{k<K}D_{n,k}g_{m,k}^{(n)}B_k^{(n)}(X;h)\right|
       \leq C(K+1)^d r^K.
\end{equation}

For later use, differentiation of
$\sum_{k\geq0}C_{n,k}r^k=(1-r)^{-(d+1)}$ gives
\begin{align}
 \eta^{d+1}\sum_{k\geq0}kC_{n,k}r^k
    &=(d+1)\frac r\eta,\notag\\
 \eta^{d+1}\sum_{k\geq0}k^2C_{n,k}r^k
    &=(d+1)(d+2)\frac{r^2}{\eta^2}
                            +(d+1)\frac r\eta,
       \label{bridge:negative-binomial-moments}\\
 \eta^{d+1}\sum_{k\geq0}k^3C_{n,k}r^k
    &=(d+1)(d+2)(d+3)\frac{r^3}{\eta^3}
       +3(d+1)(d+2)\frac{r^2}{\eta^2}
       +(d+1)\frac r\eta.\notag
\end{align}
More generally, for each fixed nonnegative integer $j$,
\begin{align}
 \sum_{k\geq0}(k+1)^jC_{n,k}r^k
       &\leq C_j\eta^{-d-j-1},
       \label{bridge:negative-binomial-bound}\\
 \sum_{k\geq K}(k+1)^jC_{n,k}r^k
       &\leq C_j\eta^{-d-j-1}e^{-\eta K/2}.
       \label{bridge:negative-binomial-tail}
\end{align}
The last inequality follows from $r^k\leq e^{-\eta k}$ by
factoring out $e^{-\eta K/2}$ and summing the remaining
polynomial times $e^{-\eta k/2}$.

Suppose first that $m\geq2$.  For a word of length $k$, its
energy in \eqref{pa:word-energy} is nonnegative and at most
$C_{m,n}(k+1)^2$.  Its mean is
$\alpha_mk^2+\beta_{m,n}k$.  The inequality used in the preceding
proposition thus yields the global bound
\begin{equation}\label{bridge:word-taylor-bound}
 \left|\frac{g_{m,k}^{(n)}}{m^k}
          -1+h(\alpha_mk^2+\beta_{m,n}k)\right|
       \leq C h^2(k+1)^4.
\end{equation}
On the finite prefix, telescoping the product in
\eqref{bridge:overlap-basis-majorant} against $r^k$ gives
\begin{equation}\label{bridge:overlap-basis-difference}
 |V_{h,k}-r^k|
  \leq m r^{k-1}\sum_{j<k}\tau((j+n/2)h)
  \leq C h^2(k+1)^3r^k.
\end{equation}
In the last bound we used $r\geq7/8$ and
$\tau(v)\leq C v^2$ on the bounded argument interval.
Combining \eqref{bridge:overlap-rank-bounds},
\eqref{bridge:word-taylor-bound}, and
\eqref{bridge:overlap-basis-difference} gives
\begin{align}
 &\left|D_{n,k}g_{m,k}^{(n)}B_k^{(n)}(X;h)
       -C_{n,k}r^k\{1-h(\alpha_mk^2+\beta_{m,n}k)\}\right|
       \notag\\
 &\hspace{40mm}\leq C h^2(k+1)^4C_{n,k}r^k,
                  \qquad k<K.\label{bridge:gaussian-outer-summand}
\end{align}
The rank-factor error is bounded using
$|V_{h,k}|\leq r^k$ and $g_{m,k}^{(n)}/m^k\leq1$;
it has only the lower power $(k+1)^2$ and is included in the
right side.

By \eqref{bridge:negative-binomial-bound}, the normalized sum of
the right side of \eqref{bridge:gaussian-outer-summand} is
$O(h^2/\eta^4)$.  We can extend the two polynomially weighted
geometric sums to all $k\geq0$.  Their normalized tails are
\[
 O\left((1+h/\eta^2)e^{-c\eta/h}\right)
       =O(e^{-c\eta/h})
\]
when $\eta\geq\sqrt h$.  The normalized analytic remainder in
\eqref{bridge:overlap-actual-remainder} obeys the same bound,
after changing $c>0$: writing $x=\eta/h$ gives
\[
 \eta^{d+1}(K+1)^d r^K
       \leq C\eta x^d e^{-U_mx}
       \leq C e^{-c x}.
\]
Since $x\geq h^{-1/2}$, these errors are $O(h^2)$ and hence
are absorbed by $O(h^2/\eta^4)$.
Formula \eqref{bridge:negative-binomial-moments} now gives the
slightly more precise identity
\begin{align}
 \eta^{n-1}S_{m,n,h}(X)
  =1-h\left\{
       \alpha_m n(n-1)\frac{r^2}{\eta^2}
       +(\alpha_m+\beta_{m,n})(n-1)\frac r\eta\right\}
       +O_{m,n}(h^2/\eta^4).
       \label{bridge:gaussian-outer-first-correction}
\end{align}
Since $r=1-\eta$, this proves
\eqref{bridge:gaussian-abelian-overlap}.

For $m=1$ the word coefficient is exactly one.  On the same
finite prefix, set $x_j=\tau((j+n/2)h)/r$.  Then
$0\leq x_j\leq2$ by \eqref{bridge:overlap-cutoff}.  For numbers
in this interval, telescoping twice gives
\begin{equation}\label{bridge:quadratic-product-remainder}
 \left|\prod_{j<k}(1-x_j)-1+\sum_{j<k}x_j\right|
       \leq\sum_{i<j<k}x_ix_j
       \leq\frac12\left(\sum_{j<k}x_j\right)^2.
\end{equation}
Indeed $|\prod_{i<j}(1-x_i)|\leq1$, and the first telescoping
identity bounds its difference from one by $\sum_{i<j}x_i$.
The Taylor estimate $\tau(v)=v^2+O(v^4)$ is uniform on the
bounded argument interval.  It implies
\begin{align*}
 \sum_{j<k}x_j
   &=\frac{h^2}{r}\sum_{j<k}(j+n/2)^2
                          +O_n(h^4(k+1)^5)\\
   &=\frac{h^2k^3}{3r}
              +O_n\bigl(h^2(k+1)^2+h^4(k+1)^5\bigr),
 \qquad
 \sum_{j<k}x_j\leq C_nh^2(k+1)^3.
\end{align*}
Consequently \eqref{bridge:quadratic-product-remainder} and
\eqref{bridge:overlap-rank-bounds} yield
\begin{align}
 &\left|D_{n,k}B_k^{(n)}(3-\eta;h)
       -C_{n,k}r^k\left(1-\frac{h^2k^3}{3r}\right)\right|
       \notag\\
 &\hspace{20mm}\leq C_n C_{n,k}r^k
       \bigl(h^2(k+1)^2+h^4(k+1)^6\bigr),
                         \qquad k<K.
       \label{bridge:scorer-outer-summand}
\end{align}
The normalized sum of the right side is
$O_n(h^2/\eta^2+h^4/\eta^6)$.  The omitted geometric tails
are bounded, after normalization, by
\[
 C_n(1+h^2/\eta^3)e^{-c\eta/h}.
\]
For $\eta\geq h^{2/3}$ this is $O(e^{-c\eta/h})$.
The normalized analytic remainder has the same estimate as
above.  Now $\eta/h\geq h^{-1/3}$, so these tails are $O(h^4)$
and are absorbed by $O(h^4/\eta^6)$.
Finally, the third moment in
\eqref{bridge:negative-binomial-moments} gives
\[
 \frac{h^2}{3r}\eta^{n-1}
       \sum_{k\geq0}k^3C_{n,k}r^k
    =\frac{(n-1)n(n+1)}3\frac{h^2}{\eta^3}
                         +O_n(h^2/\eta^2).
\]
This proves \eqref{bridge:scorer-abelian-overlap}.
\end{proof}

\begin{corollary}\label{bridge:actual-critical-matching}
Fix $m\geq2$ and $n\geq2$.  Uniformly for $u\geq1$ and
$m\sqrt h\,u\leq\eta_0$, as $h\to0^+$,
\begin{align}
 (m\sqrt h\,u)^{n-1}
       S_{m,n,h}(X_m-\sqrt h\,u)
   &=1-\frac{\alpha_m n(n-1)}{m^2u^2}
              +O_{m,n}(\sqrt h/u+u^{-4}),
       \label{bridge:actual-gaussian-left-tail}\\
 \frac{h^{(n-1)/2}S_{m,n,h}(X_m-\sqrt h\,u)}
                   {\mathcal G_{m,n}(-u)}
   &=1+O_{m,n}(\sqrt h/u+u^{-4}).
       \label{bridge:gaussian-matching-ratio}
\end{align}
For the figure-eight knot, uniformly for $u\geq1$ and
$h^{2/3}u\leq\eta_0$,
\begin{align}
 (h^{2/3}u)^{n-1}S_{1,n,h}(3-h^{2/3}u)
   &=1-\frac{(n-1)n(n+1)}{3u^3}
              +O_n(h^{2/3}/u^2+u^{-6}),
       \label{bridge:actual-scorer-left-tail}\\
 \frac{h^{2(n-1)/3}S_{1,n,h}(3-h^{2/3}u)}
                       {\mathcal H_n(-u)}
   &=1+O_n(h^{2/3}/u^2+u^{-6}).
       \label{bridge:scorer-matching-ratio}
\end{align}
In particular the Gaussian matching, with its displayed first
correction, holds whenever
$u\to\infty$ and $\sqrt h\,u\to0$; the corresponding Scorer
matching holds whenever
$u\to\infty$ and $h^{2/3}u\to0$.
\end{corollary}

\begin{proof}
Substitute $\eta=m\sqrt h\,u$ or $\eta=h^{2/3}u$ into
Theorem~\ref{bridge:analytic-abelian-overlap}.  This gives
\eqref{bridge:actual-gaussian-left-tail} and
\eqref{bridge:actual-scorer-left-tail}.  Proposition
\ref{bridge:profile-left-tails} gives the ratio statements.
Both denominators are positive for real $u>0$; after their
displayed powers of $u$ have been removed they are bounded
away from zero on $u\geq1$, by their integral representations
and their limits at infinity.  In the stated simultaneous limits,
the errors in \eqref{bridge:actual-gaussian-left-tail} and
\eqref{bridge:actual-scorer-left-tail} are respectively
$o(u^{-2})$ and $o(u^{-3})$, so the first corrections agree as
well as the leading terms.
\end{proof}

Since $\eta=1-m(X-2)$, the leading term $\eta^{-(n-1)}$ in
Theorem~\ref{bridge:analytic-abelian-overlap} is the inverse-Alexander
expression $[1-m(X-2)]^{-(n-1)}$ found at the quantum colors in
\eqref{pa:subcritical-limit}. On the real axis the analytic
interpolation therefore passes from this expression to the critical
profiles, and the two descriptions overlap.

\section{Uniform critical asymptotics with growing rank}
\label{pa:growing-rank}

We now let the rank $n=n(h)$ depend on $h$, keep $m$ fixed, and
determine how far the fixed-rank limits of
Section~\ref{a:fixed-critical} remain valid. Each statement below
involves a combined parameter, $hn^3$ in the Gaussian case and
$h^{2/3}n^{5/3}$ in the Scorer case, and its $o(1)$ terms are uniform as
long as this parameter stays below a fixed bound. Put
\begin{align}
 V_m&=\frac{(m-1)(m+1)(16m^2+1)}{180m^2},\label{pa:variance-constant}\\
 c_m&=-\frac{(4m+1)(3m^3+3m^2+8m-2)}
 {180m^2(2\alpha_m)^{3/2}}<0\qquad(m\ge2).
 \label{pa:rank-correction-constant}
\end{align}
All limits below are uniform on compact subsets of the complex
$s$-plane. On such a set the profiles in the denominators have no zeros
once $h$ is small (see \eqref{pa:centered-gaussian-rank} and its
analogue in the proof of Theorem~\ref{pa:growing-rank-scorer}).

\begin{theorem}\label{pa:growing-rank-gaussian}
Fix $m\ge2$. Let $n=n(h)\to\infty$ as $h\to0^+$, with $hn^3$
bounded. Then
\begin{equation}\label{pa:growing-rank-gaussian-limit}
 \frac{h^{(n-1)/2}S_{m,n,h}(X_m+\sqrt h\,s)}
 {\mathcal G_{m,n}(s)}
 -\exp\{c_m\sqrt h\,n^{3/2}\}\longrightarrow0
\end{equation}
locally uniformly in $s$.
\end{theorem}

\begin{theorem}\label{pa:growing-rank-scorer}
For $m=1$, suppose $n\to\infty$ and $h^{2/3}n^{5/3}$ remains
bounded.  Then
\begin{equation}\label{pa:growing-rank-scorer-limit}
 \frac{h^{2(n-1)/3}S_{1,n,h}(3+h^{2/3}s)}
 {\mathcal H_n(s)}
 -\exp\{-3h^{2/3}n^{5/3}/40\}\longrightarrow0
\end{equation}
locally uniformly in $s$.
\end{theorem}

Together with the fixed-rank theorems, these results give relative
limits equal to one when $n=o(h^{-1/3})$ in the Gaussian case and
$n=o(h^{-2/5})$ in the Scorer case. (If $n$ stays bounded along a subsequence, the compact set must avoid
the zeros of the corresponding fixed-rank profiles.) At $n\sim a h^{-1/3}$ and
$n\sim a h^{-2/5}$, respectively, the limiting ratios are
$\exp(c_ma^{3/2})$ and $\exp(-3a^{5/3}/40)$ for $a>0$.
Both profiles are positive at $s=0$, so the correction is visible
already at the center of the critical window, and these scales are
sharp.

\subsection{Cumulants of the word energy}
Let $W_1,\ldots,W_k$ be independent uniform letters in
$\{1,\ldots,m\}$, and set
\[
 b(a)=m-a,\qquad
 K(a,b)=\boldsymbol1_{a>b}+2(m-\max(a,b)),\qquad
 a_m=(m-1)/2.
\]
The word energy of \eqref{pa:word-energy} is exactly
\begin{equation}\label{pa:pair-energy}
 E_{m,n}=\sum_{i<j}K(W_i,W_j)+n\sum_i b(W_i).
\end{equation}
In particular
$\mathbb E E_{m,n}=\alpha_m k^2+(a_m n-\alpha_m)k$.
Changing one letter changes the energy by at most $C_m(k+n)$.
The bounded-difference inequality therefore gives, for real $u$,
\begin{equation}\label{pa:growing-bd}
 \log\mathbb E e^{u(E-\mathbb E E)}
 \le C_m u^2k(k+n)^2.
\end{equation}

\begin{lemma}\label{pa:rank-cumulants}
For $k\ge2$ and $h(k+n)$ sufficiently small, depending only on
$m$, one has
\begin{align}
 \operatorname{Var}E_{m,n}
 &=V_m k^3+O_m((n+1)k^2+n^2k),\label{pa:rank-variance}\\
 \log\frac{g_{m,k}^{(n)}}{m^k}
 &=-\alpha_mhk^2-(a_m n-\alpha_m)hk
   +\tfrac12V_mh^2k^3\nonumber\\
 &\quad+O_m\bigl(h^2((n+1)k^2+n^2k)
                    +h^3k(k+n)^3\bigr).
 \label{pa:rank-logg}
\end{align}
\end{lemma}

\begin{proof}
The first-order projections of $K$ are
\[
 r(a)=m-1-(a-1)^2/m-2\alpha_m,\qquad
 c(a)=m-a^2/m-2\alpha_m.
\]
Both have mean zero.  Writing
$K(a,b)=2\alpha_m+r(a)+c(b)+K_0(a,b)$ gives a kernel $K_0$
with both conditional means zero.  Terms with different unordered
index pairs in $\sum_{i<j}K_0(W_i,W_j)$ are orthogonal, even when
they share one index, and that sum is orthogonal to all one-letter
functions.  Its variance is therefore $O_m(k^2)$.  The linear
projection of the full centered energy is
\[
 \sum_{i=1}^k\{(k-i)r(W_i)+(i-1)c(W_i)
                         +n(b(W_i)-a_m)\}.
\]
Its variance has leading coefficient
\[
 \frac{1}{3m}\sum_{a=1}^m\{r(a)^2+c(a)^2+r(a)c(a)\}=V_m,
\]
by the sums of powers through degree four.  The remaining terms
are bounded by $C_m((n+1)k^2+n^2k)$, proving
\eqref{pa:rank-variance}.

To bound the remainder in the cumulant expansion,
distribute the one-letter term of \eqref{pa:pair-energy} equally
among the pairs, and write
\[
 E=\sum_{i<j}\left\{K(W_i,W_j)
       +\frac{n}{k-1}(b(W_i)+b(W_j))\right\}.
\]
There are at most $k^2/2$ bounded variables, each of size at most
$C_m(1+n/k)$, with a dependency graph of maximum degree at most
$2k$.  The dependency-graph cumulant bound of
F\'eray--M\'eliot--Nikeghbali \cite[Theorems~9.6--9.7]{FMNCumulants}
therefore gives
\[
 |\kappa_j(E)|\le j^{j-2}C_m^j k(k+n)^j\qquad(j\ge2).
\]
Since $j^{j-2}/j!\le e^j$, its Taylor series, divided by $j!$,
converges absolutely for $|h|(k+n)<c_m'$, and the terms $j\ge3$
sum to $O_m(h^3k(k+n)^3)$.  The exponential of this convergent
series agrees near zero with the entire moment generating
function; the identity theorem then justifies its use throughout
that disk.  Retaining the mean and variance, and using the exact
word representation, proves \eqref{pa:rank-logg}.
\end{proof}

\subsection{The Gaussian correction}
\begin{proof}[Proof of Theorem~\ref{pa:growing-rank-gaussian}]
Put
\[
 d=n-2,\quad k_0=\sqrt{n/(2\alpha_mh)},\quad
 \eta=\sqrt h\,n^{3/2},\quad t=\sqrt h\,k.
\]
For each fixed $0<\epsilon<M<\infty$, the following expansions
hold uniformly for $\epsilon\le k/k_0\le M$ and $s$ in a compact
set:
\begin{align}
 \log\frac{d!D_{n,k}}{k^d}
 &=\frac{n^2}{2k}+o(1),\label{pa:grow-D-gauss}\\
 \log\bigl(m^kB_k^{(n)}(X_m+\sqrt h\,s;h)\bigr)
 &=mst-\frac m3h^2k^3+o(1),\label{pa:grow-B-gauss}\\
 \log\frac{g_{m,k}^{(n)}}{m^k}
 &=-\alpha_mhk^2-a_mnhk+\tfrac12V_mh^2k^3+o(1).
 \label{pa:grow-g-gauss}
\end{align}
For each Newton factor we use the branch of the logarithm near
one; all its arguments lie in that neighborhood uniformly on
these indices. To verify the first expansion, use
\eqref{pa:rank-factor-products} and $\log(\sinh x/x)=x^2/6+O(x^4)$:
\begin{align}
 \log\frac{d!D_{n,k}}{k^d}
 &=\frac{d(d+1)}{2k}
   +O(d^3/k^2)
   +\frac{h^2}{24}\{dk^2+d(d+1)k\}\nonumber\\
 &\quad+O(h^4d(k+d)^4).
 \label{pa:uniform-rank-expansion}
\end{align}
All errors after replacing $d(d+1)$ by $n^2$ tend to zero under
$hn^3=O(1)$.  The second expansion follows by Taylor expansion of
$\log(1+m\sqrt h\,s-m\tau((j+n/2)h))$ and summation:
$\sum_{j<k}(j+n/2)^2=k^3/3+O(nk^2+n^2k+k^2)$.
The omitted terms are bounded by a constant times
\[
 hk+h^2(nk^2+n^2k)+h^4(k+n)^5
            +\sqrt h\,h^2(k+n)^3,
\]
which tends to zero in this range.  Lemma~\ref{pa:rank-cumulants}
gives the third expansion; in particular
$h^3k^4=O(hn^2)=o(1)$.
Consequently the exact summand divided by
$k^d e^{mst-\alpha_mhk^2}/d!$ has logarithm
\begin{equation}\label{pa:gaussian-correction-function}
 \eta\left\{\frac{\sqrt{2\alpha_m}}{2y}
       -\frac{a_my}{\sqrt{2\alpha_m}}
       +\frac{(V_m/2-m/3)y^3}{(2\alpha_m)^{3/2}}\right\}+o(1),
 \qquad y=k/k_0.
\end{equation}
At $y=1$ the expression in braces is exactly $c_m$.

We next compare the sum with its normalization. Set
\[
 Z_0=\frac1{d!}\int_0^\infty k^d e^{-\alpha_mhk^2}\,dk
     =h^{-(d+1)/2}\mathcal G_{m,n}(0),\qquad
 t_d=\sqrt{d/(2\alpha_m)}.
\]
The Gamma integral and Stirling's formula give
\begin{equation}\label{pa:gaussian-rank-normalization}
 Z_0\asymp_m\frac{k_0^{d+1}e^{-n/2}}{d!\sqrt n}.
\end{equation}
Fix $S>0$ and restrict to $|s|\le S$. All tail estimates below
are relative to $Z_0e^{m\Re(s)t_d}$.

For $hk\le a_0$, with $a_0>0$ sufficiently small and $h$ sufficiently
small, the elementary product bounds give
\[
 D_{n,k}\le\frac{(k+d)^d}{d!}e^{hdk/2},\qquad
 m^k|B_k|\le e^{mS\sqrt h k}.
\]
Since $n/k_0\to0$, for $k\le\epsilon k_0$ one has
$k+d\le2\epsilon k_0$ eventually. Also $hnk_0=O(\eta)$ and
$\sqrt h k_0=O(\sqrt n)$. Using $g/m^k\le1$ and counting at
most $Ck_0$ terms, \eqref{pa:gaussian-rank-normalization} bounds
this part of the sum by
\[
 C_S\sqrt n(2\epsilon)^d
       \exp\{n/2+C_S\sqrt n+C\eta\}.
\]
It tends to zero exponentially in $n$ when $\epsilon$ is a fixed
sufficiently small positive number.

For $k\ge n$ with $hk\le a_0$, the exact mean and
\eqref{pa:growing-bd}, after decreasing $a_0$ if necessary, give
\[
 g_{m,k}^{(n)}/m^k\le e^{-\alpha_mhk^2/2}.
\]
If $k\ge Mk_0$, then $d^2/k=O(\eta)$, while
$n/k\to0$ and $1/(\sqrt h k)\to0$ uniformly. Thus
$D_{n,k}\le k^d e^{d^2/k+hdk/2}/d!$ and the preceding basis
bound imply
\[
 |D_{n,k}g_{m,k}^{(n)}B_k|
 \le\frac{C}{d!}k^d e^{-\alpha_mhk^2/4}.
\]
For a sufficiently large fixed $M$, the function
$y^d e^{-ny^2/8}$ decreases for $y\ge M-1$. With mesh $1/k_0$,
\[
 \frac1{k_0}\sum_{k\ge Mk_0}
      (k/k_0)^d e^{-n(k/k_0)^2/8}
 \le\int_{M-1}^\infty y^d e^{-ny^2/8}\,dy.
\]
After division by the normalization, the upper tail is at most
$C\sqrt n e^{n/2+C_S\sqrt n}$ times this integral. Choosing
$M$ so that $\log y-y^2/8\le-y^2/16$ for $y\ge M-1$ proves
an exponentially small bound. In particular the number of lattice
points introduces no factor depending on $1/h$ after normalization.

On $\epsilon\le y\le M$, the correction in
\eqref{pa:gaussian-correction-function} is bounded, whereas
$n^{-1}(d\log y-ny^2/2)\to\log y-y^2/2$ uniformly.
The latter function has its unique maximum at $y=1$. Outside
any fixed neighborhood of one, its strict maximum gap gives an
$e^{-c n}$ bound after multiplication by
$\sqrt n e^{C_S\sqrt n}$. Inside that neighborhood the expression
in braces in \eqref{pa:gaussian-correction-function} differs from
$c_m$ by a quantity tending uniformly to zero as the neighborhood
shrinks. The boundedness of $\eta$ therefore allows this
correction to be replaced by $\eta c_m$ in the normalized sum.

It remains to compare the complex profiles with their lattice sums,
which we do by means of a global bound. Put
$\Phi_d(t)=d\log t-\alpha_mt^2$. Its maximum is at $t_d$, and
\[
 \Phi_d''(t)=-d/t^2-2\alpha_m\le-2\alpha_m,
 \qquad
 \Phi_d(t_d+u)-\Phi_d(t_d)\le-\alpha_mu^2.
\]
On bounded $u$-intervals the difference in the last expression
converges to $-2\alpha_mu^2$. Hence dominated convergence, with
majorant $e^{-\alpha_mu^2+mS|u|}$, gives
\begin{equation}\label{pa:centered-gaussian-rank}
 \frac{\mathcal G_{m,n}(s)}
      {\mathcal G_{m,n}(0)e^{mst_d}}
 \longrightarrow \exp\{m^2s^2/(8\alpha_m)\}
\end{equation}
locally uniformly for complex $s$. The shifted lattice
$u=\sqrt h k-t_d$ has mesh $\sqrt h\to0$. On bounded intervals
the same Taylor convergence proves convergence of its Riemann
sums. Outside such intervals the displayed Gaussian majorant
bounds both the integrals and the mesh-weighted sums uniformly.
Thus the lattice has the same limit, and its relative error
against the profile is $o(1)$. The limit in
\eqref{pa:centered-gaussian-rank} is nonzero, so the profiles
are eventually zero-free on every fixed compact set. It also
follows that the preceding absolute tail and correction estimates
give relative estimates for complex $s$.

Finally use $K=\lceil U_m/h\rceil$ as in
Lemma~\ref{pa:gaussian-tail}.  The statements of the exact tail
formula and its bounds for the coefficient sums hold for every integer $n$.
Since $hn=o(1)$, the bound $|1-m\tau((j+n/2)h)|\le1$ for $j<K$ used in
that proof still holds. The logarithm of the rank factor $D_{n,K}$ is
$O(n\log(1/h+n))=o(1/h)$, the bound \eqref{pa:nu-bound} for the signed
weights $\nu_r$ has logarithm $O(n)$, and \eqref{pa:mu-uniform-bound}
contributes $O(h^{-1/2})$.
Thus the analytic remainder is $O(e^{-c/h})$ for some fixed $c>0$.
For $a_0/h\le k<K$, the global coefficient estimate of
Lemma~\ref{pa:energy} supplies an $e^{-c'/h}$ bound uniformly in
$n$. The same rank and basis estimates, summed over $O(1/h)$
indices, preserve this exponential bound. Finally,
\[
 |\log Z_0|\le C_m n(1+|\log h|+\log n)=o(1/h),
 \qquad |mst_d|=O_S(\sqrt n)=o(1/h).
\]
Thus all remaining terms and the analytic remainder are
negligible relative to the profile. This proves
\eqref{pa:growing-rank-gaussian-limit}. For the smaller-rank
consequence stated after the theorems, bounded-rank subsequences
are covered by Theorem~\ref{pa:gaussian} on zero-free compact sets.
\end{proof}

\subsection{The Scorer correction}
\begin{proof}[Proof of Theorem~\ref{pa:growing-rank-scorer}]
For $m=1$ write
\[
 k_0=(n/h^2)^{1/3},\quad \eta=h^{2/3}n^{5/3},\quad
 y=k/k_0,\quad t=h^{2/3}k.
\]
Assume $n\to\infty$ and $\eta=O(1)$. On a fixed interval
$0<\epsilon\le y\le M$, expansion
\eqref{pa:uniform-rank-expansion} gives
\[
 \log\frac{d!D_{n,k}}{k^d}
 =\frac{n^2}{2k}+\frac{h^2nk^2}{24}+o(1).
\]
The Taylor formula
$\log(1-\tau(z))=-z^2-7z^4/12+O(z^6)$ gives
\[
 \log B_k^{(n)}(3+h^{2/3}s;h)
 =st-\frac{h^2k^3}{3}-\frac{h^2nk^2}{2}
                     -\frac{7h^4k^5}{60}+o(1).
\]
The cross terms involving bounded $s$ tend to zero.  For example,
$h^{2/3}h^2k^3=O(h^{2/3}n)=o(1)$.
The remaining error terms are bounded by sums of
$h^{4/3}n^{7/3}$, $h^{2/3}n^{2/3}$ and
$h^{4/3}n^{4/3}$, all tending to zero under $\eta=O(1)$.
The exact summand divided by $k^d e^{st-h^2k^3/3}/d!$ therefore
has logarithm
\begin{equation}\label{pa:scorer-correction-function}
 \eta\left\{\frac1{2y}-\frac{11y^2}{24}
                         -\frac{7y^5}{60}\right\}+o(1).
\end{equation}
At $y=1$ the expression in braces equals $-3/40$. Indeed, write
$-11y^2/24=-y^2/2+y^2/24$: the terms $1/(2y)$ and $-y^2/2$ cancel at
$y=1$, and $-3/40=1/24-7/60$ comes from the hyperbolic correction to the
rank factor and from the quartic term of the Newton product.

Set
\[
 Z_0=\frac1{d!}\int_0^\infty k^d e^{-h^2k^3/3}\,dk
 =h^{-2(d+1)/3}\mathcal H_n(0),\qquad t_d=d^{1/3}.
\]
The Gamma integral and Stirling's formula give
\[
 Z_0\asymp\frac{k_0^{d+1}e^{-n/3}}{d!\sqrt n}.
\]
Fix $S>0$ and normalize absolute estimates by
$Z_0e^{\Re(s)t_d}$ for $|s|\le S$.
Here $n/k_0\to0$ and $hnk_0=\sqrt\eta\sqrt n$.
For $k\le\epsilon k_0$, the elementary rank bound used above and
$|B_k|\le e^{Sh^{2/3}k}$ therefore give a normalized bound
\[
 C_S\sqrt n(2\epsilon)^d
             \exp\{n/3+C_S\sqrt n\}.
\]
It tends to zero exponentially for a fixed sufficiently small
$\epsilon>0$.

For $hk\le a_0$ with $a_0>0$ sufficiently small, every
unperturbed Newton factor is positive. The inequality
$\tau(v)\ge v^2$ gives the exact bound
$|B_k|\le e^{St-h^2k^3/3}$. Also
\[
 \log\frac{d!D_{n,k}}{k^d}
 \le d^2/k+C h^2d(k+d)^2\qquad(k\ge n).
\]
For $k\ge Mk_0$, the first term is $O(\eta)$ and the second is
$o(h^2k^3)$, because $n/k\to0$ uniformly. The tilt $St$ is
also $o(h^2k^3)$ on this range. Consequently
\[
 |D_{n,k}B_k|\le\frac C{d!}k^d e^{-h^2k^3/6}.
\]
For a sufficiently large fixed $M$, the decreasing upper envelope
and its lattice mesh $1/k_0$ give
\[
 \frac1{k_0}\sum_{k\ge Mk_0}
        (k/k_0)^d e^{-n(k/k_0)^3/6}
 \le\int_{M-1}^\infty y^d e^{-ny^3/6}\,dy.
\]
After normalization this is multiplied by at most
$C\sqrt n e^{n/3+C_Sn^{1/3}}$. Since
$\log y-y^3/6\le-y^3/12$ for $y$ large, the result is
exponentially small in $n$ for sufficiently large $M$.
On $\epsilon\le y\le M$, the leading exponent divided by $n$
converges uniformly to $\log y-y^3/3$, whose unique maximum is
at one. The strict maximum gap outside any fixed neighborhood
of one, followed by shrinking that neighborhood, permits the
correction in \eqref{pa:scorer-correction-function} to be
replaced by $-3\eta/40$ in the normalized sum.

To justify the complex comparison and the lattice normalization,
put $\Phi_d(t)=d\log t-t^3/3$ and
$\sigma_d=(3t_d)^{-1/2}$. Its maximum is at $t_d$, and, for all
$t>0$,
\[
 \Phi_d''(t)=-d/t^2-2t\le-3t_d.
\]
Indeed $y^{-2}+2y\ge3$ for $y=t/t_d>0$. Thus
\[
 \Phi_d(t_d+\sigma_du)-\Phi_d(t_d)\le-u^2/2,
\]
and on bounded $u$-intervals the left side tends to $-u^2/2$.
Since $\sigma_d\to0$, dominated convergence, with the uniform
bound $e^{-u^2/2+S\sigma_d|u|}$, proves
\[
 \frac{\mathcal H_n(s)}{\mathcal H_n(0)e^{st_d}}\longrightarrow1
\]
locally uniformly for complex $s$. In the $u$ coordinate the
lattice mesh is
$h^{2/3}/\sigma_d=\sqrt3 h^{2/3}d^{1/6}
=O(\eta n^{-3/2})\to0$.
Uniform convergence on bounded intervals and the same Gaussian
tail bound give the identical limit for the lattice sums. This
proves their relative error is $o(1)$, eventual absence of profile
zeros on compact sets, and the complex relative versions of the
preceding tail and correction estimates.

For $a_0/h\le k<K=\lceil1/h\rceil$, the fixed prefix of
$\lfloor a_0/h\rfloor$ positive factors supplies an
$e^{-c/h}$ loss; subsequent factors have modulus at most
$1+Sh^{2/3}$. The rank factor has logarithm
$O(n\log(1/h+n))=o(1/h)$, and the bounds for the coefficient sums contribute
$O(n+h^{-1/3})=o(1/h)$ to the logarithm. The exact resolvent
remainder at $K$ is therefore $O(e^{-c'/h})$ for some $c'>0$,
and the same holds for the sum over these $O(1/h)$ indices.
Finally $|\log Z_0|=O(n(1+|\log h|+\log n))=o(1/h)$ and
$|st_d|=O_S(n^{1/3})=o(1/h)$. These tails are negligible
relative to the profile, completing the proof.
\end{proof}

\begin{corollary}\label{pa:growing-rank-nodes}
Set $h=u_m/(r+n/2)$ and let $r\to\infty$.  For fixed $m\ge2$,
if $hn^3$ is bounded and $n\to\infty$, then
\[
 J_r^{(n)}(L_m;e^{-h/2})
 =\frac{\Gamma((n-1)/2)}{2(n-2)!\alpha_m^{(n-1)/2}}
 h^{-(n-1)/2}\exp\{c_m\sqrt h\,n^{3/2}\}(1+o(1)).
\]
For $m=1$, if $h^{2/3}n^{5/3}$ is bounded and $n\to\infty$, then
\[
 J_r^{(n)}(4_1;e^{-h/2})
 =\frac{3^{(n-4)/3}\Gamma((n-1)/3)}{(n-2)!}
 h^{-2(n-1)/3}\exp\{-3h^{2/3}n^{5/3}/40\}(1+o(1)).
\]
The relative factor is one in the respective smaller-rank ranges.
\end{corollary}

\begin{proof}
The choice of $h$ makes $x_r^{(n)}=2\cosh u_m=X_m$ exactly.
Apply the two theorems at $s=0$, evaluate the elementary Gamma
integrals, and use the exact node-interpolation formula.
\end{proof}

\section{Logarithmic monodromy and geometric periods}
\label{pa:complex}

In this section we continue the sums from real $h>0$ to the half-plane
$\Re h>0$ and study the monodromy of the continuation under
$h\mapsto h+2\pi i$. Fix integers $m\ge1$ and $n\ge2$, and put $d=n-2$.
The coordinate
\[
 Y=Q^{n/2}X
\]
removes the half-integral powers of $Q$ from the Borel coefficients.
For $0<|Q|<1$, define
\begin{equation}\label{pa:complex-polynomial}
 P_{n,k}(Q,Y)=\prod_{j=0}^{k-1}
                (1-YQ^j+Q^{n+2j})
\end{equation}
and the germ
\begin{equation}\label{pa:complex-borel}
 \widehat H_{m,n}(Q,Y,t)
 =\sum_{k\ge0}g_{m,k}^{(n)}(Q)
 \frac{(Q^{n-1};Q)_k P_{n,k}(Q,Y)}{(Q;Q)_k^2}
                    (-Q^{-d}t)^k.
\end{equation}
Here the $g_{m,k}^{(n)}$ are the polynomial coefficients defined
in Section~\ref{pa:data}. Thus, when $Q=e^{-h}$ with $h>0$,
\[
 \widehat H_{m,n}(Q,Q^{n/2}X,t)=H_{m,n,h}(X,t).
\]

\begin{lemma}\label{pa:complex-continuation}
The germ \eqref{pa:complex-borel} has a continuation that is
jointly holomorphic off its polar divisor for $0<|Q|<1$ and
$Y\in\mathbb C$, and meromorphic in $t$. Its possible poles
are simple and lie at
\[
 t_j=-Q^{d-j},\qquad j\ge0.
\]
For compact sets $C\Subset\{\Re h>0\}$ and $E\Subset\mathbb C$,
there are constants $K,C_0$ such that
\[
 |\widehat H_{m,n}(e^{-h},Y,e^{hx})|
 \le C_0(1+|e^{hx}|^K),
 \qquad h\in C,\ Y\in E,\ x\in\mathbb R.
\]
\end{lemma}

\begin{proof}
Write $\rho=|Q|$. The coefficients of $g_{m,k}^{(n)}$ are
nonnegative, and the positive moment formula gives
\[
 |g_{m,k}^{(n)}(Q)|\le g_{m,k}^{(n)}(\rho)
 \le(\rho^n;\rho)_\infty^{-1}.
\]
Consequently
\begin{equation}\label{pa:complex-symbol}
 G_{m,n}(Q,w)=\sum_{r\ge0}
 \frac{Q^{r(r+n-1)}(w;Q^{-1})_r}{(Q;Q)_r}
             g_{m-1,r}^{(n)}(Q),\qquad m\ge2,
\end{equation}
is entire in $w$ and normally convergent on compact subsets of
$0<|Q|<1$, jointly with $w$. Set $G_{1,n}=1$.
The finite recursion implies
$G_{m,n}(Q,Q^k)=g_{m,k}^{(n)}(Q)$.

Let $a+b=Y$ and $ab=Q^n$, so that
$P_{n,k}(Q,Y)=(a,b;Q)_k$. On a compact parameter set choose $K$
such that $|aQ^K|,|bQ^K|<1$ uniformly, and put
\[
 b_k=\frac{(Q^{n-1},a,b;Q)_k}{(Q;Q)_k^2},\qquad k\ge0.
\]
The function
\begin{align*}
 A_K(w)={}&G_{m,n}(Q,Q^Kw)
 \frac{(Q^{K+n-1},aQ^K,bQ^K;Q)_\infty}
      {(Q^{K+1};Q)_\infty^2}\\
 &\hspace{7mm}\times
 \frac{(Q^{K+1}w;Q)_\infty^2}
      {(Q^{K+n-1}w,aQ^Kw,bQ^Kw;Q)_\infty}
 =\sum_{j\ge0}a_{K,j}w^j
\end{align*}
is analytic on a disk of radius greater than one, with a locally
uniform bound for its coefficient norm. Cancellation of products
gives
\begin{equation}\label{pa:complex-fractions}
 \widehat H_{m,n}
 =\sum_{k<K}g_{m,k}^{(n)}b_k(-Q^{-d}t)^k
  +b_K(-Q^{-d}t)^K
       \sum_{j\ge0}\frac{a_{K,j}}{1+tQ^{j-d}}.
\end{equation}
The expression is symmetric in $a,b$, hence holomorphic in $Y$
also when the two roots coincide. It agrees with the germ near
$t=0$, and establishes the stated meromorphic continuation.

For $h$ in a compact subset of the right half-plane,
\begin{equation}\label{pa:spiral-separation}
 |1+e^{h\tau}|\ge c_C(1+e^{(\Re h)\tau}),
 \qquad \tau\in\mathbb R.
\end{equation}
The ratio tends uniformly to one at both ends of the real line.
On a bounded interval it has a positive minimum: an equality
$e^{h\tau}=-1$ would force $\tau=0$ by comparison of absolute values, and
$e^0\ne-1$. Applying \eqref{pa:spiral-separation} with
$\tau=x+d-j$ to \eqref{pa:complex-fractions} proves the bound.
\end{proof}

For $\Re h>0$, orient the logarithmic spiral
\[
 \Gamma_h=\{e^{hx}:x\in\mathbb R\}
\]
by increasing $x$, and define
\begin{equation}\label{pa:spiral-sum}
 \widehat S_{m,n}(h,Y)
 =\frac1h\int_{\Gamma_h}
       \frac{\widehat H_{m,n}(e^{-h},Y,t)}
            {(-t;e^{-h})_\infty}\,dt.
\end{equation}

\begin{theorem}\label{pa:spiral-theorem}
The integral \eqref{pa:spiral-sum} converges normally and is jointly
holomorphic for $\Re h>0$ and $Y\in\mathbb C$. It satisfies
\[
 \widehat S_{m,n}(h,Q^{n/2}X)=S_{m,n,h}(X),\qquad h>0,
\]
and its kernel has moments
\begin{equation}\label{pa:complex-moments}
 \frac1h\int_{\Gamma_h}
          \frac{t^k}{(-t;Q)_\infty}\,dt
       =(Q;Q)_kQ^{-k(k+1)/2},\qquad k\ge0.
\end{equation}
At the moving nodes
\begin{equation}\label{pa:complex-nodes}
 Y_{n,r}(Q)=Q^{-r}+Q^{n+r},\qquad r\ge0,
\end{equation}
it equals the same finite colored knot invariant as on the
positive real axis.
\end{theorem}

\begin{proof}
Parameterizing the integral by $t=e^{hx}$ gives the integrand
\[
 \frac{e^{hx}\widehat H_{m,n}(e^{-h},Y,e^{hx})}
      {(-e^{hx};e^{-h})_\infty}.
\]
On a compact $h$ set write $a_0=\min\Re h>0$.
For $x\ge0$, let $L=\lfloor x\rfloor$. The first $L+1$ factors
of the kernel denominator, bounded by
\eqref{pa:spiral-separation}, give
\[
 |(-e^{hx};e^{-h})_\infty|^{-1}
 \le C\exp(-a_0x^2/2+C_1x).
\]
The remaining product is nonzero and continuous in $h$ and
$x-L\in[0,1]$, and is uniformly bounded away from zero.
For $x\le0$ the reciprocal kernel is bounded. Lemma
\ref{pa:complex-continuation} therefore gives integrable bounds
at both ends, locally uniformly in $(h,Y)$.

The same reasoning with $\widehat H_{m,n}$ replaced by $t^k$
shows that the left side of \eqref{pa:complex-moments} is
holomorphic in $h$. The real moment identity of
Theorem~\ref{thm:kernel} extends to the connected right half-plane.
At \eqref{pa:complex-nodes}, the polynomial
$P_{n,k}$ vanishes for $k>r$, so the Borel transform terminates.
Termwise application of the moment identity proves the final
assertion.
\end{proof}

The spiral cannot be replaced by the positive ray, even for $h$ in a
thin sector around the positive axis:

\begin{proposition}\label{pa:fixed-ray}
For no $h_0,\theta_0>0$ does the undeformed positive-ray integral
exist at every $h$ with $0<|h|<h_0$ and $|\arg h|<\theta_0$,
even when the Borel transform is identically one.
\end{proposition}

\begin{proof}
Choose $0<\theta<\min(\theta_0,\pi/2)$ and set
$h_N=(\pi\cot\theta+i\pi)/N$. The kernel has a simple pole
$-e^{Nh_N}=e^{\pi\cot\theta}$ on the positive ray, while
$h_N\to0$ along the direction $\theta$.
\end{proof}

\subsection{A residue formula on the logarithmic cover}

Put $t_l=-Q^{-l}$ for $l\ge-d$, and write
\[
 \widehat H_{m,n}(Q,Y,t)
   =\frac{a_l(Q,Y)}{t-t_l}+b_l(Q,Y)+O(t-t_l).
\]
For $l\ge0$ define
\begin{align}
 \rho_l&=\frac{(-1)^lQ^{l(l-1)/2}}
                  {(Q;Q)_l(Q;Q)_\infty},
 &\sigma_l&=Q^l\left(l-\sum_{j>l}\frac{Q^j}{1-Q^j}\right),
 \label{pa:complex-kernel-residues}\\
 \alpha_l&=\rho_l a_l,
 &\beta_l&=\rho_l(b_l+\sigma_la_l).
 \label{pa:complex-double-residues}
\end{align}
For $-d\le l<0$ set
\begin{equation}\label{pa:complex-finite-residues}
 \alpha_l=0,\qquad
 \beta_l=\frac{a_l}{(Q^{-l};Q)_\infty}.
\end{equation}
Thus $\alpha_l$ and $\beta_l$ are the coefficients of orders
$-2$ and $-1$ in the Laurent expansion of
$\widehat H_{m,n}(Q,Y,t)/(-t;Q)_\infty$.

\begin{theorem}\label{pa:monodromy}
The functions
\begin{equation}\label{pa:monodromy-blocks}
 \mathcal A_{m,n}(Q,Y)=-\sum_{l=-d}^\infty l\beta_l(Q,Y),
 \qquad
 \mathcal B_{m,n}(Q,Y)=\sum_{l=0}^\infty Q^l\alpha_l(Q,Y)
\end{equation}
are single-valued and holomorphic for $0<|Q|<1$, $Y\in\mathbb C$.
Both series converge normally on compact subsets, and
\begin{equation}\label{pa:affine-log}
 \widehat S_{m,n}(h,Y)
   =\mathcal A_{m,n}(e^{-h},Y)
       +\frac{\mathcal B_{m,n}(e^{-h},Y)}h.
\end{equation}
In particular,
\begin{equation}\label{pa:deck-law}
 \widehat S_{m,n}(h+2\pi i,Y)-\widehat S_{m,n}(h,Y)
 =-\frac{2\pi i}{h(h+2\pi i)}\mathcal B_{m,n}(e^{-h},Y).
\end{equation}
Moreover there is a holomorphic function $\mathcal C_{m,n}$,
entire in $Y$, such that
\begin{equation}\label{pa:monodromy-divisor}
 \mathcal B_{m,n}(Q,Y)
 =\prod_{j\ge0}(1-YQ^j+Q^{n+2j})\,
                       \mathcal C_{m,n}(Q,Y).
\end{equation}
\end{theorem}

\begin{proof}
Removing the vanishing factor in the kernel product gives
\eqref{pa:complex-kernel-residues}. The logarithmic derivative of
its regular factor at $t_l$ is
\[
 -\sum_{j\ne l}\frac{Q^j}{1-Q^{j-l}}
 =Q^l\left(l-\sum_{j>l}\frac{Q^j}{1-Q^j}\right).
\]
For $l<0$ the kernel is regular at $t_l$, giving
\eqref{pa:complex-finite-residues}.

On a compact $Q$ annulus, \eqref{pa:complex-fractions} gives
$|a_l|+|b_l|\le C A^l$ for $l\ge0$, uniformly on compact
$Y$ sets. To see this, isolate the singular fraction, evaluate
the others at $t_l$, and use
$|1-Q^s|\ge1-|Q|^s$ for $s>0$; for $s<0$ factor out
$Q^s$. The polynomial prefix and its derivative contribute
only exponential factors in $l$. The Gaussian decay of
$\rho_l$ then proves normal convergence in
\eqref{pa:monodromy-blocks}.

For real $0<Q<1$, put
$f(t)=\widehat H_{m,n}(Q,Y,t)/(-t;Q)_\infty$ and apply the
residue theorem to $f(t)\log(-t)$, with the logarithmic cut on
the positive axis. Take outer radii $Q^{-L-1/2}$. On these
circles the Borel transform is polynomially bounded by
\eqref{pa:complex-fractions}, while the reciprocal kernel
decays faster than every power. The outer integrals vanish;
the inner circle at zero also vanishes. Hence
\[
 \int_0^\infty f(t)\,dt
 =-\sum_{l=-d}^\infty
       \operatorname{Res}_{t=t_l}\bigl(f(t)\log(-t)\bigr).
\]
At $t_l$, $\log(-t_l)=lh$ and $(\log(-t))'=1/t_l=-Q^l$.
The right side is
$-h\sum l\beta_l+\sum_{l\ge0}Q^l\alpha_l$.
This proves \eqref{pa:affine-log} for real $h>0$.
Holomorphy extends it to the right half-plane, and
single-valuedness of the two blocks in $Q$ proves
\eqref{pa:deck-law}.

At a node \eqref{pa:complex-nodes}, the Borel transform is a
polynomial, so all $a_l$ vanish and $\mathcal B_{m,n}=0$.
The product in \eqref{pa:monodromy-divisor} is a normally
convergent entire function of $Y$. Its zeros are precisely the
nodes, and they are simple: equality of two distinct nodes
would imply either $Q^{r-s}=1$ or $Q^{n+r+s}=1$.
Both are impossible for $0<|Q|<1$. Division by these simple
zeros proves the factorization, jointly holomorphically in $Q$.
\end{proof}

\begin{corollary}\label{pa:trace-deck}
Let
\[
 S^{\mathrm{tr}}_{m,n}(h,X)
    =\widehat S_{m,n}(h,e^{-nh/2}X),\qquad
 \epsilon_n=(-1)^n.
\]
Then
\begin{align}
 S^{\mathrm{tr}}_{m,n}(h+2\pi i,X)
 ={}&S^{\mathrm{tr}}_{m,n}(h,\epsilon_nX)\notag\\
 &-\frac{2\pi i}{h(h+2\pi i)}
     \mathcal B_{m,n}(e^{-h},\epsilon_ne^{-nh/2}X).
 \label{pa:trace-deck-law}
\end{align}
\end{corollary}

\begin{proof}
Under $h\mapsto h+2\pi i$, the variable $Q$ is unchanged and
$e^{-nh/2}X$ is multiplied by $\epsilon_n$. Apply
\eqref{pa:deck-law} at that value of $Y$.
\end{proof}

The sign $\epsilon_n$ comes from the half-integral power $e^{-nh/2}$ in
the coordinate $Y$. At the moving nodes
$Y_{n,r}(Q)$ the logarithmic monodromy vanishes. For $m=1,n=2$,
our Borel transform and real integral, after $Y=QX$, are the
figure-eight expressions in \cite[(3.4)--(3.5)]{KMM}.
The corresponding logarithmic decomposition is already given
in \cite[(3.39)--(3.43)]{KMM}.

\subsection{Stokes residues and the logarithmic coefficient}

The fixed-$Q$ Stokes discontinuity and the logarithmic monodromy in
$h$ are obtained from the same Borel residues through two explicit
linear pairings.

For $d=n-2$ and $0<|Q|<1$, put
\begin{equation}\label{bridge:residue-coordinates}
 \mathcal R_{m,n,j}(Q,Y)
   =\operatorname{Res}_{\xi=-Q^{d-j}}
                   \widehat H_{m,n}(Q,Y,\xi),
 \qquad
 b_{m,n,j}(Q,Y)=Q^{j-d}\mathcal R_{m,n,j}(Q,Y),
 \qquad j\geq0.
\end{equation}
(The residue is zero when $-Q^{d-j}$ is not a pole.) When $0<Q<1$,
write
$h=-\log Q$ and $w_Q(t)=(-t;Q)_\infty^{-1}$, and define
\[
 \mathscr S_{m,n,Q}(Y,z)
   =\frac1h\int_0^\infty
              \widehat H_{m,n}(Q,Y,zt)w_Q(t)\,dt,
 \qquad z\in\mathbb C\setminus(-\infty,0].
\]
Thus $\mathscr S_{m,n,Q}(Y,1)=\widehat S_{m,n}(h,Y)$ and, in
the trace coordinate, this is the same sum as
$S_{m,n,h}(Q^{-n/2}Y,z)$.  Let
\[
 \Delta_{m,n,Q}(Y,r)
  =\mathscr S_{m,n,Q}(Y,-r+i0)
        -\mathscr S_{m,n,Q}(Y,-r-i0),\qquad r>0.
\]

\begin{proposition}\label{bridge:residues-stokes-logarithmic-block}
For $0<|Q|<1$ and every $Y\in\mathbb C$,
\begin{align}
 \mathcal B_{m,n}(Q,Y)
 &=\frac1{(Q;Q)_\infty}
   \sum_{l\geq0}
     \frac{(-1)^lQ^{l(l-1)/2}}{(Q;Q)_l}
                    b_{m,n,l+d}(Q,Y).\label{bridge:block-from-residues}
\end{align}
The series converges normally on compact subsets of
$\{0<|Q|<1\}\times\mathbb C_Y$.

If $0<Q<1$, then, with $\tau_j=Q^{d-j}/r$,
\begin{equation}\label{bridge:jump-from-residues}
 \Delta_{m,n,Q}(Y,r)
   =-\frac{2\pi i}{h}
       \sum_{j\geq0}b_{m,n,j}(Q,Y)\tau_jw_Q(\tau_j).
\end{equation}
This series is locally normally convergent in $Y$ and $r>0$.
For fixed $Q$ and $Y$, the function
$r\mapsto\Delta_{m,n,Q}(Y,r)$ determines every coefficient
$b_{m,n,j}(Q,Y)$ and hence determines
$\mathcal B_{m,n}(Q,Y)$.  Adding a finite polynomial to the
Borel transform, with coefficients single-valued and holomorphic
in $Q,Y$, changes neither of these two quantities.
\end{proposition}

\begin{proof}
The finite-prefix representation \eqref{pa:complex-fractions}
shows that, on each compact parameter set, there are constants
$C,D>0$ such that
\begin{equation}\label{bridge:residue-exponential-bound}
 |b_{m,n,j}(Q,Y)|\leq CD^j\qquad(j\geq0).
\end{equation}
Indeed, choose the prefix length uniformly on that compact set.
The coefficients of its analytic tail have a common Cauchy bound
on a disk of radius greater than one.  At the pole
$\xi=-Q^{d-j}$, multiplication by the fixed power of $\xi$
in the prefix representation contributes at most an additional
exponential factor in $j$.  The polynomial prefix has no residue.
The residue at each fixed pole is jointly holomorphic in $Q,Y$,
as follows either from that representation or from a local Cauchy
integral around the pole.  Distinct poles do not collide when
$0<|Q|<1$.

At $\xi=-Q^{-l}$, $l\geq0$, the residue of $w_Q$ is
\[
 \rho_l=\frac{(-1)^lQ^{l(l-1)/2}}
                 {(Q;Q)_l(Q;Q)_\infty}.
\]
The coefficient of the double pole in
$\widehat H_{m,n}(Q,Y,\xi)w_Q(\xi)$ is consequently
$\rho_l\mathcal R_{m,n,l+d}(Q,Y)$.  Formula
\eqref{pa:monodromy-blocks} gives
\[
 \mathcal B_{m,n}(Q,Y)
    =\sum_{l\geq0}Q^l\rho_l
                        \mathcal R_{m,n,l+d}(Q,Y)
    =\sum_{l\geq0}\rho_l b_{m,n,l+d}(Q,Y),
\]
which is \eqref{bridge:block-from-residues}.  On a compact
$Q$ annulus, the factors $(Q;Q)_l$ are bounded away from zero
uniformly in $l$, and $Q^{l(l-1)/2}$ has uniform quadratic
exponential decay.  Together with
\eqref{bridge:residue-exponential-bound}, this proves the asserted
normal convergence for all colors.

For real $Q$, apply the lateral residue formula
\eqref{gen:eq:jumpres} to the same finite-prefix representation.
Since the pole is $-Q^{d-j}$ and $z=-r$, its contribution is
\[
 -\frac{2\pi i}{hr}\mathcal R_{m,n,j}(Q,Y)
                  w_Q(Q^{d-j}/r)
 =-\frac{2\pi i}{h}b_{m,n,j}(Q,Y)
                  \tau_jw_Q(\tau_j).
\]
The kernel has quadratic exponential decay along this geometric
sequence, so \eqref{bridge:residue-exponential-bound} also proves
local normal convergence of \eqref{bridge:jump-from-residues}.

The coefficients $b_{m,n,j}$ may grow exponentially, so
Theorem~\ref{gen:thm:injective}, which assumes $b\in\W_R$, does not
apply directly. The exponential bound
\eqref{bridge:residue-exponential-bound} suffices, however. Put
$f(t)=t w_Q(t)$.  Its logarithmic derivative is
\[
 \frac{d\log f(t)}{d\log t}
   =1-\sum_{a\geq0}\frac{tQ^a}{1+tQ^a}
       \longrightarrow-\infty\qquad(t\to\infty).
\]
For every $\gamma>0$, there is therefore a $T_\gamma$ such that
\begin{equation}\label{bridge:kernel-scale-bound}
 \frac{f(ct)}{f(t)}\leq c^{-\gamma}
 \qquad(t\geq T_\gamma,\ c\geq1).
\end{equation}
For each fixed $c>1$ the same ratio tends to zero.
If $j_0$ is the first nonzero coefficient, divide
\eqref{bridge:jump-from-residues} by $f(Q^{d-j_0}/r)$.
Choose $\gamma$ such that $DQ^\gamma<1$ in
\eqref{bridge:residue-exponential-bound}.  The terms with $j>j_0$
are then bounded by a summable geometric sequence, independently
of sufficiently small $r$, and each tends to zero.  Hence
\[
 \Delta_{m,n,Q}(Y,r)
   \sim-\frac{2\pi i}{h}b_{m,n,j_0}(Q,Y)
                      f(Q^{d-j_0}/r),\qquad r\downarrow0.
\]
Equivalently, all coefficients, including possible zero ones,
can be recovered successively by
\begin{equation}\label{bridge:residue-recovery}
 b_{m,n,j}(Q,Y)
 =-\frac{h}{2\pi i}\lim_{r\downarrow0}
 \frac{\displaystyle
  \Delta_{m,n,Q}(Y,r)
   +\frac{2\pi i}{h}\sum_{i<j}b_{m,n,i}(Q,Y)
                         f(Q^{d-i}/r)}
      {f(Q^{d-j}/r)}.
\end{equation}
Substitution into \eqref{bridge:block-from-residues} determines
the logarithmic coefficient.  Finally, the Laplace transform of
a polynomial is a polynomial in $z$, by the exact moment
identity.  It has no lateral discontinuity.  When its coefficients
are single-valued and holomorphic in $Q,Y$, the same is true of
its transformed coefficients, since the moments are single-valued
in $Q$.  It therefore contributes nothing to $\mathcal B_{m,n}$, in agreement
with \eqref{bridge:block-from-residues}, since a polynomial has no
residues.
\end{proof}

Thus the Stokes data in the auxiliary variable $z$, at fixed $Q$,
determine the coefficient of the deck transformation \eqref{pa:deck-law}
in $h$. Both are linear in the residues $b_{m,n,j}$, with different
weights: the jump weighs them by the kernel values $\tau_jw_Q(\tau_j)$,
which depend on $r$, while $\mathcal B_{m,n}$ is the single pairing
\eqref{bridge:block-from-residues}.

\subsection{The contour representation}
\label{r4:period-bridge}

We now write the coefficient $\mathcal B_{m,n}$ of
Theorem~\ref{pa:monodromy} as a contour integral, and then estimate it
for the figure-eight knot. As before, $m\geq1$, $n\geq2$, $d=n-2$ and
$0<|Q|<1$, and $a,b$ are determined by $a+b=Y$ and $ab=Q^n$. Recall
that
\[
 \widehat S_{m,n}(h,Y)=\mathcal A_{m,n}(Q,Y)
                         +h^{-1}\mathcal B_{m,n}(Q,Y),\qquad Q=e^{-h}.
\]
We use the entire coefficient interpolant $G_{m,n}$ from
\eqref{pa:complex-symbol}; it satisfies
$G_{m,n}(Q,Q^k)=g^{(n)}_{m,k}(Q)$, and $G_{1,n}=1$.

\begin{theorem}\label{r4:closed-block}
For any positively oriented circle $|w|=R$ with
\[
 0<R<\min\{|a|^{-1},|b|^{-1},|Q|^{1-n}\},
\]
one has the exact identity
\begin{align}
 \mathcal B_{m,n}(Q,Y)
 ={}&\frac{(Q^{n-1},a,b;Q)_\infty}{(Q;Q)_\infty^3}
       \frac1{2\pi i}\int_{|w|=R}w^{-d}G_{m,n}(Q,w)
 \frac{(Qw;Q)_\infty^2(w^{-1};Q)_\infty}
      {(Q^{n-1}w,aw,bw;Q)_\infty}\,\frac{dw}{w}.
 \label{r4:closed-block-formula}
\end{align}
The integral is independent of the admissible radius.  The expression is
symmetric in $a,b$ and extends jointly holomorphically in $Q,Y$ throughout
$0<|Q|<1$, $Y\in\mathbb C$, including the zeros of its displayed prefactor.
\end{theorem}

\begin{proof}
First suppose $|a|,|b|<1$.  On a disk of radius greater than one put
\[
 A(w)=G_{m,n}(Q,w)
 \frac{(Q^{n-1},a,b;Q)_\infty}{(Q;Q)_\infty^2}
 \frac{(Qw;Q)_\infty^2}{(Q^{n-1}w,aw,bw;Q)_\infty}
 =\sum_{j\geq0}A_jw^j.
\]
The original Borel transform is
\[
 \widehat H_{m,n}(Q,Y,t)=\sum_{j\geq0}
                  \frac{A_j}{1+tQ^{j-d}}.
\]
At $t_l=-Q^{-l}$, $l\geq0$, its residue is
$A_{l+d}Q^{-l}$.  The residue of the Laplace kernel at the same point is
\[
 \rho_l=\frac{(-1)^lQ^{l(l-1)/2}}
                {(Q;Q)_l(Q;Q)_\infty}.
\]
Consequently the double-pole coefficient in the product is
$\alpha_l=\rho_l A_{l+d}Q^{-l}$.  By \eqref{pa:monodromy-blocks},
\[
 \mathcal B_{m,n}=\sum_{l\geq0}Q^l\alpha_l
   =\frac1{(Q;Q)_\infty}
      \sum_{l\geq0}\frac{(-1)^lQ^{l(l-1)/2}}{(Q;Q)_l}A_{l+d}.
\]
Euler's expansion of $(w^{-1};Q)_\infty$ identifies this absolutely
convergent coefficient pairing with the constant term of
$w^{-d}A(w)(w^{-1};Q)_\infty$, proving
\eqref{r4:closed-block-formula}.  Euler's series converges normally on
every circle avoiding $w=0$, so the initial radius need not exceed one.

For arbitrary $Y$, choose a sufficiently small positive $R$ on a compact
parameter neighborhood.  Every denominator in the integrand is then
nonzero and its products converge normally; the same is true of $G$.
The integral is jointly holomorphic and symmetric in the roots $a,b$.
It is therefore holomorphic in $Y$, including where $a=b$.  Analytic
continuation from the nonempty domain $|a|,|b|<1$ proves the formula
everywhere.  The annulus between two admissible circles contains no
poles, which proves independence of the radius.  Since no division by $(Q^{n-1},a,b;Q)_\infty$ occurs, the formula also
holds on the zero set of this product.
\end{proof}

\subsection{The contour integral at \texorpdfstring{$X=2$}{X=2}}

For the figure-eight case $m=1$, specialize the trace to $X=2$ and set
\[
 B_n(h)=\mathcal B_{1,n}(e^{-h},2e^{-nh/2}),\qquad \Re h>0.
\]
Thus $Y=2Q^{n/2}$, so that $a=b=Q^{n/2}$ (in particular $a=b=Q$ at
rank two), and the circle $|w|=1$ is admissible in
\eqref{r4:closed-block-formula}. Inverting $w=z^{-1}$ gives
\begin{align}
 B_n(h)&=p_n(h)\frac1{2\pi i}\int_{|z|=1}F_{n,h}(z)\frac{dz}{z},
 \label{r4:rank-period}\\
 p_n(h)&=\frac{(Q^{n-1};Q)_\infty(Q^{n/2};Q)_\infty^2}
                    {(Q;Q)_\infty^3},\notag\\
 F_{n,h}(z)&=z^d\frac{(Q/z;Q)_\infty^2(z;Q)_\infty}
              {(Q^{n-1}/z;Q)_\infty(Q^{n/2}/z;Q)_\infty^2}.
 \notag
\end{align}
The inversion reverses the orientation of the circle, and the sign in
$dw/w=-dz/z$ compensates for this. The zeros of the denominator of
$F_{n,h}$ lie strictly inside the unit circle, so $F_{n,h}$ is
holomorphic on the circle; it vanishes at $z=1$.

For $n=2$ the prefactor is one and Euler's two product expansions give
\begin{equation}
 B_2(h)=g(Q)=\sum_{k\geq0}
      \frac{(-1)^kQ^{k(k+1)/2}}{(Q;Q)_k^2}
   =\frac1{2\pi i}\int_{|z|=1}
      \frac{(z;Q)_\infty}{(Q/z;Q)_\infty}\frac{dz}{z}.
 \label{r4:g-contour}
\end{equation}
The identity $B_2=g$ and the role of this $q$-series for the
figure-eight knot are known; see \cite[(3.47)--(3.48)]{KMM} and
\cite[Corollary~1.7]{GKState}. For $n>2$, formula \eqref{r4:rank-period}
multiplies the integrand by an explicit factor, and we now determine the
asymptotics of the resulting integral.

Write $\operatorname{Li}_2(z)=\sum_{k\ge1}z^k/k^2$ for $|z|\le1$,
with its principal analytic continuation when needed, and
$\operatorname{Cl}_2(\theta)=\sum_{k\ge1}\sin(k\theta)/k^2
=\operatorname{Im}\operatorname{Li}_2(e^{i\theta})$.
Let
\[
 V=2\,\operatorname{Im}\operatorname{Li}_2(e^{i\pi/3})
       =\operatorname{Vol}(S^3\setminus4_1),\qquad
 b_d=\frac{27d^2+11}{72\sqrt3}.
\]
That $V$ is the hyperbolic volume of the figure-eight knot complement
is classical; for this normalization of the critical value see
\cite[(2), \S4.1]{FantiniWheeler} and \cite[(3.70)--(3.72)]{KMM}.
Powers and square roots of $h$ below take their principal values in
$\Re h>0$.

\begin{theorem}\label{r4:sector-period}
Fix $n\geq2$ and $0<\delta<\pi/4$.  Uniformly as $h\to0$ in
\[
 \delta\leq\arg h\leq\frac\pi2-\delta,
\]
the contour integral \eqref{r4:rank-period} satisfies
\begin{equation}
 B_n(h)=\frac{(-1)^d h^{1/2-2d}}
                {d!\,\Gamma(n/2)^2\sqrt{2\pi}\,3^{1/4}}
       e^{iV/h-i\pi/4}
       \left(1-i b_dh+O_{n,\delta}(|h|^2)\right).
 \label{r4:upper-period}
\end{equation}
Uniformly in the lower sector
$-\pi/2+\delta\leq\arg h\leq-\delta$, one instead has
\begin{equation}
 B_n(h)=\frac{(-1)^d h^{1/2-2d}}
                {d!\,\Gamma(n/2)^2\sqrt{2\pi}\,3^{1/4}}
       e^{-iV/h+i\pi/4}
       \left(1+i b_dh+O_{n,\delta}(|h|^2)\right).
 \label{r4:lower-period}
\end{equation}

\end{theorem}

\begin{proof}
Put $z=e^{i\theta}$ with $0\leq\theta\leq2\pi$ and define
\[
 W(\theta)=\operatorname{Li}_2(e^{-i\theta})
                 -\operatorname{Li}_2(e^{i\theta})
              =-2i\operatorname{Cl}_2(\theta).
\]
Near $z=1$ the asymptotic expansion of the Pochhammer symbols is not
uniform, so we first show that a neighborhood of $z=1$ does not
contribute to the leading order.
For $\theta\ne0,2\pi$, expand all products except the explicit factor
$1-z$ in $(z;Q)_\infty$.  The coefficients of $z^{-k}$ and $z^k$ in
$\log F_{n,h}(z)-d\log z-\log(1-z)$ are respectively
\[
 -\frac{2e^{-kh}-e^{-(n-1)kh}-2e^{-nkh/2}}
           {k(1-e^{-kh})},\qquad
 -\frac{e^{-kh}}{k(1-e^{-kh})}.
\]
After subtracting $W(\theta)/h$, the two coefficient numerators divided
by $1-e^{-kh}$ are bounded for $|kh|\leq1$.  For $|kh|\geq1$, their
differences from the subtracted terms are bounded by
$C(e^{-c|kh|}+|kh|^{-1})$, uniformly in the closed sector.
Splitting the sum at $k=\lfloor|h|^{-1}\rfloor$ therefore gives
\begin{equation}
 |F_{n,h}(e^{i\theta})|
       \leq C|h|^{-C}\exp\bigl(\Re(W(\theta)/h)\bigr),
 \qquad 0\leq\theta\leq2\pi.
 \label{r4:global-contour-bound}
\end{equation}
At $z=1$ the left side is zero, so the same bound holds there by
continuity.  

Since
\[
 W'(\theta)=2i\log(2\sin(\theta/2)),\qquad
 W''(\theta)=i\cot(\theta/2),
\]
the two critical points are $\theta=\pi/3$ and $\theta=5\pi/3$, with
values $-iV$ and $iV$.  In the upper sector $\Re(W(\theta)/h)$ has
a unique global maximum at $\theta_*=5\pi/3$.  Outside any fixed
neighborhood of $\theta_*$ there is a gap at least $c/|h|$, uniformly
in that sector.  It follows from \eqref{r4:global-contour-bound} that
this entire complement, including $z=1$, contributes at most a
polynomial in $|h|^{-1}$ times
$e^{\Re(iV/h)-c/|h|}$.  Relative to the main term it is therefore $O(|h|^N)$ for every $N$.

On a fixed neighborhood of $\theta_*$ the Euler--Maclaurin expansion
of each logarithmic product holds uniformly, together with its
derivatives. To obtain the remainder, apply the Euler--Maclaurin formula
to
$\log(1-ze^{-s})$ on the ray $s\in h\mathbb R_{\geq0}$; the ray lies
in a fixed closed cone in the right half-plane, and $z$ remains away
from one.  Rotating that ray to the positive real axis crosses no
singularity.  This yields, for fixed $c\geq0$,
\[
 \log(Q^cz;Q)_\infty
 =-h^{-1}\operatorname{Li}_2(z)
  +(\tfrac12-c)\log(1-z)
  -\tfrac h2(c^2-c+1/6)\frac z{1-z}+O(|h|^2),
\]
together with the analogous expansions to higher order and their
derivatives.
Also, with $p_1=d(3d+4)/8$,
\begin{equation}
 p_n(h)=\frac{h^{-2d}}{d!\,\Gamma(n/2)^2}
                   \bigl(1+p_1h+O(|h|^2)\bigr).
 \label{r4:rank-prefactor}
\end{equation}
For even $n=2r$, this follows directly from the exact identity
$p_n=[(Q;Q)_d(Q;Q)_{r-1}^2]^{-1}$.  For odd $n=2r+1$ use
\[
 p_n=\frac{(Q^{1/2};Q)_\infty^2/(Q;Q)_\infty^2}
                  {(Q;Q)_{2r-1}(Q^{1/2};Q)_r^2}
\]
and the eta transformation, which gives
\[
 \frac{(Q^{1/2};Q)_\infty^2}{(Q;Q)_\infty^2}
   =\frac h\pi e^{-h/8}
      \frac{(\widetilde Q^2;\widetilde Q^2)_\infty^2}
           {(\widetilde Q;\widetilde Q)_\infty^4},
 \qquad \widetilde Q=e^{-4\pi^2/h}.
\]
The quotient of dual products is $1+O(e^{-c/|h|})$ in the stated
sector.  This proves \eqref{r4:rank-prefactor}, with a remainder that is
uniform in the sector.

After multiplying \eqref{r4:rank-prefactor} into the local integrand,
the leading angular amplitude, apart from the fixed normalization,
is
\[
 a_d(\theta)=(-1)^d(2\sin(\theta/2))^{2d+1}.
\]
Its first relative correction is
$-i(p_1+1/12)\cot(\theta/2)$.
At $\theta_*$ one has
\begin{gather*}
 a_d(\theta_*)=(-1)^d,\quad
 \frac{a'_d}{a_d}(\theta_*)=-\frac{(2d+1)\sqrt3}{2},\quad
 \frac{a''_d}{a_d}(\theta_*)=\frac{(2d+1)(6d-1)}4,\\
 W''(\theta_*)=-i\sqrt3,\qquad
 W'''(\theta_*)=-2i,\qquad
 W''''(\theta_*)=-2i\sqrt3.
\end{gather*}
The real angular contour is a uniformly decaying Gaussian direction
in the upper sector, since
$\Re(-i\sqrt3/h)\leq-c/|h|$.  The Gaussian factor with the positive
orientation of the original circle is
\[
 \frac1{2\pi}\sqrt{\frac{2\pi h}{i\sqrt3}}
        =\frac{\sqrt h}{\sqrt{2\pi}\,3^{1/4}}e^{-i\pi/4}.
\]
The Gaussian moments give the relative order-$h$ coefficient
\[
 -\frac{a''_d}{2a_dW''}
 +\frac{a'_dW'''}{2a_d(W'')^2}
 +\frac{W''''}{8(W'')^2}
 -\frac{5(W''')^2}{24(W'')^3}
 +i\sqrt3(p_1+1/12)
 =-\frac{i(27d^2+11)}{72\sqrt3}.
\]
Taylor expansion with two further orders and Gaussian domination
gives the stated relative $O(|h|^2)$ error.  The exponentially small
complement was already estimated above.  Complex conjugation of the
exact products proves the lower-sector formula.
\end{proof}

At rank two, the action and the first two coefficients agree with the
geometric perturbative periods in \cite[(3.70)--(3.72)]{KMM};
Theorem~\ref{r4:sector-period} adds a remainder estimate that is uniform
in the sectors. For $n>2$ the integral is the rank-two one with the
explicit additional factors of \eqref{r4:rank-period}.

In Appendix~\ref{app:torsion} we check that the rank-two leading factor
equals \(\sqrt{-h/(4\pi\tau_\mu)}\), where $\tau_\mu$ is the
meridian-normalized adjoint torsion and the square root is determined
by the orientation of the contour.

\begin{corollary}\label{r4:rank-raising}
In the figure-eight case, the integrands and prefactors satisfy
\begin{align*}
 \frac{p_{n+2}(h)}{p_n(h)}
 &=\frac1{(1-Q^{n-1})(1-Q^n)(1-Q^{n/2})^2},\\
 \frac{F_{n+2,h}(z)}{F_{n,h}(z)}
 &=z^2(1-Q^{n-1}/z)(1-Q^n/z)(1-Q^{n/2}/z)^2.
\end{align*}
Thus raising the rank by two multiplies the integrand of the same
contour integral by a Laurent polynomial in $z$; in particular the
critical values of the classical phase do not change.
\end{corollary}

\begin{proof}
Cancel the shifted infinite products in \eqref{r4:rank-period}.
As $h\to0$ the inserted factor tends to $z^{-2}(z-1)^4$, which does not
vanish at the critical points $z=e^{\pm i\pi/3}$, so neither saddle point
is removed.
\end{proof}

\subsection{The real expansion at \texorpdfstring{$X=2$}{X=2}}

On the positive real $h$-axis the behavior at $X=2$ is much simpler:
the remainder formula of Theorem~\ref{pa:source} gives a complete
asymptotic expansion of $S_{m,n,h}(2)$ in powers of $h$.

\begin{proposition}\label{r4:parabolic-abelian}
For fixed $m\geq1$, $n\geq2$, and every fixed integer $L\geq1$,
as $h\to0^+$ one has
\begin{equation}
 S_{m,n,h}(2)-\sum_{k=0}^{L-1}
       D_{n,k}(h)g_{m,k}^{(n)}(e^{-h})B_k^{(n)}(2;h)
          =O_{m,n,L}(h^{2L}).
 \label{r4:parabolic-remainder}
\end{equation}
In particular $S_{m,n,h}(2)$ has a full Poincar\'e expansion at this
trace, whose first nonconstant terms are
\begin{equation}
 S_{m,n,h}(2)=1-\frac{m(n-1)n^2}{4}h^2
       +\frac{m(m-1)(n-1)n^3}{8}h^3+O_{m,n}(h^4).
 \label{r4:parabolic-first}
\end{equation}
\end{proposition}

\begin{proof}
Write $\tau(v)=2\cosh v-2$ and choose $a>0$ so small that
$m\tau(a)<1$.  Take $K=\lfloor a/h\rfloor$.  For sufficiently
small $h$ there is a fixed $r<1$ such that
\[
 m\tau((K+n/2)h)\leq r,
 \qquad Kh\geq a/2.
\]
The word model gives $0<g_{m,k}^{(n)}\leq m^k$, and the exact
sinh product for $D_{n,k}$ implies
$D_{n,k}\leq C(k+1)^d$ for $k\leq K$.  At $X=2$,
\[
 |B_k^{(n)}(2;h)|=\prod_{j=0}^{k-1}\tau((j+n/2)h).
\]
If $L\leq k\leq K$, the first $L$ factors are bounded by
$C_Lh^{2L}$ and every remaining factor by $r/m$.  It follows that
\[
 \sum_{k=L}^{K-1}D_{n,k}g_{m,k}^{(n)}
                         |B_k^{(n)}(2;h)|
 \leq C_Lh^{2L}\sum_{k\geq L}(k+1)^d r^{k-L}
 =O(h^{2L}).
\]
The exact analytic remainder \eqref{pa:tail} at $D=X-2=0$ has
$\sum_l|\mu_l(0)|=1$ and
$\sum_r|\nu_r|\leq(1-e^{-a/2})^{-2(n+1)}$ by
\eqref{pa:nu-bound}.  Hence \eqref{pa:tail-bound} bounds the
remaining tail by $C(K+1)^dr^K=O(e^{-c/h})$.
This proves \eqref{r4:parabolic-remainder}.

Every term in the fixed finite sum is analytic at $h=0$ and its
$k$th Newton factor vanishes to order $2k$.  Thus the estimates for
all fixed $L$ give the asserted full expansion.  For the displayed
coefficients, use
\[
 D_{n,1}=(n-1)+O(h^2),\quad
 g_{m,1}^{(n)}=\sum_{j=0}^{m-1}e^{-njh}
       =m-\tfrac12nm(m-1)h+O(h^2),
\]
and $B_1^{(n)}(2;h)=-n^2h^2/4+O(h^4)$.
All $k\geq2$ terms contribute $O(h^4)$ by the proved remainder.
\end{proof}

\appendix

\section{Boundary values and subcritical growth}
\label{app:boundary}

Fix $0<p<1$ and write $w_p(t)=(-t;p)_\infty^{-1}$.
We first control Cauchy integrals and their Euler derivatives uniformly
as their poles approach the positive axis. We then treat meromorphic
Borel transforms with growth below the quadratic decay rate of the
Laplace kernel.

\begin{lemma}\label{lem:lateral-tail}
For an integer \(K\ge0\) let \(f_K(t)=t^Kw_p(t)\), and for
$\tau\in\mathbb C\setminus[0,\infty)$ put
\[
 J_K(\tau)=\int_0^\infty\frac{f_K(t)}{1-t/\tau}\,dt.
\]
Fix $t_0>0$. At every real $\tau\ge t_0$ the function $J_K$ has
boundary values from above and from below,
\[
 J_K(\tau\pm i0)
 =\operatorname{PV}\!\int_0^\infty\frac{f_K(t)}{1-t/\tau}\,dt
  \mp i\pi\tau f_K(\tau),
\]
and both the principal value and $\tau f_K(\tau)$ are bounded uniformly
for $\tau\ge t_0$. The same uniform bound holds for $J_K(a+ib)$ when
$a\ge t_0/2$ and $0<|b|\le a/4$. Consequently, if
$\sum_m|\beta_m|<\infty$ and the points $\tau_m$ tend inside this
region to real limits $\ge t_0$, all from above or all from below, then
$\sum_m\beta_mJ_K(\tau_m)$ tends to the sum of the corresponding
boundary values.
\end{lemma}

\begin{proof}
The kernel $w_p$ and its derivative decay faster than any power of $t$.
Indeed,
\[
 \frac{d}{dt}\log w_p(t)
 =-\sum_{j\ge0}\frac{p^j}{1+tp^j}
\]
is $O((1+\log t)/t)$ for large $t$, by splitting the sum where
$tp^j=1$, so that
\[
 |f_K'(t)|
 \le C\,t^{K-1}(1+\log t)\,w_p(t)\qquad(t\ge1).
\]
For real $\tau\ge t_0$ split the integral at $\tau/2$ and $3\tau/2$. On
$[0,\tau/2]$ we have $|1-t/\tau|^{-1}\le2$, so this piece is at most
$2\int_0^\infty f_K(t)\,dt$. On the middle interval write
$f_K(t)=f_K(\tau)+(f_K(t)-f_K(\tau))$. The principal value of the
constant term vanishes by symmetry about $\tau$, and the remaining
integrand is bounded by $\tau\sup|f_K'|$, so this part is at most
$\tau^2\sup_{\tau/2\le t\le3\tau/2}|f_K'(t)|$, which is bounded for
$\tau\ge t_0$ by the decay of $f_K'$. On $[3\tau/2,\infty)$ we have
$|1-t/\tau|^{-1}\le3\tau/t$, so this piece is at most
$3\tau\int_{3\tau/2}^\infty f_K(t)\,dt/t$, which is also bounded. The
residue term has modulus $\pi\tau f_K(\tau)$, bounded for the same
reason.

For $\tau=a+ib$ with $a\ge t_0/2$ and $0<|b|\le a/4$ split instead at
$a/2$ and $3a/2$ and subtract $f_K(a)$. On the middle interval
$|t-a|/|t-a-ib|\le1$ and $|\tau|\le2a$, so the difference contributes at
most $2a^2\sup|f_K'|$. The constant term contributes
$-\tau f_K(a)\int_{a/2}^{3a/2}(t-a-ib)^{-1}\,dt$, and the last integral
is a difference of logarithms that is bounded uniformly for either
sign of $b$ and converges as $b\to0$. The outer intervals are bounded
as before, with the constants $4$ and $6a/t$ in place of $2$ and
$3\tau/t$. This proves the uniform bound and the existence of the
limits as $b\to0$ from either side. The statement about sums follows by
dominated convergence.
\end{proof}
\subsection{Polynomial coefficient factors}
\label{app:jets}

Here we extend Lemma~\ref{lem:lateral-tail} to analytic tails with
polynomial factors $P(k)$. In this subsection \(0<p<1\) is arbitrary,
\(h_p=-\log p\), and \(w_p(t)=(-t;p)_\infty^{-1}\).

\begin{lemma}\label{lem:kernel-derivatives}
Fix nonnegative integers \(K,j\) and put \(F_K(t)=t^K w_p(t)\). There is
a constant $C$, depending on \(p,K,j\), such that for \(t\ge1\)
\begin{align}
 |(t\partial_t)^jF_K(t)|
 &\le C t^K(1+\log t)^j w_p(t),\label{eq:kernel-jet-bound}\\
 |\partial_t(t\partial_t)^jF_K(t)|
 &\le C t^{K-1}(1+\log t)^{j+1} w_p(t).
\end{align}
Each $(1+t\partial_t)^jF_K$ is analytic near $t=0$ and has a bounded
first derivative on $[0,1]$. Every function on the right of
\eqref{eq:kernel-jet-bound} is integrable on $[1,\infty)$, and remains
bounded there after multiplication by any fixed power of \(t\).
\end{lemma}
\begin{proof}
Logarithmic differentiation gives
\[
 (t\partial_t)\log w_p(t)=-\sum_{a\ge0}\frac{tp^a}{1+tp^a}.
\]
For \(t\ge1\) split the sum at \(a=\lfloor\log t/h_p\rfloor\). There are
\(O(1+\log t)\) terms before the split, each at most one, and the
terms after it are bounded by a geometric series. Each further
application of $t\partial_t$ to $u/(1+u)$, $u=tp^a$, gives a function
that is bounded for \(u\ge0\) and is \(O(u)\) as \(u\to0\), so the same
splitting shows that \((t\partial_t)^i\log w_p(t)=O(1+\log t)\) for every
fixed $i\ge1$. Writing $(t\partial_t)^jF_K$ as $F_K$ times a polynomial
in these logarithmic derivatives and in $K$ (by the formula for
repeated derivatives of an exponential) gives both inequalities. Near
$t=0$ the product $(-t;p)_\infty$ is holomorphic and nonzero, which
gives the regularity there. The last assertion follows from
\eqref{eq:kernel-growth}.
\end{proof}

\begin{proposition}\label{prop:lateral-jets}
Suppose a Borel transform is the sum of a polynomial and a tail
$C v^K\sum_{m\ge0}e_m/(1-\kappa Q^m v)$, where the polynomial
coefficients, $C$ and the $e_m$ depend holomorphically on the color,
$\sum_m|e_m|$ is locally uniformly bounded in the color, and the
remaining parameters are fixed. Suppose its poles lie on the negative
axis of the Borel plane and it is integrated against $w_p(t)\,dt/h_p$. This
includes the tails $c_k=\lambda^k f(Q^k)$ with $\lambda>0$ of
Theorem~\ref{thm:analytic-symbol}. Then, for colors in compact sets and
moduli in compact subsets of $(0,\infty)$, every fixed Euler derivative
of the sum with respect to $z$ has upper and lower boundary values on
the negative axis, locally uniformly. If \(\theta_z=z\partial_z\) and
\(P\) is a polynomial, then
\begin{equation}\label{eq:jet-lateral-identity}
 S_{P(k)c_k}=P(\theta_z)S_c,\qquad
 \Delta_{P(k)c_k}(X,r)=P(r\partial_r)\Delta_c(X,r).
\end{equation}
The derivatives in the second identity are taken along the negative
axis. For poles of higher order the boundary integrals are Hadamard
finite parts together with their residue terms.
\end{proposition}
\begin{proof}
Off the negative axis, differentiating the coefficients and using
normal convergence gives \(S_{k^jc_k}=\theta_z^j S_c\). What has to be
shown is that the boundary values and their derivatives exist,
uniformly along the whole progression of poles. After the finitely many
initial terms are removed, the integral is a sum, with coefficients
that are absolutely summable, of terms of the form
\[
 J_\eta(\tau;F)=\int_0^\infty
             \frac{F(t)}{1-e^{i\eta}t/\tau}\,dt,
 \qquad \tau\ge\tau_0>0,
\]
with $F(t)=t^Kw_p(t)$, multiplied by a power of the modulus that is
bounded on compact sets. The substitution \(t=\tau u\) shows that
\(\tau\partial_\tau J_\eta(\tau;F)=J_\eta(\tau;(1+t\partial_t)F)\), the
denominator being independent of \(\tau\) after the substitution; for
\(\eta\ne0\) this follows by differentiating under the integral sign.
A term $z^KJ_\eta(\tau;F)$ with $\tau$ proportional to $z^{-1}$ is
mapped by $\theta_z$ to $z^K(K-\tau\partial_\tau)J_\eta(\tau;F)$. Thus
every fixed derivative is a finite linear combination of integrals of
the same Cauchy type with the weights of
Lemma~\ref{lem:kernel-derivatives}.

For each such weight split the integral at \(\tau/2\) and \(3\tau/2\).
The outer pieces are bounded uniformly by the integrable envelopes of
Lemma~\ref{lem:kernel-derivatives}. On the middle interval subtract the
value of the weight at \(\tau\). The integral of the constant is an
explicit logarithm with finite limits as $\eta\to0$ from either side,
and the remaining part is at most
\(C\tau^2\sup_{[\tau/2,3\tau/2]}|F'|\), uniformly in $\eta$.
Lemma~\ref{lem:kernel-derivatives} makes these bounds uniform for
\(\tau\ge\tau_0\), and the residue term \(\pi\tau F(\tau)\) is bounded
as well. When the pole is slightly off the real axis, the same argument
works with subtraction at its real part and the inequality
\(|t-a|/|t-a-ib|\le1\), as in Lemma~\ref{lem:lateral-tail}.

Hence each derivative has a summable majorant along the progression of
poles, and the boundary values may be taken term by term, as may the
finitely many derivatives in the modulus. The boundary values of the
Cauchy kernel give the jump at a simple pole, and differentiating them
gives both the finite parts at higher poles and
\eqref{eq:jet-lateral-identity}. The jump is taken upper minus lower,
as in the main text. Finite sums of such tails with positive
$\lambda$'s and arbitrary initial coefficients are covered by
linearity.
\end{proof}

If $\lambda$ is not real, a rotation of the auxiliary variable moves
the corresponding singular ray to the negative axis. The integration
variable $t$ stays positive and the kernel $w_p(t)$ does not change, so
the same argument gives the boundary values on that ray.

\subsection{Meromorphic transforms with subcritical growth}

\begin{lemma}\label{merged:subcritical-laterals}
Let $0<p,\rho<1$, let $H$ be holomorphic near zero and meromorphic
in the plane, and suppose that its poles belong to finitely many
progressions $\alpha_\nu\rho^{-\mathbb N_0}$, with uniformly bounded
orders. Choose a ray $\mathcal R=e^{i\theta}\mathbb R_{>0}$ containing
some of these poles and a wedge about $\mathcal R$ containing no
other pole directions. Suppose that, outside disjoint disks
$|\xi-\alpha|<\epsilon|\alpha|$ about the possible poles in this
wedge,
\begin{equation}\label{merged:eq:subcritical-growth}
 |H(\xi)|\le C(1+|\xi|)^a
       \exp\{\kappa\log^2(2+|\xi|)\},\qquad
 \kappa<\frac1{2h_p}.
\end{equation}
Then $\mathcal L_pH$ has lateral limits on $\mathcal R$, locally
uniformly for $z=e^{i\theta}r$, $r>0$. With the limit from decreasing
argument minus the limit from increasing argument, their difference is
\begin{equation}\label{merged:eq:subcritical-residues}
 \Delta_H(z)=\frac{2\pi i}{h_pz}
 \sum_{\alpha\in\mathcal R}
   \operatorname{Res}_{\xi=\alpha}
      \bigl(H(\xi)w_p(\xi/z)\bigr).
\end{equation}
The series converges normally for $r$ in compact subsets of
$(0,\infty)$. If $H$ and the bound in
\eqref{merged:eq:subcritical-growth} depend holomorphically and
locally uniformly on an auxiliary parameter, with a fixed possible
pole set, the lateral limits and the residue sum have the same
parameter dependence.
\end{lemma}

\begin{proof}
After coincident points are identified, a finite union of
$\rho$-progressions has positive relative separation. Indeed, ratios
of two points from the same pair of progressions have the form
$c\rho^j$, and only finitely many of these ratios lie in any fixed
annulus about one. Thus $\epsilon$ may be decreased so that the
closed relative disks are disjoint and remain inside the wedge.
Choose two paths $\Gamma_-$ and $\Gamma_+$ running from zero to
infinity on $\mathcal R$ and passing around these disks on the sides
of decreasing and increasing argument, respectively.
If $M$ bounds the pole orders, the maximum principle applied to
$(\xi-\alpha)^M H(\xi)$ on a relative disk bounds $H$ inside it
by the boundary bound times
$(\epsilon|\alpha|/|\xi-\alpha|)^M$.
Along a ray at a fixed nonzero angle from $\mathcal R$, this last
factor is bounded uniformly in $\alpha$. Thus
\eqref{merged:eq:subcritical-growth}, with a changed constant and
polynomial factor, also holds along every such ray in the wedge.

On each smaller sector $|\arg t|\le\theta_0<\pi$, the kernel obeys
\begin{equation}\label{merged:eq:complex-kernel-decay}
 |w_p(t)|\le C_{\theta_0}
 \exp\left\{-\frac{\log^2|t|}{2h_p}
                         +C_{\theta_0}\log|t|\right\},
 \qquad |t|\ge1.
\end{equation}
To verify this estimate, put $J=\lfloor\log|t|/h_p\rfloor$.
For $j\le J$ the inequality
$|1+tp^j|\ge c_{\theta_0}(1+|t|p^j)$ gives
\[
 \sum_{j=0}^J\log|1+tp^j|
 \ge\sum_{j=0}^J(\log|t|-jh_p)-C_{\theta_0}(J+1).
\]
For $j>J$, the product is bounded away from zero uniformly in
$\arg t$ and $|t|p^{J+1}\in[p,1]$: its first finitely many factors
stay away from zero, and the logarithms of the remaining factors
are uniformly summable. The resulting lower bound for the full
logarithm is $\log^2|t|/(2h_p)-O(\log|t|)$, which proves
\eqref{merged:eq:complex-kernel-decay}.

For $z$ in a compact part of the ray and a sufficiently small
neighborhood on either side, $\xi/z$ stays in a fixed smaller
sector about the positive axis along both paths. Combining
\eqref{merged:eq:subcritical-growth} and
\eqref{merged:eq:complex-kernel-decay}, the integrand in
\[
 \frac1{h_pz}\int_{\Gamma_\pm}H(\xi)w_p(\xi/z)\,d\xi
\]
is bounded at infinity by
$C\exp\{-c\log^2|\xi|+C'\log|\xi|\}$ for some $c>0$.
The length of a detour is $O(|\alpha|)$, and the possible pole
moduli lie on finitely many geometric progressions. Hence these
integrals converge normally. They are also integrable at zero,
where $H$ and $w_p$ are holomorphic.

For $\arg z$ on the appropriate side of $\theta$, deform the ray
$z\mathbb R_{>0}$ to $\Gamma_\pm$ in the corresponding punctured
wedge. The poles of $w_p(\xi/z)$ lie on
$-z\mathbb R_{>0}$ and do not meet this deformation. The estimates
above make the connecting integrals at infinity tend to zero;
one may choose outer radii with a fixed relative separation from
the pole set. The connecting integral at zero tends to zero by
holomorphy. Thus the path integral equals
$\mathcal L_pH(z)$ on its side. Normal convergence lets $z$ tend to
the ray and proves the lateral limits.

Subtract the two path integrals. Each pair of detours encloses its
pole with positive orientation, and the residue theorem gives
\eqref{merged:eq:subcritical-residues}. Its normal convergence
follows either from the same contour bounds or from Cauchy's
estimate on each relative disk: the residue is bounded by its
boundary length times the supremum of
$|H(\xi)w_p(\xi/z)|$, which has the summable bound just proved.
All bounds are locally uniform in the auxiliary parameters under
the stated assumptions. Differentiation of the path integrals, or
Morera's theorem with normal convergence, gives the final assertion.
\end{proof}

\section{Peripheral character and rank-two torsion}\label{app:torsion}

We compare the leading rank-two factor with meridian-normalized
adjoint torsion, using Yoon's formula
\cite[Theorem~2, Remark~1, and Lemma~2]{YoonTwist}.

For $m=1,n=2$, let $x=e^u$, $X=x+x^{-1}$ and $Y=QX$.
The closed contour integrand, including its product prefactor, has
the local classical phase
\begin{align*}
 W_u(z)={}&2\operatorname{Li}_2(1)
 -\operatorname{Li}_2(x)-\operatorname{Li}_2(x^{-1})
 -\operatorname{Li}_2(z)-\operatorname{Li}_2(z^{-1})\\
 &+\operatorname{Li}_2(x/z)+\operatorname{Li}_2(x^{-1}/z),
\end{align*}
with branches fixed locally away from the product divisors.
Its derivative is exactly
\[
 \frac{\partial W_u}{\partial\log z}
 =\log\frac{(1-z)(1-x/z)(1-x^{-1}/z)}{1-z^{-1}}
 =\log(2\cosh u-z-z^{-1}).
\]
The same expression follows directly from the finite Newton factors.
In rank two,
\[
 B_{k+1}^{(2)}(X;h)
 =\bigl(2\cosh u-2\cosh((k+1)h)\bigr)B_k^{(2)}(X;h).
\]
With $t=kh$, the limiting logarithm of the displayed factor is
$\log(2\cosh u-e^t-e^{-t})$.  Thus the period uses the coordinate
$z=e^t$.  Its stationary equation is
$z+z^{-1}=2\cosh u-1$.

To state the torsion normalization independently, put $\mu=e^{u/2}$.
Yoon's negative twist index $-1$ is the figure-eight case; the
Chebyshev values $S_{-1}=-1$, $S_{-2}=-T$ give its defining equation
and meridian-normalized, sign-refined homological torsion as
\[
 F(\mu,T)=(T-1)(T-\mu^2-\mu^{-2})+1=0,\qquad
 \tau_\mu=\tfrac12\partial_TF.
\]
The substitution $T=1+z$ identifies the stationary equation with
$F=0$.  Meridian regularity at the complete characters can be checked
in the trace coordinate $\xi=\mu+\mu^{-1}$.  In these coordinates,
\[
 \widehat F(\xi,T)=(T-1)(T-\xi^2+2)+1.
\]
At $\xi=2$ and $T=(3\pm i\sqrt3)/2$, its derivative with respect to
$T$ is $\pm i\sqrt3\ne0$. Hence $\xi$ is a local coordinate on
the character curve, so the differential of the meridian trace is nonzero.
Yoon's regularity criterion \cite[Remark~2]{YoonTwist} therefore
applies at these characters and in a neighborhood of each.
On this curve, $\mu^2+\mu^{-2}=T+(T-1)^{-1}$, and therefore
\[
 \tau_\mu=\frac12\bigl((T-1)-(T-1)^{-1}\bigr)
          =\frac12(z-z^{-1}).
\]
On the upper-sector dominant saddle $z=e^{-i\pi/3}$ this is
$\tau_\mu=-i\sqrt3/2$.  Since the contour is parameterized by
$t=i\theta$, its Hessian at the saddle is
\[
 W''(\theta_*)=2\tau_\mu,\qquad \theta_*=\frac{5\pi}3.
\]
The rank-two leading factor in \eqref{r4:upper-period} can therefore
be written
\[
 \sqrt{-\frac{h}{4\pi\tau_\mu}}
   =\frac{\sqrt h}{\sqrt{2\pi}\,3^{1/4}}\,e^{-i\pi/4},
\]
where the square root is determined by the positively oriented circle.
In rank $n$ the additional factors in \eqref{r4:rank-period} multiply
this $\mathrm{SL}_2$ factor by
\[
 \frac{(-1)^d h^{-2d}}{d!\,\Gamma(n/2)^2}.
\]

\section{The modular form of the node divisor}\label{app:modular-divisor}

By \eqref{pa:monodromy-divisor}, the node divisor $E_n$ divides
$\mathcal B_{m,n}$. In this appendix we express $E_n$ through the dual
nome. Let
\[
 E_n(h,u)=(Q^{n/2}e^u,Q^{n/2}e^{-u};Q)_\infty,
 \qquad Q=e^{-h},\quad\widetilde Q=e^{-4\pi^2/h},
\]
where $X=e^u+e^{-u}$ and $Y=Q^{n/2}X$.

\begin{proposition}\label{r4:node-modular}
For $\Re h>0$ the following identities hold between meromorphic
functions of $u$:
\begin{align}
 E_{2r}(h,u)
 &=\frac{e^{(u^2/2-\pi^2/3)/h+h/12}}
        {(Qe^u,Qe^{-u};Q)_{r-1}}
   \frac{\sin(\pi u/h)}{\sinh(u/2)}
   (\widetilde Qe^{2\pi iu/h},
        \widetilde Qe^{-2\pi iu/h};\widetilde Q)_\infty,
 \label{r4:even-divisor}\\
 E_{2r+1}(h,u)
 &=\frac{2e^{(u^2/2-\pi^2/3)/h-h/24}}
        {(Q^{1/2}e^u,Q^{1/2}e^{-u};Q)_r}
   \cos(\pi u/h)
   (-\widetilde Qe^{2\pi iu/h},
        -\widetilde Qe^{-2\pi iu/h};\widetilde Q)_\infty.
 \label{r4:odd-divisor}
\end{align}
Here $r\geq1$; the odd identity also holds for $r=0$, and this case is
the starting point of the proof.

If, on a fixed closed cone $|\arg h|\leq\theta_0<\pi/2$, the pair
$(h,u)$ satisfies
\begin{equation}
 4\pi^2\Re(1/h)-2\pi|\Im(u/h)|\geq\eta/|h|
 \label{r4:dual-gap}
\end{equation}
for a fixed $\eta>0$, then, as $h\to0$, either dual product in
\eqref{r4:even-divisor}--\eqref{r4:odd-divisor} equals
$1+O(e^{-\eta/|h|})$, uniformly.
\end{proposition}

\begin{proof}
Let $\theta_Q(z)=(Q,z,Q/z;Q)_\infty$.  Poisson summation of its
bilateral Gaussian series, followed by the Jacobi triple product,
gives
\begin{align*}
 \theta_Q(e^u)={}&-2\sqrt{\frac{2\pi}{h}}
  e^{u^2/(2h)+u/2+h/8-\pi^2/(2h)}\sin(\pi u/h)\\
 &\hspace{6mm}\times(\widetilde Q;\widetilde Q)_\infty
  (\widetilde Qe^{2\pi iu/h},
       \widetilde Qe^{-2\pi iu/h};\widetilde Q)_\infty.
\end{align*}
The square root is the one positive for $h>0$ and holomorphic in
$\Re h>0$.  Also
\[
 (Q;Q)_\infty
   =\sqrt{\frac{2\pi}{h}}
      e^{-\pi^2/(6h)+h/24}
          (\widetilde Q;\widetilde Q)_\infty.
\]
These are the usual theta and eta modular transformations;
see \cite[\S20.7(viii), \S23.18]{DLMFModular}.
Dividing the first identity by
$(Q;Q)_\infty(1-e^u)$ proves \eqref{r4:even-divisor} when $r=1$.
Replacing $u$ by $u-h/2$ before dividing by $(Q;Q)_\infty$
proves the rank-one base of \eqref{r4:odd-divisor}.  Removing the
finitely many initial product factors proves the general identities.

Under \eqref{r4:dual-gap} the first two dual arguments have modulus
at most $e^{-\eta/|h|}$, while
$|\widetilde Q|\leq e^{-4\pi^2\cos\theta_0/|h|}$.
The absolute sum of their product factors is therefore
$O(e^{-\eta/|h|})$, giving the asserted uniform estimate.
\end{proof}

For example, $|u|\leq u_0<2\pi\cos\theta_0$ implies
\eqref{r4:dual-gap} with
$\eta=4\pi^2\cos\theta_0-2\pi u_0>0$.
These formulas apply to the divisor for every $m$, because
$\mathcal B_{m,n}(Q,Y)=E_n(h,u)\mathcal C_{m,n}(Q,Y)$.

\section*{Declaration of AI Use}
Claude Opus 5.5 was used solely for language editing and polishing of this manuscript. GPT-6.2 Astra was used to assist in checking the accuracy of the mathematical arguments and computations.


\begin{thebibliography}{99}

\bibitem{CLZDoubleTwist}
Q.~Chen, K.~Liu, and S.~Zhu,
\emph{Cyclotomic expansions for the colored HOMFLY-PT invariants of double twist knots},
preprint, arXiv:2110.03616 (2021),
\url{https://arxiv.org/abs/2110.03616}.

\bibitem{CMPPush}
I.~Corwin, K.~Matveev, and L.~Petrov,
\emph{The $q$-Hahn PushTASEP},
Int. Math. Res. Not. IMRN \textbf{2021} (2021), no.~3, 2210--2249.

\bibitem{DVZ}
L.~Di Vizio and C.~Zhang,
\emph{On $q$-summation and confluence},
Ann. Inst. Fourier (Grenoble) \textbf{59} (2009), no.~1, 347--392,
\href{https://doi.org/10.5802/aif.2433}{doi:10.5802/aif.2433}.

\bibitem{Dreyfus}
T.~Dreyfus,
\emph{Building meromorphic solutions of $q$-difference equations using a Borel--Laplace summation},
Int. Math. Res. Not. IMRN \textbf{2015} (2015), no.~15, 6562--6587,
\href{https://doi.org/10.1093/imrn/rnu137}{doi:10.1093/imrn/rnu137}.

\bibitem{Fang}
H.~Fang,
\emph{Quantum Color Newton Rings},
preprint, arXiv:2610.21102 (2026),
\url{https://arxiv.org/abs/2610.21102}.

\bibitem{FangZhou}
H.~Fang and T.~Zhou,
\emph{Cyclotomic Newton expansions and a rank-uniform integer-valued Newton completion},
preprint, arXiv:2608.29585 (2026),
\url{https://arxiv.org/abs/2608.29585}.

\bibitem{FantiniWheeler}
V.~Fantini and C.~Wheeler,
\emph{Summability for state integrals of hyperbolic knots},
preprint, arXiv:2410.20973v1 (2024),
\url{https://arxiv.org/abs/2410.20973}.

\bibitem{FMNCumulants}
V.~F{\'e}ray, P.-L.~M{\'e}liot, and A.~Nikeghbali,
\emph{Mod-$\phi$ convergence: Normality zones and precise deviations},
SpringerBriefs in Probability and Mathematical Statistics,
Springer, Cham, 2016.
Theorem numbering follows arXiv:1304.2934v4 (2015), Theorems~9.6--9.7.

\bibitem{GGM}
S.~Garoufalidis, J.~Gu, and M.~Mari\~no,
\emph{The resurgent structure of quantum knot invariants},
Comm. Math. Phys. \textbf{386} (2021), 469--493,
\href{https://doi.org/10.1007/s00220-021-04076-0}{doi:10.1007/s00220-021-04076-0}.

\bibitem{GGMW}
S.~Garoufalidis, J.~Gu, M.~Mari\~no, and C.~Wheeler,
\emph{Resurgence of Chern--Simons theory at the trivial flat connection},
Comm. Math. Phys. \textbf{406} (2025), article~20,
\href{https://doi.org/10.1007/s00220-024-05149-6}{doi:10.1007/s00220-024-05149-6}.

\bibitem{GKState}
S.~Garoufalidis and R.~Kashaev,
\emph{From state integrals to $q$-series},
Math. Res. Lett. \textbf{24} (2017), no.~3, 781--801.

\bibitem{GLL}
S.~Garoufalidis, A.~D.~Lauda, and T.~T.~Q.~L\^e,
\emph{The colored HOMFLYPT function is $q$-holonomic},
Duke Math. J. \textbf{167} (2018), 397--447,
\url{https://arxiv.org/abs/1604.08502}.

\bibitem{GLAnalytic}
S.~Garoufalidis and T.~T.~Q.~L\^e,
\emph{An analytic version of the Melvin--Morton--Rozansky conjecture},
preprint, arXiv:math/0503641v2 (2005),
\url{https://arxiv.org/abs/math/0503641}.

\bibitem{GLAsymptotics}
S.~Garoufalidis and T.~T.~Q.~L\^e,
\emph{Asymptotics of the colored Jones function of a knot},
Geom. Topol. \textbf{15} (2011), no.~4, 2135--2180.

\bibitem{GLSurvey}
S.~Garoufalidis and T.~T.~Q.~L\^e,
\emph{A survey of $q$-holonomic functions},
Enseign. Math. \textbf{62} (2016), nos.~3--4, 501--525,
\href{https://doi.org/10.4171/LEM/62-3/4-7}{doi:10.4171/LEM/62-3/4-7};
\url{https://arxiv.org/abs/1601.07487}.

\bibitem{GasperRahman}
G.~Gasper and M.~Rahman,
\emph{Basic hypergeometric series}, 2nd ed.,
Encyclopedia of Mathematics and its Applications, vol.~96,
Cambridge University Press, Cambridge, 2004.

\bibitem{HM}
K.~Hikami and H.~Murakami,
\emph{Colored Jones polynomials with polynomial growth},
Commun. Contemp. Math. \textbf{10} (2008), suppl.~1, 815--834,
\url{https://arxiv.org/abs/0711.2836}.

\bibitem{IS}
M.~E.~H.~Ismail and D.~Stanton,
\emph{$q$-Taylor theorems, polynomial expansions, and interpolation of entire functions},
J. Approx. Theory \textbf{123} (2003), no.~1, 125--146,
\href{https://doi.org/10.1016/S0021-9045(03)00076-5}{doi:10.1016/S0021-9045(03)00076-5}.

\bibitem{IMMM}
H.~Itoyama, A.~Mironov, A.~Morozov, and An.~Morozov,
\emph{HOMFLY and superpolynomials for figure eight knot in all symmetric and antisymmetric representations},
J. High Energy Phys. \textbf{2012} (2012), no.~7, article~131,
\href{https://doi.org/10.1007/JHEP07(2012)131}{doi:10.1007/JHEP07(2012)131}.

\bibitem{KMM}
R.~Kashaev, V.~Mangazeev, and M.~Mari\~no,
\emph{Analytic quantum invariants},
preprint, arXiv:2610.24247 (2026),
\url{https://arxiv.org/abs/2610.24247}.

\bibitem{Kawagoe}
K.~Kawagoe,
\emph{The colored HOMFLY-PT polynomials of the trefoil knot, the figure-eight knot, and twist knots},
J. Geom. Phys. \textbf{213} (2025), article~105488.
Formula numbering follows arXiv:2107.08678v1.

\bibitem{Masbaum}
G.~Masbaum,
\emph{Skein-theoretical derivation of some formulas of Habiro},
Algebr. Geom. Topol. \textbf{3} (2003), no.~1, 537--556.

\bibitem{MorozovTselousov}
A.~Morozov and N.~Tselousov,
\emph{Differential expansion for antiparallel triple pretzels: the way the factorization is deformed},
Eur. Phys. J. C \textbf{82} (2022), article~912,
\href{https://doi.org/10.1140/epjc/s10052-022-10851-7}{doi:10.1140/epjc/s10052-022-10851-7}.

\bibitem{MurakamiTran}
H.~Murakami and A.~T.~Tran,
\emph{On the asymptotic behavior of the colored Jones polynomial of the figure-eight knot associated with a real number},
in \emph{Low Dimensional Topology and Number Theory}
(M.~Morishita, H.~Nakamura, and J.~Ueki, eds.),
Springer Proceedings in Mathematics \& Statistics, vol.~456,
Springer, Singapore, 2025, pp.~175--209.
Formula numbering follows arXiv:2109.04664v4.

\bibitem{DLMFHeine}
NIST Digital Library of Mathematical Functions,
\emph{Chapter 17: $q$-hypergeometric and related functions},
Eq.~(17.6.7) (Heine's second transformation),
\url{https://dlmf.nist.gov/17.6.E7}, accessed November~8, 2026.

\bibitem{DLMFModular}
NIST Digital Library of Mathematical Functions,
\emph{Theta and eta modular transformations},
Sections~20.7(viii) and~23.18,
\url{https://dlmf.nist.gov/20.7.viii} and
\url{https://dlmf.nist.gov/23.18}, accessed November~8, 2026.

\bibitem{RamisGrowth}
J.-P.~Ramis,
\emph{About the growth of entire functions solutions of linear algebraic $q$-difference equations},
Ann. Fac. Sci. Toulouse Math. (6) \textbf{1} (1992), no.~1, 53--94,
\href{https://doi.org/10.5802/afst.739}{doi:10.5802/afst.739}.

\bibitem{RSZ}
J.-P.~Ramis, J.~Sauloy, and C.~Zhang,
\emph{Local analytic classification of $q$-difference equations},
Ast\'erisque \textbf{355} (2013), vi+151~pp.,
\url{https://arxiv.org/abs/0903.0853}.

\bibitem{RT}
N.~Yu.~Reshetikhin and V.~G.~Turaev,
\emph{Ribbon graphs and their invariants derived from quantum groups},
Comm. Math. Phys. \textbf{127} (1990), 1--26,
\href{https://doi.org/10.1007/BF02096491}{doi:10.1007/BF02096491}.

\bibitem{Witten}
E.~Witten,
\emph{Analytic continuation of Chern--Simons theory},
in \emph{Chern--Simons Gauge Theory: 20 Years After},
AMS/IP Studies in Advanced Mathematics, vol.~50,
American Mathematical Society and International Press, 2011, pp.~347--446.

\bibitem{WongAu}
K.~H.~Wong and T.~K.-K.~Au,
\emph{Asymptotic behavior of colored HOMFLY polynomial of figure eight knot},
preprint, arXiv:1711.04437v1 (2017),
\url{https://arxiv.org/abs/1711.04437}.

\bibitem{YoonTwist}
S.~Yoon,
\emph{A vanishing identity of adjoint Reidemeister torsions of twist knots},
Algebr. Geom. Topol. \textbf{22} (2022), no.~1, 227--249.
Formula numbering follows arXiv:2002.12576v2.

\end{thebibliography}
\end{document}